\documentclass[11pt,a4paper]{article}
\usepackage[breaklinks,colorlinks,
  citecolor=blue,
  urlcolor=blue]{hyperref}

\usepackage{cite}

\usepackage[left=2.2cm,right=2.2cm,top=2cm,bottom=2cm]{geometry}
\usepackage{fixmath}
\usepackage{amsmath}
\usepackage{empheq}
\usepackage{bbm}
\usepackage{amssymb}
\usepackage{mathrsfs}  
\usepackage{amsthm}
\usepackage{bm}
\usepackage{circuitikz}
\usepackage{subcaption}

\usepackage{nomencl}
\makenomenclature

\usepackage{array}

\usepackage{tikz-cd}

\usepackage{url}

\theoremstyle{plain}
\newtheorem{theorem}{Theorem}[section]
\newtheorem{proposition}[theorem]{Proposition}
\newtheorem{lemma}[theorem]{Lemma}

\newtheorem{definition}[theorem]{Definition}
\newtheorem{example}[theorem]{Example}

\newtheorem{remark}[theorem]{Remark}

\newcommand{\giantzero}{\mbox{\normalfont\Huge0}}
\newcommand{\bigzero}{\mbox{\normalfont\LARGE 0}}

\newcommand{\rvline}{\hspace*{-\arraycolsep}\vline\hspace*{-\arraycolsep}}

\newcommand{\rowsvdots}{\multicolumn{1}{@{}c@{}}{\vdots}}

\begin{document}

\title{{\bf Learnability of Linear Port-Hamiltonian Systems}}

\author{Juan-Pablo~Ortega$^{1}$ and Daiying Yin$^{1}$}
%
\maketitle

\begin{abstract}
A complete structure-preserving learning scheme for single-input/single-output (SISO) linear port-Hamiltonian systems is proposed. The construction is based on the solution, when possible, of the unique identification problem for these systems, in ways that reveal fundamental relationships between classical notions in control theory and crucial properties in the machine learning context, like structure-preservation and expressive power. In the canonical case, it is shown that the set of uniquely identified systems  can be explicitly characterized as a smooth manifold endowed with global Euclidean coordinates, which allows concluding that the parameter complexity necessary for the replication of the dynamics is only $\mathcal{O}(n)$ and not $\mathcal{O}(n^2)$, as suggested by the standard parametrization of these systems. Furthermore, it is shown that linear port-Hamiltonian systems can be learned while remaining agnostic about the dimension of the underlying data-generating system. Numerical experiments show that this methodology can be used to efficiently estimate linear port-Hamiltonian systems out of input-output realizations, making the contributions in this paper the first example of a structure-preserving machine learning paradigm for linear port-Hamiltonian systems based on explicit representations of this model category.\\
	
\end{abstract}

\makeatletter
\addtocounter{footnote}{1} \footnotetext{%
Juan-Pablo~Ortega and Daiying Yin are with the Division of Mathematical Sciences, Nanyang Technological University, Singapore. Their email addresses are {\texttt{Juan-Pablo.Ortega@ntu.edu.sg}} and {\texttt{YIND0004@e.ntu.edu.sg}}.}
\makeatother

%

\tableofcontents

\addcontentsline{toc}{section}{Glossary of main symbols}
\makenomenclature

\nomenclature{${\displaystyle \mathbb{J}_n=
\begin{bmatrix}
	0 & \mathbb{I}_n\\
	-\mathbb{I}_n & 0
\end{bmatrix}}$ }{Canonical symplectic matrix}
\nomenclature{$F:\mathcal{Z}\times\mathcal{U}\rightarrow\mathcal{Z}$}{State equation}
\nomenclature{$H:\mathbb{R}^{2n} \longrightarrow \mathbb{R} $}{Hamiltonian function}
\nomenclature{$PH_n$}{The space of $2n$-dimensional linear normal form port-Hamiltonian systems \eqref{ph_definition} }
\nomenclature{$PH_n^{can}$}{The subspace of $PH_n$ consisting of canonical linear normal form port-Hamiltonian systems}
\nomenclature{$\Theta_{PH_n}$}{The space of parameters $(Q,B)$ for $PH_n$}
\nomenclature{$\theta_{PH_n}$}{The map that sends parameters in  for $\Theta_{PH_n}$ to the corresponding state space system in $PH_n$}
\nomenclature{$\mathcal{PH}_n$}{The space of input-output dynamics/filters induced by systems in $PH_n$}
\nomenclature{$\mathcal{PH}_n^{can}$}{The space of input-output dynamics/filters induced by systems in $PH_n^{can}$}
\nomenclature{$Sp(2n,\mathbb{R})$}{Symplectic group}
\nomenclature{$\mathfrak{sp} (2n,\mathbb{R})$}{Lie algebra of the symplectic group}
\nomenclature{$Q$}{Quadratic form that determines a linear Hamiltonian system}
\nomenclature{$B$}{Input matrix of a port-Hamiltonian system in normal form}
\nomenclature{$CH_n/OH_n$}{The space of $2n$-dimensional controllable/observable Hamiltonian representations}
\nomenclature{$\mathcal{CH}_n/\mathcal{OH}_n/\mathcal{CH}^0_n$}{The space of filters induced by $CH_n$/$OH_n$/$CH_n$ with zero initial condition}
\nomenclature{$\Theta_{CH_n}/\Theta_{OH_n}$}{The space of parameters $({\bf d},{\bf v})$ for $CH_n$ and/or $OH_n$, which are the same}
\nomenclature{$\theta_{CH_n}/\theta_{OH_n}$}{The map that send parameters in  for $\Theta_{CH_n}/\Theta_{OH_n}$ to the corresponding state space system in $CH_n/OH_n$}
\nomenclature{$\sim_{filter}$}{The equivalence relation of inducing the same filter}
\nomenclature{$\sim_{sys}$}{The equivalence relation of system automorphism}
\nomenclature{$\sim_{\star}$}{An equivalence relation defined on $\Theta_{CH_n}$}
\nomenclature{$\mathcal{G}_n\rightrightarrows PH_n$}{Port-Hamiltonian groupoid, see Proposition \ref{groupoid_orbit}}
\nomenclature{$\mathcal{H}_n\rightrightarrows\Theta_{CH_n}$}{Reduced port-Hamiltonian groupoid, see Proposition \ref{diagonalized port-Hamiltonian groupoid}}
\nomenclature{$\Theta_{CH_{n}}^{can}$}{The subset of $\Theta_{CH_{n}}$ that corresponds to canonical systems}
\nomenclature{$S_n$}{Permutation group of $n$-elements}
\nomenclature{$\mathbb{T}^n$}{The $n$-torus}
\nomenclature{$\mathcal{X}^n_{\uparrow}$}{The set of $n$-tuples of distinct real numbers in increasing order}
\nomenclature{$PH_{m,n}$}{The subspace of $PH_{m}$ containing all $(Q^{\prime},B^{\prime})=\bigg(O\begin{bmatrix}
		Q&0\\
		0&\mathbb{I}_{2m-2n}
	\end{bmatrix}O^{T},O\begin{bmatrix}
		B\\
		0
	\end{bmatrix}\bigg)$, $(Q,B)\in PH_n$, $O\in O(2m,\mathbb{R})$}
\nomenclature[26]{$\Theta_{CH_{m,n}}$}{The space of parameters $({\bf d}^{\prime},{\bf v}^{\prime})$ for $PH_{m,n}$}
\printnomenclature

\section{Introduction}
%
%
%
%
Machine learning has experienced substantial development in recent years due to significant advances in algorithmics and a fast growth in computational power. The universal approximation properties of neural networks \cite{cybenko, hornik} and other similar families make it possible for them to learn any function with very few prior assumptions. A typical modus operandi in supervised machine learning is first to choose a neural network architecture, to perform forward propagation using available data, to compute some loss function, and then to carry out backward propagation, that is, gradient descent, to recursively optimize the parameters. This paradigm has proved to be very successful in the learning of numerous complicated tasks, including time-series forecasting \cite{LSTM}, computer vision \cite{Krizhevsky2012}, and natural language processing\cite{Devlin2018}.\par
 
In physics and engineering, machine learning is called to play an essential role in predicting and integrating the equations associated with physical dynamical systems. Physical systems are primarily formulated in terms of ordinary, time-delay, and partial differential equations that can be deduced mostly from variational principles. Consequently, some researchers propose to learn adequately discretized versions of their corresponding vector fields 
(see, for instance, \cite{Raissi2018},\cite{Qin2019}, \cite{Long2018}, and references therein). In addition to vector fields learning, researchers have proposed ``model-free" methods like {\it transformers} \cite{Shalova2020, acciaio2022metric}, {\it reservoir computing} \cite{Jaeger04, Ott2018, Pathak:PRL, pathak2018}, {\it recurrent neural networks} \cite{Bailer-Jones1998}, {\it  convolutional neural networks} \cite{Mukhopadhyay2020}, or {\it LSTMs} \cite{Wang2017}.  \par 

Various universal approximation properties theoretically explain the empirical success of some of these approaches (see, for instance, \cite{RC6, RC7, RC8, RC20}) of some of these learning paradigms. Nevertheless, for physics-related problems, like in mechanics or optics, it is natural to build into the learning algorithm any prior knowledge that we may have about the system based on physics' first principles. This may include specific forms of the laws of motion, conservation laws, symmetry invariance, as well as other underlying geometric and variational structures. This observation regarding the construction of structure-preserving schemes, has been profusely exploited with much success before the emergence of machine learning in the field of numerical integration \cite{gonzalez2000time, marsden_west_2001, LeimkuhlerBook, mclachlan2006geometric}. Many examples in that context show how the failure to maintain specific conservation laws can lead to physically inconsistent solutions.

The translation of this idea to the context of machine learning has lead to the emergence of a new domain collectively known as {\it physics-informed machine learning} (see \cite{raissi2017physics, wu2018physics, karniadakis2021physics} and references therein). In the specific case of Hamiltonian systems, the two main structural constraints are that the flow is symplectic and the energy, i.e., the Hamiltonian, is conserved along the flow. Additionally, symmetries are frequently present, which carries the emergence of additional conserved quantities in the form of the so-called momentum maps via Noether's Theorem \cite{Abraham1978, Marsden1994, Ortega2004}. These are all examples of qualitative properties to be preserved by the learning algorithms. Needless to say that the above-mentioned ``model-free" approaches generically fail to preserve all these structures.  With these in mind, several attempts in the literature have been made to develop tailor-made learning algorithms for Hamiltonian systems. For example,  in \cite{greydanus2019hamiltonian, celledoni2023learning} neural methods are proposed to learn the Hamiltonian function directly. In \cite{Chen2020}, a symplectic recurrent neural network is proposed that uses symplectic integration while matching the predictions and observations and leads to a structure-preserving paradigm. Other structure-preserving methods include the so-called SympNet \cite{Jin2020}, the generating function neural networks (GFNN) in \cite{Chen2021}, and the symplectic reversible neural networks in \cite{Valperga2022}. SympNet constructs a universal approximating family of symplectic maps, while GFNN applies a modified KAM theory to control long-term prediction error. Symplectic reversible neural networks are also proposed as a family of universal approximating maps that concern, in particular, reversible symplectic dynamics. In \cite{Zhong2020Symplectic}, a parametric framework of learning Hamiltonian state dynamics with control is proposed, with the assumption that the Hamiltonian is separable. Under the same assumption, \cite{Tong2021} proposes to learn with a parametrized Hamiltonian in a Taylor series form.\par 
 
 This paper's focus differs from the references mentioned above in two ways. First, these methods are designed to learn the state evolution of Hamiltonian systems, whereas our approach focuses on {\it learning the input-output dynamics of port-Hamiltonian systems while remaining agnostic about the physical state space}. As will be introduced later on, these systems have an underlying Dirac structure that describes the geometry of numerous physical systems with external inputs \cite{vanderSchaft2014} and includes the dynamics of the observations of Hamiltonian systems as a particular case. Even though various learning schemes for these systems have been already proposed in the literature \cite{nageshrao2015adaptive, cherifi2020overview, desai2021port, beckers2022gaussian} most works on the learning of Hamiltonian systems deal with autonomous (separable) Hamiltonian systems on which one assumes access to the entire phase space and not only to its observations. Second, instead of a general nonlinear system for which only approximation error can be possibly estimated, we consider, as a first approach {\it exclusively linear systems}, in which case, we can obtain explicit representations of linear port-Hamiltonian systems in normal form and characterize the symmetries and quotient spaces associated to the invariance by system automorphisms. Thereby, we propose a structure-preserving learning paradigm with a provable minimal parameter space.
 
 The contributions in this paper are contained in several results that we briefly introduce in the following lines. In Section \ref{preliminary}, we define the notion of linear port-Hamiltonian systems in normal form and present some necessary introductory concepts. We start in Theorem \ref{main_theorem} by introducing system morphisms that allow us to represent any linear port-Hamiltonian system in normal form as the image of another linear system of the same dimension in which the state equation is in controllable canonical form.  An obvious observation is that since the original port-Hamiltonian system and the new linear system are linked by a system morphism, the image of the input/output relations of the latter are input/output relations of the former. In particular, the new system can be used to learn to reproduce the input/output dynamics of the original port-Hamiltonian system (for a subspace of initial conditions) and {\it this learning paradigm is structure-preserving by construction}. Similarly, Theorem \ref{main_theorem} also contains another type of system morphisms that link any linear port-Hamiltonian system in normal form to some linear system of the same dimension in observable canonical form. Consequently, the input-ouput relations of the original port-Hamiltonian system with respect to any initial condition can be captured by the observable Hamiltonian representation. Both representations are based on classical techniques from control theory, the Cayley-Hamilton theorem, and are ultimately corollaries of the Williamson normal form \cite{williamson1936algebraic, williamson1937normal, Ikramov2018}. We show that the controllable and observable representations are closely related to each other, and both system morphisms become isomorphisms for canonical port-Hamiltonian systems. However, for the purpose of learning a general port-Hamiltonian system that may not be canonical, we reveal that there is a trade-off between the structure-preserving property and the expressive power. These results establish a strong link between classical notions in the control theory, that is controllability and observability, and those in machine learning, namely, structure-preservation and expressive power. \par
 
 Based on these explicit constructions and using the parametrizations that come with them, we tackle in Section \ref{quotient} the unique identifiability of input-output dynamics of linear port-Hamiltonian systems in normal form. Such a characterization is obviously needed to solve the model estimation problem since in applications we only have access to input/output data and different state space systems can induce the same filter that produces that data. This fact has important implications when it comes to the learning of port-Hamiltonian systems out of finite-sample realizations of a given data-generating process because such degeneracy makes impossible its exact recovery. Said differently, it is not the space of port-Hamiltonian systems that needs to be characterized but its quotient space with respect to the equivalence relation defined by the constraint on inducing the same input/ouput system. We shall see in Subsection \ref{sub0} that the presence of non-canonical systems in $PH_n $ makes it in general difficult to directly characterize that quotient space and we shall settle for the closest to it that we can get, namely, the quotient space by system automorphisms that, as it will be justified, approximates the general case in a certain sense and admits an explicit characterization as a {Lie groupoid} orbit space (Subsection \ref{sub2}). In Subsection \ref{sub3}, we restrict our identification analysis to {canonical port-Hamiltonian systems} and show, first, that in that situation eliminating the system isomorphisms completely identifies the set of input/output systems, and second, that the corresponding quotient spaces can be characterized as orbit spaces with respect to a  {group} (as opposed to groupoids in the general unrestricted case) action, where the group is explicitly given by a semi-direct product. Moreover, (see Subsection \ref{sub6}) this orbit space can be explicitly endowed with a smooth manifold structure that has global Euclidean coordinates that can be used at the time of constructing estimation algorithms. 
 Consequently, canonical port-Hamiltonian dynamics can be identified fully and explicitly in either the controllable or the observable Hamiltonian representations and learned by estimating a unique set of parameters in a smooth manifold that is obtained as a group orbit space.
 
Another learning-related problem that we tackle is that, in applications, one is obliged to remain agnostic as to the dimension of the underlying data-generating port-Hamiltonian system. This leads to the difficulty of choosing the dimension of the controllable/observable Hamiltonian representations. We solve this issue by proving in Theorem \ref{higher dimension} that, for $m\geq n$, any $2n$-dimensional linear port-Hamiltonian system in normal form can be regarded as the restriction of a $2m$-dimensional one to some subspace. This fact, together with some subsequent results, guarantees theoretically that we can choose a sufficiently large $m$ in practice and parametrize the observable Hamiltonian representation in dimension $2m$ and use it for learning without assuming any knowledge about the dimension of the data generating system. The paper concludes with some numerical examples in Section \ref{numerics} that illustrate the viability of the method that we propose in systems with various levels of complexity and dimensions as well as the computational advantages associated with the use of the parameter space in which unique identification is guaranteed.
 
\section{Preliminaries}\label{preliminary}

In this section, we introduce various notions and preliminary results necessary to understand the context and the contributions of the paper.

\subsection{State-space systems and morphisms} \label{morphism_session}
A continuous time state-space system is given by the following two equations 
\begin{equation}
\label{state space system}
	\left\{
	\begin{aligned}
		\dot{{\bf z}}&=F({\bf z}, u),\\
		y&=h({\bf z}),
	\end{aligned} \right.
\end{equation} 
where $u \in \mathcal{U}$ is the {\it input}, ${\bf z} \in \mathcal{Z}$ is the {\it internal state} and $F:\mathcal{Z}\times\mathcal{U}\rightarrow\mathcal{Z}$ is called the {\it state map}. The first equation is called the {\it state equation} while the second one is usually referred to as the {\it observation equation}. The solutions of \eqref{state space system} (when available and unique) yield an input/output map that is by construction causal and time-invariant. State-space systems will be sometimes denoted using the triplet $(\mathcal{Z}, F, h)$.

\begin{definition}
\label{morphism}
A map $f:\mathcal{Z}_1\rightarrow\mathcal{Z}_2$ is called a {\it system morphism} (see \cite{RC16}) between the continuous-time state-space systems $(\mathcal{Z}_1, F_1,h_1)$ and $(\mathcal{Z}_2, F_2,h_2)$ if it satisfies the following two properties:
\begin{description}
	\item[(i)] {\it System equivariance}: $f(F_1({\bf z}_1,u))=F_2(f({\bf z}_1),u)$, for all ${\bf z}_1\in\mathcal{Z}_1$ and $u\in\mathcal{U}$.
	\item[(ii)] {\it Readout invariance}: $h_1({\bf z}_1)=h_2(f({\bf z}_1))$ for all ${\bf z}_1\in\mathcal{Z}_1$.
\end{description}
\end{definition}
As a direct consequence of this definition, composition of system morphisms is again a system morphism. In the case $f$ is invertible and $f^{-1}$ is also a morphism, we say that $f$ is a system isomorphism. An elementary but very important fact is that if
$f:\mathcal{Z}_1\rightarrow\mathcal{Z}_2$ is a linear system-equivariant map between $(\mathcal{Z}_1, F_1,h_1)$ and $(\mathcal{Z}_2, F_2,h_2)$ ($\mathcal{Z}_1 $ and $\mathcal{Z}_2 $ are in this case vector spaces) then, for any solution
 ${\bf z}_1\in C ^1(I, \mathcal{Z}_1)$ of the state equation associated to $F_1 $ and to the input $u \in C ^1(I, \mathcal{U})$, with $I \subset \mathbb{R} $ an interval, its image $f \circ {\bf z} _1 \in C ^1(I, \mathcal{Z}_2)$  is a solution for the state space system associated to $F_2$ with the same input. Indeed, for any $t \in I $ we have, by the linearity and the system equivariance of $f$:
\begin{equation*}
	\frac{d}{dt}[f({\bf z}_1(t))]=Df({\bf z}_1(t))\cdot \dot{{\bf z}}_1(t)=f(\dot{{\bf z}}_1(t))
	=f(F_1({\bf z}_1(t),u(t)))=F_2(f({\bf z}_1(t)),u(t)).
\end{equation*}
This fact has as an important consequence that, in general, input/output systems {\it are not uniquely identified} since all the system-isomorphic state-space systems yield the same input/output map. 

\subsection{Hamiltonian and port-Hamiltonian systems}

Hamiltonian systems are dynamical systems whose behavior is governed by Hamilton's variational principle. Even though these autonomous systems can be in general formulated on any symplectic manifold \cite{Abraham1978}, we will restrict in this paper to the case in which the phase space is the even dimensional vector space $\mathbb{R}^{2n}$ endowed with the Darboux canonical symplectic form. In this case, the {\it Hamiltonian system} determined by the {\it Hamiltonian function} $H \in C^1(\mathbb{R}^{2n} ) $ is given by the differential equation
\begin{equation}
\label{hamiltonian_def}
	\dot{{\bf z}}=\mathbb{J}\frac{\partial H}{\partial {\bf z}},
\end{equation}
where  ${\displaystyle \mathbb{J}=
\begin{bmatrix}
	0 & \mathbb{I}_n\\
	-\mathbb{I}_n & 0
\end{bmatrix}}$ is the so-called the {\it canonical symplectic matrix}. Note that $-\mathbb{J}=\mathbb{J}^{T}=\mathbb{J}^{-1}$ and hence endows $\mathbb{R}^{2n}$ also with a complex structure. In this paper, we will denote the canonical symplectic matrix as $\mathbb{J}$, unless the context requires to specify the dimension, in which case we denote it by $\mathbb{J}_n$.

A {\it linear} Hamiltonian system is determined by a quadratic Hamiltonian function $H({\bf z})=\frac{1}{2}{\bf z}^{T}Q{\bf z}$, where ${\bf z}\in\mathbb{R}^{2n}$ and $Q \in \mathbb{M}_{2n}$ is a square matrix that without loss of generality can be assumed to be symmetric. In this case, Hamilton's equations (\ref{hamiltonian_def}) reduce to 
\begin{equation}
\label{linear hamiltonian system}
	\dot{{\bf z}}=\mathbb{J}Q{\bf z}.
\end{equation}

\noindent {\it Port-Hamiltonian systems} (see \cite{vanderSchaft2014}) are state-space systems that generalize autonomous Hamiltonian systems to the case in which external signals or inputs control in a time-varying way the dynamical behavior of the Hamiltonian system.  The family of input-state-output port-Hamiltonian systems are those port-Hamiltonian systems with no algebraic constraints on the state-space variables, and where the flow and effort variables of the resistive, control and interaction ports are split into conjugated pairs. In such cases, the implicit representation may be proved (see \cite{vanderSchaft2014}) to be equivalent to the following explicit form:
\begin{equation}
\label{general ph system}
\left\{
	\begin{aligned}
		\dot{{\bf x}}&=[J({\bf x})-R({\bf x})]\frac{\partial H}{\partial {\bf x}}({\bf x})+g({\bf x})u,\\
		y&=g^{T}({\bf x})\frac{\partial H}{\partial {\bf x}}({\bf x}),
	\end{aligned} \right.
\end{equation}
where $(u,y)$ is the input-output pair (corresponding to the control and output conjugated ports), $J({\bf x})$ is a skew-symmetric interconnection structure and $R({\bf x})$ is a symmetric positive-definite dissipation matrix. 
Our work concerns {\it linear} port-Hamiltonian systems in the {\it normal form} which we define now: a linear port-Hamiltonian system (\ref{general ph system}) is in normal form if the skew-symmetric matrix $J$ is constant and equal to the canonical symplectic matrix $\mathbb{J}$, the Hamiltonian matrix $Q$ is symmetric positive-definite, and the energy dissipation matrix $R=0$, in which case \eqref{general ph system} takes the form:
	\begin{equation}
\label{ph_definition}
		\left\{
		\begin{aligned}
			\dot{{\bf z}}&=\mathbb{J}Q{\bf z}+Bu,\\
			y&=B^{T}Q{\bf z},
		\end{aligned} \right.
	\end{equation} 
with ${\bf z}\in\mathbb{R}^{2n}$, $u,y\in\mathbb{R}$, and where $B\in\mathbb{R}^{2n}$ specifies the interconnection structure simultaneously at the input and output levels. By definition, such systems are fully determined by the pair $(Q,B)$, and hence we define by  
\begin{equation}
\label{set ph params}
\Theta_{PH_n}:=\left\{(Q,B)|0<Q\in\mathbb{M}_{2n}, Q=Q^T, B\in\mathbb{R}^{2n}\right\}
\end{equation}
the space of {\it paramters} of (\ref{ph_definition}). Let $\theta_{PH_n}:\Theta_{PH_n}\rightarrow PH_n$ the map that associates to the parameter $(Q,B)\in \theta_{PH_n}$ the corresponding port-Hamiltonian state space system. For convenience,  {\it we shall often use $(Q,B)$ to denote elements in $PH_n$ unless there is a risk of confusion}. Note that the condition $Q> 0$ implies that the origin is a Lyapunov stable equilibrium of  \eqref{linear hamiltonian system}. All these systems have the existence and uniqueness of solutions property and hence determine a family of {\it input/output systems}, also known as {\it filters}, that will be denoted by $\mathcal{PH}_n$. More specifically, the elements in $\mathcal{PH}_n$ are maps $U_{(Q,B)}: C ^1([0,1])\times \mathbb{R}^{2n} \longrightarrow C ^1([0,1]) $ given by
\begin{equation*}
\begin{array}{cccc}
U_{(Q,B)}: &C ^1([0,1])\times \mathbb{R}^{2n} &\longrightarrow &C ^1([0,1])\\
 	&\left(u, \mathbf{x}_0\right)&\longmapsto &U_{(Q,B)}(u, \mathbf{x} _0)_t=B ^{T}Qe^{\mathbb{J}Q t} \left[\int _0 ^t e^{- \mathbb{J} Q s}B u (s)\, d s + \mathbf{x} _0\right], \quad \mbox{$t \in [0,1]$.}
\end{array}
\end{equation*}
Note that $PH_n$ includes as a special case linear observations of autonomous linear Hamiltonian systems (case $B=0 $). Note that as a manifold $\Theta_{PH_n}= \mathcal{S}_{2n}^+ \times \mathbb{R}^{2n} $, where $\mathcal{S}_{2n}^+$ denotes the space of symmetric positive-definite matrices (SPD). We recall that $\mathcal{S}_{2n}^+$ has a natural differentiable manifold structure whose tangent space at any point is the vector space of symmetric matrices $\mathcal{S}_{2n}$ (see \cite{Quang:Murino}, and references therein).

Port-Hamiltonian systems are also closely linked to the so-called {\it affine Hamiltonian input-output systems} that have been considered as a natural extension of Hamiltonian systems with external forces and studied extensively in the literature (see \cite{Crouch1987} for the deterministic case and \cite{bismut1982mecanique, lazaro2008stochastic} for stochastic extensions), which take the form
  \begin{equation}
 	\label{affine_Hamiltonian_input_output_systems}
 	\left\{
 	\begin{aligned}
 		\dot{{\bf x}}&=X_{H}({\bf x})+X_{g}({\bf x})u,\\
 		\tilde{y}&=g({\bf x}),
 	\end{aligned} \right.
 \end{equation} 
where $X_H$ and $X_g$ are the Hamiltonian vector fields of $H, g \in  \in C^1(\mathbb{R}^{2n} ) $. In the linear case, (\ref{affine_Hamiltonian_input_output_systems}) reduces to

 \begin{equation}
	\label{linear_affine_Hamiltonian_input_output_systems}
	\left\{
	\begin{aligned}
		\dot{{\bf z}}&=\mathbb{J}Q{\bf z}-\mathbb{J}Bu,\\
		\tilde{y}&=B^T{\bf z},
	\end{aligned} \right.
\end{equation}

The relation between (\ref{linear_affine_Hamiltonian_input_output_systems}) and (\ref{ph_definition}) is that $\dot{\tilde{y}}=B^{T}\dot{{\bf z}}=B^{T}\mathbb{J}Q{\bf z}=(-\mathbb{J}B)^{T}Q{\bf z}$, showing that the time derivative of the affine Hamiltonian input-output system has a port-Hamiltonian structure. Note that in the last equality, we used that $B ^T \mathbb{J} B=0  $ since $\mathbb{J}  $ is antisymmetric.

Consider now a general linear single-input/single-output system that takes the form 
\begin{equation}
\label{linear systems general}
	\left\{
	\begin{aligned}
		\dot{{\bf x}}&=A{\bf x}+Bu,\\
		y&=C^{T}{\bf x},
	\end{aligned} \right.
\end{equation}
where $A\in\mathbb{M}_{n}$, $B,C\in\mathbb{R}^{n}$. 
Very often in control theory, it is the so-called transfer matrix rather than the input/output system that is studied. The transfer matrix $G(s)$ of  (\ref{linear systems general}) is defined as $G(s) = C(\mathbb{I}s-A)B$, which converts the differential equations in the time domain to an algebraic equation in the Laplace frequency domain. It can be proved that the transfer matrix of systems (\ref{ph_definition}) satisfies $G(s)=-G(-s)$ and that of systems (\ref{linear_affine_Hamiltonian_input_output_systems}) satisfies $G(s)=G(-s)$. The converse statements also hold for canonical realizations (see definition in next section)  \cite{Brockett1972}, \cite{maschke1992}. These facts exhibit a strong indication that systems (\ref{ph_definition}) and (\ref{linear_affine_Hamiltonian_input_output_systems}) carry intrinsic symmetries to be explicitly characterized.

\subsection{Controllability and observability}

Given a general linear system like \eqref{linear systems general}, we recall that its {\it controllability} and {\it observability matrices} are defined by 
\begin{equation*}
\begin{bmatrix}
		B\ |\ AB\ |\ \dots\ |\ A^{n-1}B
	\end{bmatrix}
\quad \mbox{and} \quad
	\begin{bmatrix}
		C^{T}\\C^{T}A\\\vdots \\C^{T}A^{n-1}
	\end{bmatrix}, \quad \mbox{respectively.}
\end{equation*}

The system is called {\it controllable} (respectively, {\it observable}) if its controllability (respectively, observability) matrix has full rank. Any linear controllable (respectively, observable) system can be transformed into the so-called controllable (respectively, observable) canonical forms by using appropriate linear system isomorphisms (see \cite{Polderman1998}). Conversely, systems in these canonical forms are automatically controllable (respectively, observable). In the next section, we characterize the controllable/observable/canonical systems in the linear port-Hamiltonian category.

Controllability and observability are intertwined concepts in the linear port-Hamiltonian category. Indeed, it can be proved (see \cite{Medianu2013}) that a linear port-Hamiltonian system without dissipation is controllable and $\det(Q)\neq0$, then it is also observable. Conversely, if it is observable, then this implies that $\det(Q)\neq0$ and it is also controllable (see \cite{Medianu2013}). As it is customary in systems theory, we say a linear port-Hamiltonian system in normal form is {\it canonical} if it is both controllable and observable. In view of the results that we just recalled, if $\det(Q)\neq 0$, then either controllability or observability is equivalent to the system being canonical. Furthermore, it can be shown that being canonical is a generic property, that is, the set of canonical systems forms an open and dense subset. We shall denote by $PH_n^{can}\subset PH_n $ the subset of $PH_n $ made of canonical linear port-Hamiltonian systems. Later on in the paper, the significance of these observations will become apparent.

\subsection{The symplectic Lie group and its Lie algebra} 

A square matrix $S \in \mathbb{M}_{2n}$ in dimension $2n $ is called {\it symplectic} if  it satisfies $S^{T}\mathbb{J}S=\mathbb{J}$. The set of all symplectic matrices forms a Lie group denoted by $Sp(2n,\mathbb{R})$. It is well-known see that if $S \in Sp(2n,\mathbb{R})$ then $\det S=\pm 1 $  and hence $Sp(2n,\mathbb{R})$ is a subgroup of the general linear group $GL(2n, \mathbb{R})$. The Lie algebra $\mathfrak{sp} (2n,\mathbb{R}) $ of $Sp(2n,\mathbb{R})$ is given by the matrices $A \in \mathbb{M}_{2n} $ that satisfy the identity 
 $A^{T}\mathbb{J}+\mathbb{J}A=0$. Equivalently, $A \in \mathfrak{sp} (2n,\mathbb{R}) $ if and only if $A=\mathbb{J}R$, where $R \in \mathbb{M}_{2n}$ is symmetric. We will refer to the elements in $Sp(2n,\mathbb{R})$ as {\it symplectic matrices} and to those in $\mathfrak{sp} (2n,\mathbb{R}) $ as {\it infinitesimally symplectic}.

Notably, the eigenvalues of the elements in $\mathfrak{sp} (2n,\mathbb{R}) $ appear in specific patterns that are spelled out in the following classical proposition (see \cite[Section 3.1]{Abraham1978}).
\begin{proposition}
\label{eigenvalues}
The characteristic polynomial of any matrix in $A \in \mathfrak{sp} (2n,\mathbb{R}) $ is even. Thus, if $\lambda$ is an eigenvalue of $A$ then so are $-\lambda$, $\bar{\lambda}$, and $-\bar{\lambda}$.
\end{proposition}
The importance of this group in our developments is that the (constant) vector field associated with the Hamilton's equations \eqref{linear hamiltonian system} is an element in $\mathfrak{sp} (2n,\mathbb{R})$. Its flow determines a one-parameter subgroup of elements in $Sp(2n,\mathbb{R})$. We also introduce the unitary group $U(n,\mathbb{C})$, which consists of matrices $U\in\mathbb{M}_n(\mathbb{C})$ with $UU^*=U^*U=\mathbb{I}_n$, where $U^*$ denotes the conjugate transpose of $U$. We denote by \cite{de2006symplectic} $U(n)$ the image of $U(n,\mathbb{C})$ in $Sp(2n,\mathbb{R})$ by the monomorphism 
\begin{equation}
\label{monomorp un}
	A+iB\rightarrow\begin{bmatrix}
		A&-B\\
		B&A
	\end{bmatrix}.
\end{equation} 
The so called {\it 2-out-of-3 property} \cite{Arnold1989} implies that $U(n)=O(2n,\mathbb{R})\cap GL(n,\mathbb{C})\cap Sp(2n,\mathbb{R})$, and it is indeed the intersection of any two out of the three groups.

\subsection{Williamson's normal form}
\label{Williamson normal form section}
The following classical result can be found in \cite{williamson1936algebraic, williamson1937normal, Ikramov2018, de2006symplectic}.
\begin{theorem}
\label{Williamson's normal form}
Let $M \in \mathbb{M}_{2n}$ be a positive-definite symmetric real matrix. Then
\begin{description}
\item[(i)] There exists a symplectic matrix $S\in Sp(2n, \mathbb{R})$ such that ${\displaystyle M=S^{T}\begin{bmatrix}
		D&0\\
		0&D
	\end{bmatrix}}S$, with $D= {\rm diag}({\bf d})$ a $n$-dimensional diagonal matrix with positive entries and ${\bf d}=\left(d_1, \ldots, d_n\right)^T$
.\\
	\item[(ii)] The values $d_1, \ldots, d_n$ are independent, up to  reordering, on the choice of the symplectic matrix $S$ used to diagonalize $M$ .
	\item[(iii)] Assume $S$ and $S^{\prime}$ are two elements of $Sp(2n,\mathbb{R})$ such that $M=S^{T}\begin{bmatrix}
		D&0\\
		0&D
	\end{bmatrix}S=S^{\prime T}\begin{bmatrix}
	D&0\\
	0&D
\end{bmatrix}S^{\prime}$, where $D$ is as above, then $S(S^{\prime})^{-1}\in U(n)$.
\end{description}
\end{theorem}

Later in this paper, we always use the notation $D={\rm diag}({\bf d})$ to denote that $D$ is a diagonal matrix with diagonal entries given by the vector ${\bf d}=(d_1,\dots,d_n)^T$. The elements $d_i$ in the above theorem are called the {\it symplectic eigenvalues} of $M$ since they are also the eigenvalues of $\mathbb{J}M$.

\begin{remark}
\normalfont
	The above theorem can be generalized to {\it positive-semidefinite} real symmetric matrices. Indeed, it can first be shown that if the kernel of $M$ is a symplectic subspace of $\mathbb{R}^{2n}$ of dimension $2m$, then the statement of Theorem \ref{Williamson's normal form} still holds true holds with the only added feature that exactly $m$ of the diagonal entries in $D$ are equal to $0$ (see \cite{Son2022}). More generally, without the symplecticity assumption, all that it can be said is that there exists  $S\in Sp(2n,\mathbb{R})$ such that ${\displaystyle M=S^{T}\begin{bmatrix}
		D_1&0\\
		0&D_2
\end{bmatrix}}S$ where $D_1$ and $D_2$ may contain diagonal zero entries (see\cite{Idel2017, Egusquiza2022}). 
\end{remark}

\section{Controllable and observable Hamiltonian representations}

In this section, we state two representation results for linear port-Hamiltonian systems in normal form, which are the main building blocks in our learnability results. More precisely, we define two subfamilies of linear systems of the type \eqref{linear systems general}, that are respectively called controllable/observable Hamiltonian representations, that are by construction controllable/observable (Definition \ref{definition_of_controllable}). We subsequently show in Theorem \ref{main_theorem} that morphisms can be established between the elements in these families and those in the category ${PH}_n $ of normal form port-Hamiltonian systems.

As it will be spelled out later on in detail, the existence of these morphisms immediately guarantees that the complexity of the family of filters $\mathcal{PH}_n $ is actually not $\mathcal{O}(n^2)$, as it could be guessed from \eqref{ph_definition}, but $\mathcal{O}(n)$. However, the expressive power of our proposed representations is limited for non-canonical port-Hamiltonian systems. For example, the observable representation is guaranteed to capture all possible input-output dynamics of port-Hamiltonian systems (full expressive power), but it does not always produce port-Hamiltonian dynamics (fails to be structure-preserving). In the controllable case, structure preservation is guaranteed, but there is, in general, no full expressive power. Fortunately, for canonical port-Hamiltonian systems, all the morphisms that we shall introduce become isomorphisms, meaning that they are both structure-preserving and have full expressive power. Roughly speaking, the more canonical a port-Hamiltonian system is, the better the corresponding representations behave in terms of structure-preserving properties and expressive power. \par

The representations introduced below can be seen as a reparametrization of the elements $(Q,B) \in {PH}_n $ in terms  of a diagonal matrix $D= \operatorname{diag} (\mathbf{d})\in\mathbb{M}_{n}$, ${\bf d}\in\mathbb{R}^{n}$, and a vector ${\bf v}\in\mathbb{R}^{2n}$, where $D$ is obtained from Williamson's Theorem \ref{Williamson's normal form}  as ${\displaystyle Q=S^{T}\begin{bmatrix}
		D&0\\
		0&D
\end{bmatrix}S}$ and ${\bf v}=S^{-1}B$. 
This makes it obvious that the learning problem for port-Hamiltonian systems has parameter complexity of at most $\mathcal{O}(n)$ even if the Hamiltonian matrix has complexity $\mathcal{O}(n^2)$. 

We emphasize that even in the canonical situation, the availability of the controllable/observable representations does not yet provide a well-specified learning problem for this category since the invariance of these systems under system automorphisms implies the existence of symmetries (or degeneracies) in the parametrizations, which will be the focus of the next section. 

The proofs of all our results are provided in the appendices.

\begin{definition}\label{definition_of_controllable}
	Given ${\bf d}=\left(d_1, \ldots, d_n\right)^T\in {\Bbb R}^n$, with $d_i>0$, and ${\bf v}\in\mathbb{R}^{2n}$, we say that a $2n$-dimensional linear state space system is a {\it controllable Hamiltonian} (respectively, {\it observable Hamiltonian}) representation if it takes the form
	\begin{equation}
		\label{proposed_system_learner}
		\left\{
		\begin{aligned}
			\dot{{\bf s}}&=g^{ctr}_1({\bf d})\cdot {\bf s}+\left(
			0 , 0 ,\cdots, 0, 1\right)^{T}\cdot u,\\
			y&=g^{ctr}_2({\bf d}, {\bf v})\cdot {\bf s},
		\end{aligned} \right. 
		\left(\mbox{resp., } \left\{
		\begin{aligned}
			\dot{{\bf s}}&=g^{obs}_1({\bf d})\cdot {\bf s}+g^{obs}_2({\bf d},{\bf v})\cdot u,\\
			y&=\left(
			0 , 0 ,\cdots, 0, 1\right)\cdot {\bf s},
		\end{aligned} \right. \right)
	\end{equation} 
	where $g^{ctr}_1({\bf d}) \in \mathbb{M}_{2n }$ and $g^{ctr}_2({\bf d}, {\bf v}) \in \mathbb{M}_{1,2n}$ (respectively, $g^{obs}_1({\bf d}) \in \mathbb{M}_{2n }$ and $g^{obs}_2({\bf d}, {\bf v}) \in \mathbb{R}^{2n}$) are constructed as follows:
	\begin{description}
		\item [(i)] Given ${\bf d} \in {\Bbb R}^n$, let  $\left\{a _0, a _1, \ldots , a_{2n-1}\right\}$ be the real coefficients that make ${\lambda}^{2n}+\sum_{i=0}^{2n-1}a_i\cdot{\lambda}^{i}=({\lambda}^2+d_1^2)({\lambda}^2+d_2^2)\dots({\lambda}^2+d_n^2)$ an equality between the two polynomials in $\lambda$. Let $a_{2n}=1$ by convention. Note that the entries $a_{i}$ with an odd index $i $ are zero. Define: \begin{equation*}
			g^{ctr}_1({\bf d}):=\begin{bmatrix}
				0&1&0&\dots&0 \\
				0&0&1&\dots&0 \\
				\vdots&\vdots&\ddots&\vdots&\vdots\\
				0&0&0&\dots&1 \\
				-a_0&-a_1&-a_2&\dots&-a_{2n-1}
			\end{bmatrix}_{2n\times 2n},
		\end{equation*}
		(respectively, $g^{obs}_1({\bf d})=g^{ctr}_1({\bf d}) ^{\top}$).
		\item[(ii)] Given ${\bf d}$ and ${\bf v}$, then \begin{equation*}
			\begin{aligned}
				g^{ctr}_2({\bf d},{\bf v})&:=\begin{bmatrix}
					0\:\:c_{2n-1}\:\:0\:c_{2n-3}\:\hdots\:\:0\:\:c_{1}
				\end{bmatrix}, \quad \mbox{(resp., $g^{obs}_2({\bf d},{\bf v})= g^{ctr}_2({\bf d},{\bf v}) ^{\top}$)} 
			\end{aligned}
		\end{equation*} where 
		\begin{equation*}
			c_{2k+1}={\bf v}^{T}\begin{bmatrix}
				F_k&0\\
				0&F_k
			\end{bmatrix} {\bf v},
		\end{equation*} for $k=0,\dots,n-1$, and
		\begin{equation*}
			F_k=\begin{bmatrix}
				f_1&&\\
				&f_2&&\giantzero\\
				&&\ddots\\
				&\giantzero&&f_{n-1}\\
				&&&&f_n\\
			\end{bmatrix}
		\end{equation*}
		with $f_l=d_l\cdot\sum_{\substack{j_1,\dots,j_k\neq l\\1\leq j_1<\dots<j_k\leq n}}\big(d_{j_1}d_{j_2}\cdots d_{j_k}\big)^2$, $l=1,\dots,n$.
	\end{description}
	We denote $CH_n$ (respectively, $OH_n$) the set of all systems of the form (\ref{proposed_system_learner}), and we call them {\it controllable Hamiltonian} (respectively, {\it observable Hamiltonian}) representations. The symbol $\mathcal{CH}_n$ (respectively, $\mathcal{OH}_n$) denotes the set of input/output systems induced by the state space systems in $CH_n$ (respectively, $OH_n$). We emphasize that the elements of both $CH_n$ and $OH_n$ can be parameterized with the set 
	$$\Theta_{CH_n}=\Theta_{OH_n}:=\left\{({\bf d},{\bf v})|d_i>0, {\bf v}\in\mathbb{R}^{2n}\right\}.$$
Sometimes later on in the paper  we shall write $a_i({\bf d})$ and $c_j({\bf d},{\bf v})$ to indicate that $a_i$ and $c_j$ are functions of ${\bf d}$ and ${\bf v}$.
\end{definition}

\medskip

\begin{remark}
	\label{adjoint}
	\normalfont
	Observe that the controllable and the observable Hamiltonian representations of port-Hamiltonian systems are closely related to each other. The controllable Hamiltonian matrix $g^{ctr}_1$ is the transpose of the observable Hamiltonian matrix $g^{obs}_1$. Moreover, as can be directly observed from the construction, the input and readout matrices of the two representations, that is, $g^{ctr}_2$ and $g^{obs}_2$, are transpose of each other.
\end{remark}

Consider now the maps $\theta_{CH_n}:\Theta_{CH_n} \rightarrow CH_n $ and $\theta_{OH_n}:\Theta_{OH_n} \rightarrow OH_n $ that associate to each parameter values the corresponding state-space system. Note that the elements in $CH _n $ (respectively, in $OH _n $) of the form (\ref{proposed_system_learner}) are in canonical controllable (respectively, observable) form in the sense of \cite{sontag:book} and they are hence controllable (respectively, observable). Our main result below establishes a relationship between port-Hamiltonian systems and controllable (respectively, observable) Hamiltonian representations as defined above, which will be used later on for considerations on the structure preservation and expressiveness in the modeling of $PH _n $.

\begin{theorem}
	\label{main_theorem}
	\begin{description}
		\item [(i)] There exists, for each $S\in Sp(2n,\mathbb{R})$,  a map
		\begin{equation*}
			\begin{array}{cccc}
				\varphi_S: &CH _n& \longrightarrow & {PH}_n\\
				&\theta_{CH _n}({\bf d},{\bf v}) &\longmapsto & \theta_{PH_n}\left(S^{T}\begin{bmatrix}
					D&0\\
					0&D
				\end{bmatrix}S, S^{-1}{\bf v}\right),
			\end{array}
		\end{equation*}
		with $D={\rm diag}({\bf d})$, such that the controllable Hamiltonian system $\theta_{CH _n}({\bf d},{\bf v}) \in CH_n $  and the port-Hamiltonian image $\varphi_S\left(\theta_{CH _n}({\bf d},{\bf v})\right) \in {PH}_n$ are linked by a linear system morphism $f_S^{({\bf d},{\bf v})}:\mathbb{R}^{2n}\rightarrow\mathbb{R}^{2n}$.
		\item [(ii)] Given a port-Hamiltonian system $\theta_{PH _n}(Q,B)\in {PH}_n$,  there exists an explicit linear system morphism $f^{(Q,B)}:\mathbb{R}^{2n}\rightarrow\mathbb{R}^{2n}$ between the state space of $\theta_{PH _n}(Q,B)\in {PH}_n$ and that of an observable Hamiltonian system $\theta_{OH _n}({\bf d},{\bf v})\in OH _n$, where $({\bf d},{\bf v}) \in \Theta_{OH _n}$ is determined by the Williamson's normal form decomposition of $Q$ determined by $S\in Sp(2n,\mathbb{R})$, that is,  ${\displaystyle Q=S^{T}\begin{bmatrix}
				D&0\\
				0&D
		\end{bmatrix}}S$, $D={\rm diag}({\bf d})$ and ${\bf v}=S\cdot B$.
	\end{description}
\end{theorem}

\medskip

\begin{remark}\label{non_uniqueness}
	\normalfont
	We emphasize that given $(Q,B)\in \Theta_{PH _n}$, the pair $({\bf d}, {\bf v})\in \Theta_{CH_n}/\Theta_{OH _n}$ is not uniquely determined by Williamson's decomposition. This can be seen from Theorem \ref{Williamson's normal form} because the element $S \in Sp(2n,\mathbb{R})$ in its statement is not unique and the entries $d_i$ of $\mathbf{d} $ are independent of $S$ up to their ordering.
\end{remark}

\begin{remark}[{\bf Controllability, observability, and invertibility}]~
\label{Controllability, observability, and invertibility}
\normalfont
\begin{description}
	\item [(i)] In the proof of the theorem above (available in the Appendix), we define the linear system morphism $f_S^{({\bf d},{\bf v})}:\mathbb{R}^{2n}\rightarrow\mathbb{R}^{2n}$ as ${\bf z}=f_S^{({\bf d},{\bf v})}({\bf s}):=L{\bf s}$ and an explicit construction of the matrix $L$ is provided. It turns out that, the matrix $L$ is invertible if and only if the image port-Hamiltonian system (\ref{ph_definition}) is controllable, or equivalently, observable. Indeed, using the same notation as in the proof of Theorem \ref{main_theorem}, we have 
\begin{equation*}
		L=S^{-1}\begin{bmatrix}
			L_1{\bf v}&L_2{\bf v}&\cdots&L_{2n}{\bf v}
	\end{bmatrix}=\begin{bmatrix}
	S^{-1}L_1{\bf v}&S^{-1}L_2{\bf v}&\cdots&S^{-1}L_{2n}{\bf v}
\end{bmatrix},
	\end{equation*} where \begin{align*}
	S^{-1}L_{2n-k}{\bf v}&=S^{-1}\left[\left(\mathbb{J}_n\begin{bmatrix}
		D&0\\
		0&D
	\end{bmatrix}\right)^{k}+a_{2n-1}\cdot \left(\mathbb{J}_n\begin{bmatrix}
	D&0\\
	0&D
\end{bmatrix}\right)^{k-1}+\dots+a_{2n-k}\cdot \mathbb{I}_{2n}\right]\cdot{\bf v}\\
							&=S^{-1}\left((\mathbb{J}_nS^{-T}QS^{-1})^{k}+a_{2n-1}\cdot (\mathbb{J}_nS^{-T}QS^{-1})^{k-1}+\dots+a_{2n-k}\cdot \mathbb{I}_{2n}\right)\cdot SB\\
							&=S^{-1}\left((S\mathbb{J}_nQS^{-1})^{k}+a_{2n-1}\cdot (S\mathbb{J}_nQS^{-1})^{k-1}+\dots+a_{2n-k}\cdot \mathbb{I}_{2n}\right)\cdot SB\\
							&=\left((\mathbb{J}_nQ)^{k}+a_{2n-1}\cdot (\mathbb{J}_nQ)^{k-1}+\dots+a_{2n-k}\cdot \mathbb{I}_{2n}\right)\cdot B.
\end{align*} 
Therefore, $L$ can be transformed by elementary column operations into the controllability matrix of (\ref{ph_definition}) and hence $L$ being invertible, i.e. the two systems being isomorphic, is equivalent to the controllability matrix of (\ref{ph_definition}) having full rank (regardless of the choice of $S\in Sp(2n,\mathbb{R})$), which is again equivalent to (\ref{ph_definition}) being canonical. Additionally, the condition for $f_S^{({\bf d},{\bf v})}$ to be invertible can also be formulated in terms of $D$ and $\mathbf{v}$ directly, which we will discuss in Subsection \ref{sub3}.
	\item [(ii)] Systems in $CH_n$ are by construction in controllable canonical form, and are therefore always controllable. If the image system (\ref{ph_definition}) by $\varphi _S$ that we want to learn is controllable (or equivalently, observable), then by the previous point $L$ is necessarily an invertible matrix which means that (\ref{proposed_system_learner}) and (\ref{ph_definition}) are isomorphic systems by construction. As a consequence, (\ref{proposed_system_learner}) is not only controllable but also observable. 
\end{description}
\end{remark}

\medskip

\begin{remark}[{\bf Application to structure-preserving system learning}]~\\
\normalfont
As a corollary of the previous result, we can use controllable Hamiltonian representations to learn port-Hamiltonian systems in an efficient and structure-preserving fashion. Indeed, given a realization of a port-Hamiltonian system, a system of the type $\theta_{CH _n}({\bf d},{\bf v})\in CH_n$ can be estimated using an appropriate loss (see Section \ref{numerics}). A representation of this type is more advantageous than the original port-Hamiltonian one for two reasons:
\begin{description}
\item [(i)] The {\it model complexity} of the controllable Hamiltonian representation is only of order $\mathcal{O}(n)$, as opposed to $\mathcal{O}(n^2)$ for the original port-Hamiltonian one.
\item [(ii)] This learning scheme is automatically {\it structure-preserving}. Indeed, once a system $\theta_{CH _n}({\bf d},{\bf v})\in CH_n$ has been estimated for a given realization, we have shown that there exists a family of linear morphisms, each of which is between the state space of $\theta_{CH _n}({\bf d},{\bf v})\in CH_n$ and some $\theta_{PH _n}(Q,B)\in PH_n$, such that any solution of (\ref{proposed_system_learner}) is automatically a solution of some system in ${PH}_n$. Hence, even in the presence of estimation errors for  $({\bf d},{\bf v})\in\Theta_{CH_n}$, the solutions of $\theta_{CH _n}({\bf d},{\bf v})$ still correspond to a port-Hamiltonian system and hence this structure is {\it preserved} by the learning scheme. 
\end{description}
\end{remark}

\medskip

\begin{remark}[{\bf System learning and expressive power}]~\\
\normalfont
Expressive power is an important property of any machine learning paradigm. As a continuation of the previous remarks, we emphasize that there is an important relation between the controllability of a system in ${PH}_n$ and the expressive power of the corresponding representation in $CH_n$. Indeed, if (\ref{ph_definition}) is controllable, by point (ii) in Remark \ref{Controllability, observability, and invertibility}, the corresponding preimage system $\theta_{CH _n}({\bf d},{\bf v})\in CH_n$ can capture all possible solutions of (\ref{ph_definition}), which amounts to the learning scheme based on $\Theta_{CH_{n}}$ having full expressive power. To see this, let ${\bf z}_0$ be an initial state of the controllable system $\theta_{PH _n}(Q,B)\in {PH}_n$ in (\ref{ph_definition}). Since in that case we can find an invertible system isomorphism $f _S^{({\bf d},{\bf v})} $ that links it to some $\theta_{CH_n}({\bf d},{\bf v})\in \Theta_{CH_n}$, there exists some corresponding initial state ${\bf s}_0= \left(f _S^{({\bf d},{\bf v})}\right)^{-1}({\bf z}_0)$. Then, by Theorem \ref{main_theorem} and the uniqueness of the solutions of ODEs, the solution of (\ref{proposed_system_learner}) with initial state ${\bf s}_0$ is a representation of the solution of (\ref{ph_definition}) with initial state ${\bf z}_0$. However, if (\ref{ph_definition}) fails to be controllable (i.e. $f _S^{({\bf d},{\bf v})} $ not invertible), then such an initial condition ${\bf s}_0$ may not exist. As a rule of thumb, the more controllable a system of the type (\ref{ph_definition}) is, the higher the rank of $f _S^{({\bf d},{\bf v})}  $ is, and then the more expressive the corresponding controllable  Hamiltonian representations are.
\end{remark}

\medskip

\begin{remark}[{\bf Expressive power and structure-preservation of the observable Hamiltonian representation}]~\\
\normalfont
We emphasize that systems in $OH_n$ always has {\it full expressive power} guaranteed by the system morphism in Theorem \ref{main_theorem}. This implies that any input-output dynamics generated by the original port-Hamiltonian system will be captured by the any of the observable Hamiltonian representations in the statement. However, unlike in the controllable case, the system morphism is between $\theta_{PH _n}(Q,B)\in {PH}_n$ and $\theta_{OH_n}({\bf d},{\bf v})\in\Theta_{OH_n}$. Therefore, unless $(Q,B)$ is canonical, in which case the morphism becomes an isomorphism,  we {\it cannot, in general, assert the structure-preserving property of this representation}.
\end{remark}

\medskip

\begin{remark}[{\bf Positive semi-definite Hamiltonians}]~\\
	\normalfont
	The above results can be easily generalized to positive semi-definite Hamiltonians (PSD) Hamiltonians with the aid of the generalized Williamson's theorem in the references \cite{Son2022,Idel2017, Egusquiza2022} that we briefly discussed in Section \ref{Williamson normal form section}. In general, the number of unknown parameters in the vector $\mathbf{d} $ is doubled (because of the matrices $D _1 $ and $D_2 $  that appear in this case), and their relation with the coefficients $\left\{a _0, a _1, \ldots , a_{2n-1}\right\}$ has to be modified accordingly, that is, ${\lambda}^{2n}+\sum_{i=0}^{2n-1}a_i\cdot{\lambda}^{i}=({\lambda}^2+d_1d_{n+1})({\lambda}^2+d_2d_{n+2})\dots({\lambda}^2+d_nd_{2n})$, where some of the $d_i$'s could be $0$. The expression for $g^{ctr}_1({\bf d})$ remains the same, whereas the expression of $\begin{bmatrix}
		F_k&0\\
		0&F_k
	\end{bmatrix}$ in $g^{ctr}_2({\bf d},{\bf v})$ becomes $\begin{bmatrix}
		F_{k,0}&0\\
		0&F_{k,1}
	\end{bmatrix}$, where \begin{equation*}
		F_{k,p}=\begin{bmatrix}
			f_{1,p}&&\\
			&f_{2,p}&&\giantzero\\
			&&\ddots\\
			&\giantzero&&f_{n-1,p}\\
			&&&&f_{n,p}\\
	\end{bmatrix} \end{equation*} and $f_{l,p}=d_{np+l}\cdot\sum_{\substack{j_1,\dots,j_k\neq l\\1\leq j_1<\dots<j_k\leq n}}d_{j_1}d_{j_2}\cdots d_{j_k}d_{j_{n+1}}d_{j_{n+2}}\cdots d_{j_{n+k}}$ for $p=0,1$. In this paper, we mainly deal with positive definite $Q$, since the nondegeneracy of a positive semi-definite $Q$ destroys the symmetries studied later on in Section {\ref{quotient}}.
\end{remark}

\medskip

\begin{remark}[{\bf Symmetries of the Hamiltonian representations}]~\\
\normalfont
The parameterizations of the systems in $CH_n$ and $OH_n$ exhibit obvious symmetries. For example, the functions $g^{ctr}_1({\bf d})$ and $g^{obs}_1({\bf d})$ are invariant under the permutation of the diagonal entries $d_i$. Moreover, $g^{ctr}_2({\bf d},{\bf v})$ (similarly for $g^{obs}_2({\bf d},{\bf v})$) contains entries $c_{2k+1}$ of the form ${\bf v}^{T}\begin{bmatrix}
	F_k&0\\
	0&F_k
\end{bmatrix} {\bf v}=\sum_{i=1}^{n}F^{(i)}_k\cdot\bigg[v_i^2+v_{n+i}^2\bigg]$, which is in particular invariant under the rotation of the planes spanned by the $i$-th and $(n+i)$-th entries of ${\bf v}$. These observations will be central in the next section, in which we shall show that these and other symmetries of the representations in $CH_n$ or $OH_n$ are closely related to the system automorphism group of the space ${PH}_n$.
\end{remark}

\section{Unique identification of linear port-Hamiltonian systems}
\label{quotient}

In this section, we study the unique identifiability of input-output dynamics of linear port-Hamiltonian systems in normal form. Such a characterization is obviously needed to solve the model estimation problem. The rationale is that, in applications, we only have access to input/output data, and different state space systems in $PH_n$ can induce the same filter that produces that data. This fact has important implications when it comes to the learning of port-Hamiltonian systems out of finite-sample realizations of a given data-generating process $(Q,B)\in PH_n$ because such degeneracy makes impossible the exact recovery of $(Q,B)\in{PH}_n$ in that context, no matter how good the properties of the algorithm used for that task are or how much data we have at our disposal.  This observation indicates that it is not in the space ${PH}_n$ that we should look at for unique identification but the quotient space associated to ${PH}_n$ with respect to certain equivalence relation $\sim_{filter} $ that uniquely identifies port-Hamiltonian filters, that is, $\mathcal{PH}_n\cong PH_n/\sim_{filter}$. However, as we shall see later on in Subsection \ref{sub0}, the presence of non-canonical systems in $PH_n $ makes it in general difficult to directly characterize the quotient space $PH_n/\sim_{filter}$.

As we pointed out after Definition \ref{morphism}, all the system-isomorphic state-space systems yield the same filter, while a filter can be realized by state-space systems that are not system isomorphic. This means that the equivalence relation of system automorphism $\sim_{sys}$ is strictly stronger than $\sim_{filter}.$ Motivated by this fact, we study in Subsection \ref{sub0} how $\sim_{filter}$ and $\sim_{sys}$ are related in terms of controllable Hamiltonian representations (which by Theorem \ref{main_theorem} automatically induce filters in $\mathcal{PH}_n$), and in Subsection \ref{sub1}, we lower our expectations and characterize $PH_n/\sim_{sys}$ as an approximation to $PH_n/\sim_{filter}$. The term {\it approximation} in this sentence is justified because  $\sim_{filter}$ and $\sim_{sys}$ {\it coincide} when restricted to the subset of canonical port-Hamiltonian systems $PH_n^{can}$, which is open and dense in $PH_n$. Therefore unique identifiability can be achieved in there by studying $\sim_{sys}$.

On the other hand, recall that in the previous section, we established a link between $PH_n$ and the representation spaces $CH_n$ and $OH_n$ which, as we saw in Definition \ref{definition_of_controllable}, are both parametrized by the set 
\begin{equation}
\label{r neutral}
\Theta_{CH_{n}}=\Theta_{OH_{n}}=\left\{({\bf d},{\bf v})\mid {\bf v}\in\mathbb{R}^{2n}, \mathbf{d} \in \mathbb{R}^n, d_i>0, i \in \left\{1, \ldots, n\right\}\right\}.
\end{equation} 
Now, it is a natural question to ask what is the equivalence relation that corresponds to $\sim_{sys}$ on the parameter space $\Theta_{CH_{n}}$, and if it is possible to explicitly characterize the quotient space $PH_n/\sim_{sys}$ on $\Theta_{CH_n}$ in a certain sense. All these questions are addressed step-by-step in the following subsections. 

In Subsection \ref{sub0}, we provide sufficient and necessary conditions for two controllable Hamiltonian representations being $\sim_{filter}$-equivalent and $\sim_{sys}$-equivalent, respectively. In Subsection \ref{sub1}, we define an equivalence relation $\sim_{\star}$ on $\Theta_{CH_{n}}$ and we show that ${{PH}_n/\sim_{sys}}\cong{\Theta_{CH_{n}}/\sim_{\star}}$ (see Theorem \ref{characterization}). In Subsection \ref{sub2}, we characterize the equivalence classes ${{PH}_n/\sim_{sys}}$ and ${\Theta_{CH_n}/\sim_{\star}}$ as {\it Lie groupoid} orbit spaces.

In Subsection \ref{sub3}, we restrict our identification analysis to {\it canonical port-Hamiltonian systems} $PH_n^{can}$. We first show that the parameter subset $\Theta^{can}_{CH_{n}}\subset \Theta_{CH_{n}} $ that corresponds to $PH_n^{can}$ is open and dense in $ \Theta_{CH_{n}}$ as it is determined by certain generic non-resonance and nondegeneracy conditions. If we define on $\Theta_{CH_{n}}$ the equivalence relation $\sim_{sys}$ of system automorphisms of the corresponding controllable/observable Hamiltonian representations (see Definition \ref{sys_well_defined}), then it can be proved that, restricted to the canonical subset $\Theta^{can}_{CH_{n}}$, the equivalence relation $\sim_{\star}$ coincides with  $\sim_{sys}$, and hence 
\begin{equation*}
\label{chain isomorph}
{\mathcal{PH}^{can}_n}\cong{{PH}^{can}_n/\sim_{sys}}\cong{\Theta^{can}_{CH_n}/\sim_{\star}}\cong{\Theta^{can}_{CH_n}/\sim_{sys}},
\end{equation*}
where ${\mathcal{PH}^{can}_n} $ is the space of filters induced by systems in $PH_n^{can}$.

In Subsection \ref{sub5}, we prove that the fact that we restricted the above equivalence relations to canonical subsets allows us to characterize the corresponding quotients as orbit spaces with respect to a  {\it group} (as opposed to groupoids in the general unrestricted case) action, where the group is given by a semi-direct product $S_n\rtimes_{\phi}\mathbb{T}^{n}$ that will be specified in detail later on. Finally, in Subsection \ref{sub6}, we show that the orbit space $\Theta^{can}_{CH_n}/(S_n\rtimes_{\phi}\mathbb{T}^{n})$ can be explicitly identified as a smooth manifold $\mathcal{X}^{n}_{\uparrow}\times \mathbb{R}_{+}^{n}$ and endowed with global Euclidean coordinates, and hence 
\begin{equation*}
	\label{chain isomorph2}
{\mathcal{PH}^{can}_n}\cong{PH}^{can}_n/\sim_{sys}~\cong~\Theta_{CH_n}^{can}/\sim_{\star}~\cong~\Theta_{CH_n}^{can}/\sim_{sys}~\cong~\Theta_{CH_n}^{can}/(S_n\rtimes_{\phi}\mathbb{T}^{n})~\cong~\mathcal{X}^{n}_{\uparrow}\times \mathbb{R}_{+}^{n}.
\end{equation*}
 Consequently, canonical port-Hamiltonian dynamics can be identified fully and explicitly in either the controllable or the observable Hamiltonian representations (\ref{proposed_system_learner}) and learned by estimating a unique set of parameters in a smooth manifold that is obtained as a group orbit space.

\subsection{The unique identification problem for filters in $\mathcal{PH}_n$}
\label{sub0}

 In the context of model estimation/machine learning, we would like to characterize and identify the filters that constitute the elements in $\mathcal{PH}_n$. In Section \ref{morphism_session}, we have seen that two systems that are system isomorphic induce the same input-output dynamics, which indicates that these isomorphisms are redundancies/symmetries in $PH_n$. Our aim is to quotient out the symmetries given by system automorphisms and to investigate whether the quotient space uniquely identifies the filters in $\mathcal{PH}_n$. 
 
 \begin{definition}~
 	\begin{description}
 		\item[(i)] The fact that two systems $\theta_{PH _n}(Q_1,B_1)$ and $\theta_{PH _n}(Q_2,B_2)$ in ${PH}_n$ induce the same filter defines an equivalence relation in $PH_n$, which we denote by $(Q_1,B_1)\sim_{filter}(Q_2,B_2)$. Consequently, we have by definition $\mathcal{PH}_n=PH_n/\sim_{filter}$, which we call the unique identifiability space.
 		\item[(ii)] We observe that $\theta_{PH _n}(Q_1,B_1)$ and $\theta_{PH _n}(Q_2,B_2)$ in ${PH}_n$ are linearly system isomorphic according to Definition \ref{morphism} if and only if there exists an {\it invertible} matrix $L$ such that
 		\begin{equation}
 			\label{automorphism_condition}
 			\left\{
 			\begin{aligned}
 				L\mathbb{J}Q_1&=\mathbb{J}Q_2L\\
 				LB_1&=B_2\\
 				B^{T}_1Q_1&=B^{T}_2Q_2L.
 			\end{aligned} \right.
 		\end{equation}
 		It is straightforward to check that system isomorphisms determine an equivalence relation on ${PH}_n$. If $\theta_{PH _n}(Q_1,B_1)$ and $\theta_{PH _n}(Q_2,B_2)$ are system isomorphic, we write $(Q_1,B_1)\sim_{sys} (Q_2,B_2)$. We denote by ${PH}_n/\sim_{sys}$ the quotient space. The equivalence class in ${PH}_n/\sim_{sys}$ that contains the element $\theta_{PH _n}(Q,B)$ is denoted by $[Q,B]\in {PH}_n/\sim_{sys}$. 
 	\end{description}
 \end{definition}
  It is a natural question to ask about the relation between  $PH_n/\sim_{sys}$ and $\mathcal{PH}_n$, and if they are the same. Indeed, two distinct elements in $PH_n$ that are $\sim_{sys}$-equivalent always induce the same filter in $\mathcal{PH}_n$ (See Subsection \ref{morphism_session}), whereas as we see in the next Example \ref{sys_doesnt_mean_same_filter}, a filter in $\mathcal{PH}_n$ could be realized by two elements in $PH_n$ that are not $\sim_{sys}$-equivalent since filters identify exclusively the canonical part (that is, the minimal realization) see \cite{kalman1963}. Said differently, by going to the quotient space $PH_n/\sim_{sys}$, we remove some redundancies of the set $PH_n$ that yield the same input-output dynamics, but not all. 

\begin{example}
\label{sys_doesnt_mean_same_filter}
\normalfont
	Consider two systems $\theta_{PH _n}(Q_1,B_1), \theta_{PH _n}(Q_2,B_2) \in PH_n$ where 
	\begin{equation*}
		Q_1 = \begin{bmatrix}
			1&0&0&0\\
			0&1&0&0\\
			0&0&1&0\\
			0&0&0&1
		\end{bmatrix}, Q_2= \begin{bmatrix}
		1&0&0&0\\
		0&2&0&0\\
		0&0&1&0\\
		0&0&0&3
	\end{bmatrix}, B_1=B_2=\begin{bmatrix}
	1\\0\\0\\0
\end{bmatrix}.
	\end{equation*}
Both systems induce the same filter $y(u)_t=\int_{0}^{t}e^{t-s}u(s)ds+e^t\left(
0 , 0 ,\cdots, 0, 1\right)^{T}\cdot z_0$, where $z_0$ is the initial state. However, these two systems cannot be system isomorphic, since by \eqref{automorphism_condition} in that case there would exist an invertible $L$ such that $L\mathbb{J}Q_1=\mathbb{J}Q_2L$, and hence $\mathbb{J}Q_1$ would have the same set of eigenvalues as $\mathbb{J}Q_2$, which is not the case.
\end{example}

Despite the difficulties of uniquely identifying all filters in $\mathcal{PH}_n$, we can still uniquely identify the ``large'' subset $\mathcal{CH}^0_n\subset \mathcal{PH}_n$ which consists of filters induced by the set of controllable Hamiltonian representations initially at rest, that is, with zero initial condition. Recall that any controllable Hamiltonian representation automatically has the port-Hamiltonian structure by Theorem \ref{main_theorem} and therefore $\mathcal{CH}^0_{n}$ lies inside $\mathcal{PH}_n$ (Moreover, $\mathcal{CH}_n$ contains all the filters induced by canonical systems in $PH_n$). Our result below splits into two parts, where the first one characterizes when two elements in $CH_n$ are $\sim_{sys}$-equivalent, and the second one characterizes when they induce the same filter in $\mathcal{CH}^0_n$. In this way, it will be clear that inducing the same filter is a strictly weaker condition than being system isomorphic. We  present a lemma and a definition before we state our main result in this section.

\begin{lemma}\label{sys_on_CHn_and_OH_n}
	For $({\bf d}_1,{\bf v}_1)$ and $({\bf d}_2,{\bf v}_2)\in \Theta_{CH_{n}}=\Theta_{OH_{n}}$, $\theta_{CH_n}({\bf d}_1,{\bf v}_1)\sim_{sys}\theta_{CH_n}({\bf d}_2,{\bf v}_2)$ if and only if $\theta_{OH_n}({\bf d}_1,{\bf v}_1)\sim_{sys}\theta_{OH_n}({\bf d}_2,{\bf v}_2)$
\end{lemma}
\begin{proof}
	The proof is basically a restatement of the fact that $g_1^{ctr}({\bf d})=g_1^{obs}({\bf d})^T$ and $g_2^{ctr}({\bf d},{\bf v})=g_2^{obs}({\bf d},{\bf v})^T$.
\end{proof}

\begin{definition}\label{sys_well_defined}
We denote by $\sim_{sys}$ the equivalence relation of system isomorphism on $CH_n$ and $OH_n$. We shall also denote $({\bf d}_1,{\bf v}_1)\sim_{sys}({\bf d}_2,{\bf v}_2)$ if $\theta_{CH_n}({\bf d}_1,{\bf v}_1)\sim_{sys}\theta_{CH_n}({\bf d}_2,{\bf v}_2)$ for $({\bf d}_1,{\bf v}_1), ({\bf d}_2,{\bf v}_2)\in \Theta_{CH_{n}}$. The choice of the symbol $\sim_{sys}$ is justified because system isomorphisms for controllable/observable Hamiltonian representations are indeed equivalent as we showed in Lemma \ref{sys_on_CHn_and_OH_n},
\end{definition}

\begin{theorem}
	\label{equivalence_conditon_on_D,v}
	Given $({\bf d}_1,{\bf v}_1)$ and $({\bf d}_2,{\bf v}_2)$ in $\Theta_{CH_{n}}$, then 
	
	\begin{description}
		\item[(I)] $\theta_{CH_n}({\bf d}_1,{\bf v}_1)\sim_{sys}\theta_{CH_n}({\bf d}_2,{\bf v}_2)$ if and only if $a_i({\bf d}_1)=a_i({\bf d}_2)$ and $c_j({\bf d}_1,{\bf v}_1)=c_j({\bf d}_2,{\bf v}_2)$ for all $i,j$. In other words, there exists a permutation matrix $P_{\sigma}\in\mathbb{M}_{n}$ such that, for $D={\rm diag}({\bf d})$ and $P=\begin{bmatrix}
			P_{\sigma}&0\\0&P_{\sigma}
		\end{bmatrix}$, the following conditions hold true:
		\begin{description}
			\item [(i)] $P\begin{bmatrix}
				D_{1}&0\\0&D_{1}
			\end{bmatrix}P^{T}=\begin{bmatrix}
				D_{2}&0\\0&D_{2}
			\end{bmatrix}$
			\item [(ii)] ${\bf v}_1^{T}\begin{bmatrix}
				(F_1)_k&0\\
				0&(F_1)_k
			\end{bmatrix} {\bf v}_1={\bf v}_2^{T}\begin{bmatrix}
				(F_2)_k&0\\
				0&(F_2)_k
			\end{bmatrix} {\bf v}_2$, $k=0,\dots,n-1$
		\end{description}
		The matrices $F_i$ are defined in Theorem \ref{main_theorem}.
		\item[(II)] $\theta_{CH_n}({\bf d}_1,{\bf v}_1) \sim_{filter} \theta_{CH_n}({\bf d}_2,{\bf v}_2)$ under the zero initial state assumption if and only if $e_i({\bf d}_1,{\bf v}_1)=e_i({\bf d}_2,{\bf v}_2)$ for all $i=1,\dots,n$, where the scalar functions $e_i$ are defined recursively as
		\begin{align}\label{condition_for_same_filter}
			\begin{split}
				e_1 &= c_1\\
				e_2 &= c_3-a_{2n-2}\cdot e_1\\
				e_3&=c_5-a_{2n-2}\cdot e_2-a_{2n-4}\cdot e_1\\
				\vdots\\
				e_n &= c_{2n-1}-a_{2n-2}\cdot e_{n-1}-a_{2n-4}\cdot e_{n-2}-\cdots-a_2\cdot e_1
			\end{split}
		\end{align}
	\end{description}
	
\end{theorem}

\begin{remark}
	\normalfont
	In light of the above two propositions, it is clear that two controllable Hamiltonian representations being system isomorphic is a {\it strictly stronger} requirement than them inducing the same filter, since if $a_i({\bf d}_1)=a_i({\bf d}_2)$ and $c_j({\bf d}_1,{\bf v}_1)=c_j({\bf d}_2,{\bf v}_2)$ for all $i,j$, then according to (\ref{condition_for_same_filter}), $e_i({\bf d}_1,{\bf v}_1)=e_i({\bf d}_2,{\bf v}_2)$ trivially holds for all $i\in \mathbb{N}$. 
\end{remark}

\subsection{Equivalence classes of port-Hamiltonian systems by system isomorphisms}
\label{sub1}

We have seen that $PH_n/\sim_{sys}$ is not the set of port-Hamiltonian filters due to the presence of non-canonical systems. However, it is still informative to study the quotient space $PH_n/\sim_{sys}$ because it removes all redundancies in the {\it canonical} part, which prove to be crucial when shall we later restrict our attention to canonical systems in Sections \ref{sub5} and \ref{sub6}. 

In this section, we introduce a manageable characterization of the quotient space $PH_n/\sim_{sys}$ by using parameter spaces. First, motivated by Williamson's theorem, we consider the space $\Theta_{CH_n}$ defined before as the set of all pairs of the form $({\bf d},{\bf v})$, where ${\bf d}=(d_1,d_2,\dots,d_n)^{T}$ with $d_i> 0$, and ${\bf v}=(v_1,v_2,\dots,v_{2n})^T\in \mathbb{R}^{2n}$. Inspired by the representation results, we now define an equivalence relation $\sim_\star$ on $\Theta_{CH_n}$ as below whose equivalence classes are denoted by $[{\bf d},{\bf v}]$. The importance of the next definition is that, as we shall prove in Theorem \ref{characterization}, the relation $\sim_{\star}$ on $\Theta_{CH_n}$ plays the same role as $\sim_{sys}$ on $PH_n$.

\begin{definition}\label{equivalence_relation}
	The pairs $({\bf d}_1,{\bf v}_1)$ and $({\bf d}_2,{\bf v}_2)$ in $\Theta_{CH_n}$ are $\sim_{\star} $-equivalent, that is,  $({\bf d}_1,{\bf v}_1)\sim_{\star}({\bf d}_2,{\bf v}_2)$, if there exists a permutation matrix $P_{\sigma}\in\mathbb{M}_{n}$ and an invertible matrix $A$ such that, for $D_i={\rm diag}({\bf d}_i)$, $i \in \left\{1,2\right\} $ and $P=\begin{bmatrix}
		P_{\sigma}&0\\0&P_{\sigma}
	\end{bmatrix}$, the following conditions hold true:
	\begin{description}
		\item [(i)] $P\begin{bmatrix}
			D_{1}&0\\0&D_{1}
		\end{bmatrix}P^{T}=\begin{bmatrix}
			D_{2}&0\\0&D_{2}
		\end{bmatrix}$
		\item[(ii)] $A^{T}\begin{bmatrix}
			D_{1}&0\\0&D_{1}
		\end{bmatrix}A{\bf v}_1=\begin{bmatrix}
			D_{1}&0\\0&D_{1}
		\end{bmatrix}{\bf v}_1$
		\item[(iii)] $A\mathbb{J}\begin{bmatrix}
			D_{1}&0\\0&D_{1}
		\end{bmatrix}=\mathbb{J}\begin{bmatrix}
			D_{1}&0\\0&D_{1}
		\end{bmatrix}A$
	\item[(iv)] ${\bf v}_2=PA{\bf v}_1$.
	\end{description}
\end{definition}

\begin{proposition}
	The relation $\sim_\star$ defined in Definition \ref{equivalence_relation} is an equivalence relation on $\Theta_{CH_n}$.
\end{proposition}

In the next subsection, we shall give meaning to $\sim_\star$ in terms of groupoid orbits. Now, we aim to characterize the $\sim_{sys}$ equivalence relation on $PH_n$ as the $\sim_{\star}$ equivalence relation on the space $\Theta_{CH_n}$ of $({\bf d},{\bf v})$-pairs, that is, we shall prove that $\Theta_{CH_n}/\sim_{\star}\cong PH_n/\sim_{sys}$. This will be proved in three steps. First, we show that for an arbitrary $S\in Sp(2n,\mathbb{R})$, the map $\varphi_S$ defined in Theorem \ref{main_theorem} composed with $\theta_{CH _n}$ is compatible with the equivalence relations $\sim_{\star}$ and $\sim_{sys}$, that is, $({\bf d}_1,{\bf v}_1)\sim_{\star}({\bf d}_2,{\bf v}_2)$ if and only if $\varphi_S(\theta_{CH_n}({\bf d}_1,{\bf v}_1))\sim_{sys}\varphi_S(\theta_{CH_n}({\bf d}_2,{\bf v}_2))$. Then, we show that the unique map $\psi_S$ induced by $\varphi_S\circ \theta_{CH _n}$ on the quotient spaces does not depend on the choice of $S$ and hence the family of maps $\psi_S$ parameterized by $S\in Sp(2n,\mathbb{R})$ induces a unique map  $\Phi:\Theta_{CH_n}/\sim_{\star}\rightarrow PH_n/\sim_{sys}$ which is a homeomorphism. 

\begin{theorem}[{\bf Characterization of $ {PH_n/\sim_{sys}}$ as $ {\Theta_{CH_n}/\sim_{\star}}$}]
\label{characterization}
	Given any arbitrary $S\in Sp(2n,\mathbb{R})$, the map $\varphi_S\circ \theta_{CH _n}$ induces on the quotient spaces a map $\Phi:\Theta_{CH_n}/\sim_{\star}\rightarrow PH_n/\sim_{sys}$ which does not depend on $S\in Sp(2n,\mathbb{R})$ and is given by $\Phi([{\bf d},{\bf v}]_{\star})=\left[\begin{bmatrix}
	D&0\\
	0&D
\end{bmatrix},{\bf v}\right]_{sys}$, where $D={\rm diag}({\bf d})$. Moreover, $\Phi$ is a  homeomorphism with respect to the quotient topologies.
\end{theorem}

\subsection{The quotient spaces as groupoid orbit spaces}
\label{sub2}

 Recall that from a category theory point of view, a group can be seen as a category with a single object where all morphisms are invertible. Groupoids are a natural generalization of this notion and refer to categories with possibly more than one object, where again all morphisms are invertible (see \cite{mackenzie2005general} for a comprehensive introduction). As it is customary, groupoids will be denoted with the symbol $s,t: {\cal G} \rightrightarrows M $ (or simply ${\cal G} \rightrightarrows M $), where $s$ and $t $ are the {\it source} and the {\it target} maps, respectively. Given $m \in M $, the {\it groupoid orbit} that contains this point is given by $\mathcal{O} _m= t \left(s ^{-1}(m)\right) \subset M$. The {\it orbit space} associated to ${\cal G} \rightrightarrows M $ is denoted by $M/ {\cal G} $.

In this section, we provide an alternative point of view for Theorem \ref{characterization} in terms of groupoid orbits. More precisely, we show first that the set of equivalence classes $PH_{n}/\sim_{sys}$ (resp. $\Theta_{CH_n}/\sim_{\star}$) is the orbit space $PH_n/\mathcal{G}_n$ (resp. $\Theta_{CH_n}/\mathcal{H}_n$) of a groupoid $\mathcal{G}_n\rightrightarrows PH_n$ (resp. $\mathcal{H}_n\rightrightarrows\Theta_{CH_n}$) which we construct in the following paragraphs. In a second step we show that the statement in Theorem \ref{characterization} is equivalent to saying that the orbit spaces $PH_{n}/\sim_{sys}$ and $\Theta_{CH_n}/\mathcal{H}_n$ of the two groupoids coincide.

\begin{definition}\label{groupoid}~
	\begin{enumerate}
		\item Let $\mathcal{G}_n:=\{\left(L,(Q,B)\right)|L\in GL(2n,\mathbb{R}), (Q,B)\in PH_n \text{ such that (i) } \mathbb{J}^TL\mathbb{J}QL^{-1}\text{ is symmetric}$ $ \text{positive-definite (ii) } B=\mathbb{J}^TL^T\mathbb{J}LB\}$.
		\item Let the target and source maps $\alpha,\beta:\mathcal{G}_n\rightarrow {PH}_n$ be defined as $\alpha(L,(Q,B))$ $:=(\mathbb{J}^TL\mathbb{J}QL^{-1},LB)$ and $\beta(L,(Q,B)):=(Q, B)$.
		\item Define the set of composable pairs as $\mathcal{G}^{(2)}_n:=\{((L_1,(Q_1,B_1)),(L_2,(Q_2,B_2)))~|~\beta((L_1,(Q_1,B_1)))$ $=\alpha((L_2,(Q_2,B_2)))\}$.
		\item Let the multiplication map $m:\mathcal{G}^{(2)}_n\rightarrow \mathcal{G}_n$ be defined as $m((L_1,(Q_1,B_1)),(L_2,(Q_2,B_2)))=(L_1L_2,(Q_2,B_2))$.
		\item Let the identity section $\epsilon:{PH}_n\rightarrow \mathcal{G}_n$ be defined as $\epsilon(Q,B):=(\mathbb{I}_{2n},(Q,B))$.
		\item Let the inversion map $i:\mathcal{G}_n\rightarrow \mathcal{G}_n$ be defined as $i(L,(Q,B)):=(L^{-1},(\mathbb{J}^TL\mathbb{J}QL^{-1},LB))$.
	\end{enumerate}
\end{definition}

\begin{proposition}\label{groupoid_orbit}
	The definition above determines a Lie groupoid $\mathcal{G}_n\rightrightarrows PH_n$ with $\mathcal{G}_n$ the total space, ${PH}_n$ the base space, and structure maps $\alpha,\beta,m,\epsilon,i$. We refer to $\mathcal{G}_n\rightrightarrows PH_n$ as the port-Hamiltonian groupoid. The orbit space of this groupoid $PH_n/\mathcal{G}_n$ coincides with $PH_n/\sim_{sys}$.
\end{proposition}

\begin{definition}\label{groupoid}~
	\begin{enumerate}
		\item Let $\mathcal{H}_n:=\big\{\left((P_{\sigma},A),({\bf d},{\bf v})\right)|P_{\sigma}\in \mathbb{M}_n \text{ is a permutation matrix}, A\in GL(2n,\mathbb{R}), ({\bf d},{\bf v})\in \Theta_{CH_n}, \text{ such that (i) } A^{T}\begin{bmatrix}
	    	D&0\\0&D
	    \end{bmatrix}A{\bf v}=\begin{bmatrix}
	    	D&0\\0&D
	    \end{bmatrix}{\bf v} \text{ and (ii) }
       A\mathbb{J}\begin{bmatrix}
	    	D&0\\0&D
	    \end{bmatrix}=\mathbb{J}\begin{bmatrix}
	    	D&0\\0&D
	    \end{bmatrix}A,
        \text{ where } D={\rm diag}({\bf d})\big\}.$
		\item Let the target and source maps $\alpha,\beta:\mathcal{H}_n\rightarrow \Theta_{CH_n}$ be defined as $\alpha((P_{\sigma},A),({\bf d},{\bf v}))$ $:=({\bf d},{\bf v})$ and $\beta((P_{\sigma},A),({\bf d},{\bf v})):=(P_{\sigma}{\bf d}, PA{\bf v})$, where $P=\begin{bmatrix}
			P_{\sigma}&0\\0&P_{\sigma}
		\end{bmatrix}$.
		\item Define the set of composable pairs as $\mathcal{H}^{(2)}_n:=\big\{\big(((P_{\sigma,1},A_1),({\bf d}_1,{\bf v}_1)),((P_{\sigma,2},A_2),({\bf d}_2,{\bf v}_2))\big)~|\\\beta((P_{\sigma,2},A_2),({\bf d}_2,{\bf v}_2))$ $=\alpha((P_{\sigma,1},A_1),({\bf d}_1,{\bf v}_1))\big\}$.
		\item Let the multiplication map $m:\mathcal{H}^{(2)}_n\rightarrow \mathcal{H}_n$ be defined as $m\big(((P_{\sigma,1},A_1),({\bf d}_1,{\bf v}_1)),((P_{\sigma,2},A_2),({\bf d}_2,{\bf v}_2))\big)\\=((P_{\sigma,2}P_{\sigma,1}, P^T_{\sigma,1}A_2P_{\sigma,1}A_1),({\bf d}_1,{\bf v}_1))$.
		\item Let the identity section $\epsilon:\Theta_{CH_n}\rightarrow \mathcal{H}_n$ be defined as $\epsilon({\bf d},{\bf v}):=((\mathbb{I}_{n},\mathbb{I}_{2n}),({\bf d},{\bf v}))$.
		\item Let the inversion map $i:\mathcal{H}_n\rightarrow \mathcal{H}_n$ be defined as $i((P_{\sigma},A),({\bf d},{\bf v})):=((P^T_{\sigma},P_{\sigma}A^{-1}P^T_{\sigma}),(P_{\sigma}{\bf d},PA{\bf v}))$.
	\end{enumerate}
\end{definition}

\begin{proposition}
\label{diagonalized port-Hamiltonian groupoid}
	The definition above determines a Lie groupoid $\mathcal{H}_n\rightrightarrows\Theta_{CH_n}$ with $\mathcal{H}_n$ the total space, $\Theta_{CH_n}$ the base space, and structure maps $\alpha,\beta,m,\epsilon,i$. We refer to $\mathcal{H}_n\rightrightarrows\Theta_{CH_n}$ as the reduced port-Hamiltonian groupoid. The orbit space of this groupoid $\Theta_{CH_n}/\mathcal{H}_n$ coincides with $\Theta_{CH_n}/\sim_{\star}$.
\end{proposition}

Theorem \ref{characterization} can now be restated in terms of the elements that we just introduced.
\begin{theorem}
	The orbit spaces of the Lie groupoids $\mathcal{G}_n\rightrightarrows {PH}_n$ and $\mathcal{H}_n\rightrightarrows\Theta_{CH_n}$ are isomorphic. 
\end{theorem}

\subsection{Characterization of canonical port-Hamiltonian systems}
\label{sub3}

In Subsections \ref{sub1} and \ref{sub2} we have provided a characterization of ${PH}_n/\sim_{sys}$ in terms of $\Theta_{CH_n}/\sim_{\star}$ and groupoid orbit spaces. Recall from Subsection \ref{sub0} that the difficulty of the unique identifiability of filters in $\mathcal{PH}_n$ comes from the possible presence of non-canonical systems, without which the equivalence relation $\sim_{filter}$  coincides with $\sim_{sys}$. Hence, it is worth studying what the quotient spaces above look like when restricted to the subset that contains only canonical port-Hamiltonian systems. In this section, we take a step in that direction. 

Recall that a port-Hamiltonian system in ${PH}_n $ of the form (\ref{ph_definition}) is controllable (or equivalently, observable/canonical) if and only if \begin{equation}\label{controllable_condition1}
	\det\left(\begin{bmatrix}
		B~|~\mathbb{J}QB~|~\dots~|~(\mathbb{J}Q)^{2n-1}B
	\end{bmatrix}\right)\neq0.
\end{equation} 
Using the Williamson decomposition of $Q$ into $D$ and $S$, and $v:=S\cdot B$, this is equivalent to \begin{equation}\label{controllable_condition2}
\det\left(\begin{bmatrix}
	{\bf v}~\bigg|~\mathbb{J}\begin{bmatrix}
		D&0\\0&D
	\end{bmatrix}{\bf v}~\bigg|~\dots~\bigg|~\left(\mathbb{J}\begin{bmatrix}
	D&0\\0&D
\end{bmatrix}\right)^{2n-1}{\bf v}
\end{bmatrix}\right)\neq0.
\end{equation}

Denote by ${PH}^{can}_n$ (respectively, $\Theta_{CH_n}^{can}$) the subset of ${PH}_n$ (respectively, $\Theta_{CH_n}$) made of systems that satisfy (\ref{controllable_condition1}) (respectively,  (\ref{controllable_condition2})). As an immediate consequence, it holds true that 
$$\mathcal{PH}^{can}_n\cong  PH^{can}_n/\sim_{filter}\cong PH^{can}_n/\sim_{sys}.$$ We now characterize the space of pairs $({\bf d},{\bf v}) \in \Theta_{CH_n}$  that correspond to canonical port-Hamiltonian systems in normal form. The calculation of the determinant in (\ref{controllable_condition2}) yields $\big(\prod_{i=1}^{n}d_i\big)\cdot\big(\prod_{1\leq j<k\leq n}^{}(d_j+d_k)^2(d_j-d_k)^2\big)\cdot\big(\prod_{l=1}^{n}(v_l^2+v_{n+l}^2)\big)$ up to a sign. Therefore, 
\begin{equation*}
\Theta_{CH_n}^{can}= \left\{({\bf d},{\bf v}) \in \Theta_{CH_n}\mid \ \mbox{the entries of {\bf d} are all different and $v_l^2+v_{n+l}^2> 0$ for all $l\in \left\{1, \ldots, n\right\}$} \right\}.
\end{equation*}
We shall refer to the statement on the entries of {\bf d} being all different as the {\it non-resonance condition} and to $v_l^2+v_{n+l}^2> 0$ for all $l\in \left\{1, \ldots, n\right\}$ as the {\it nondegeneracy} condition.
There might be a concern about whether different choices of the matrix $S$ lead to different vectors ${\bf v}$ and hence the notion of nondegeneracy would be ill-defined. This is indeed not a problem since, as we show in Remark \ref{torus_motivation} below, once the non-resonance condition is assumed, different vectors ${\bf v}$ are obtained by rotating the planes spanned by each and every pair of $l$-th and $n+l$-th entries, which preserves the value of $v_l^2+v_{n+l}^2$. Thus, the nondegeneracy condition is actually based on the non-resonance condition.

\medskip

\begin{remark}[\bf Williamson's decomposition in the canonical case]\label{torus_motivation}
		\normalfont
	We have mentioned in Theorem \ref{Williamson's normal form} {\bf (iii)} that two symplectic matrices $S$ and $S^{\prime}$ that Williamson decompose the same $Q$ differ by a unitary matrix. We now note that for an element $Q$ that satisfies the non-resonance condition, $S$ and $S^{\prime}$ do not only differ by an arbitrary $U\in U(n)$, (see \eqref{monomorp un} for the definition of $U(n)$) but by a special one $R$ that has the form 
	\begin{equation}
		\label{torus}
		R = \left[
		\begin{array}{ccc|ccc}
			\cos\theta_1&&\bigzero&-\sin\theta_1&&\bigzero\\
			&\ddots&&&\ddots&\\
			\bigzero&&\cos\theta_n&	\bigzero&&-\sin\theta_n\\
			\hline
			\sin\theta_1&&\bigzero&\cos\theta_1&&\bigzero\\
			&\ddots&&&\ddots&\\
			\bigzero&&\sin\theta_n&	\bigzero&&\cos\theta_n\\
		\end{array}
		\right].\end{equation} 
	This fact accounts for part of a symmetry that we shall spell out later on. The proof of this fact is purely computational: the assumption that the diagonal entries of $D$ are all positive and distinct, the fact that $U$ satisfies the equation $U\begin{bmatrix}
		D&0\\0&D
	\end{bmatrix}U^T=\begin{bmatrix}
	D&0\\0&D
\end{bmatrix}$ and, at the same time, $U\in U(n)=SO(2n,\mathbb{R})\cap Sp(2n,\mathbb{R})$, guarantees the claim.
\end{remark}

\medskip

\begin{remark}[\bf Being canonical is a generic property]
	\normalfont
	It is well-known that the set of canonical systems, as a subset of all linear systems, corresponds to a Zariski open set, which is open and dense in the usual topology \cite{Krzysztof1983}. In particular, this also holds for linear port-Hamiltonian systems. Therefore, ${PH}^{can}_n$ is open and dense in ${PH}_n$. On the other hand, using the characterization provided above, it is clear that $ \Theta_{CH_n}^{can}$ is also open and dense in $\Theta_{CH_n}$.
\end{remark}

The isomorphism in Theorem \ref{characterization} naturally restricts to canonical subsets, that is $ {PH}^{can}_{n}/\sim_{sys}{\cong}~\Theta_{CH_n}^{can}/\sim_\star$. On the other hand, we will see below another isomorphism result involving ${PH}^{can}_{n}/\sim_{sys}$.

\begin{proposition}[{\bf Characterization of $ {{PH}^{can}_{n}/\sim_{sys}}$ as $ {\Theta_{CH_n}^{can}/\sim_{sys}}$}]
\label{characterization2}
	The map $\Phi:\Theta^{can}_{CH_n}/\sim_{sys}\rightarrow PH^{can}_n/\sim_{sys}$ defined by $\Phi([{\bf d},{\bf v}]_{sys})=\left[\begin{bmatrix}
		D&0\\
		0&D
	\end{bmatrix},{\bf v}\right]_{sys}$, where $D={\rm diag}({\bf d})$, is an isomorphism.
\end{proposition}

We just proved that both $\Theta_{CH_n}^{can}/\sim_\star$ and $\Theta_{CH_n}^{can}/\sim_{sys}$ are isomorphic to ${PH}^{can}_{n}/\sim_{sys}$, and even via the same ismorphism $\Phi$. Therefore, the equivalence relations ${\sim_{\star}}$ and ${\sim_{sys}}$ {\it  coincide} when restricted to $\Theta_{CH_n}^{can}$. 

To summarize, we have proved in this subsection that
\begin{equation*}
	{\mathcal{PH}^{can}_n}\cong{PH}^{can}_n/\sim_{sys}~\cong~\Theta_{CH_n}^{can}/\sim_{\star}~\cong~\Theta_{CH_n}^{can}/\sim_{sys}.
\end{equation*}
In the next subsection, we  continue the investigation of the above chain of isomorphisms.

\subsection{The unique identifiability space for canonical port-Hamiltonian systems as a group orbit space}
\label{sub5}

In Subsection \ref{sub2}, it is proved that the quotient space $PH_n/\sim_{sys}$ can be treated as a Lie groupoid orbit space. We now show that the restricted quotient space to canonical port-Hamiltonian systems, that is, $PH^{can}_n/\sim_{sys}$, is isomorphic to the orbit space of a certain group action on $\Theta^{can}_{CH_{n}}$, where the group is a semi-direct product of the $n$-permutation group and the $n$-torus, that is, $S_n\rtimes_{\phi}\mathbb{T}^n$. The intuition behind this fact is that restricting to the subset of canonical systems $PH^{can}_n$ removes the degeneracies in $PH_n$, which allows to reduce the symmetry of the Lie groupoid $\mathcal{G}_n\rightrightarrows PH_n$ to that of the Lie group $S_n\rtimes_{\phi}\mathbb{T}^n$.

We start by defining the group action. First, let the permutation group $S_n$ act on ${\Bbb R}^n $ by permuting the entries $d_i$ of the vector ${\bf d} \in {\Bbb R}^n$. For each $i \in  \left\{1, \ldots, n\right\}$ the circle $S^1$ acts on the plane spanned by the $i$-th and $(n+i)$-th entries of ${\bf v}$ by rotations. More precisely, we define the action of $S_n$ on elements ${\bf d}$ and ${\bf v}$ as \begin{equation*}
	\begin{aligned}
		&\Gamma_{\sigma}\big((d_1,\dots,d_n)^T\big)=(d_{\sigma(1)},\dots,d_{\sigma(n)})^T=P_{\sigma}\cdot(d_1,\dots,d_n)^T\,
	\end{aligned}
\end{equation*} 
where $P_\sigma $ is the corresponding permutation matrix and 
\begin{equation*}
	\begin{aligned}
		&\Gamma_{\sigma}\big((v_1,\dots,v_{2n})^T\big)=(v_{\sigma(1)},\dots,v_{\sigma(n)},v_{n+\sigma(1)},\dots,v_{n+\sigma(n)})^T=\begin{bmatrix}
			P_\sigma&0\\
			0&P_\sigma
		\end{bmatrix}\cdot (v_1,\dots,v_{2n})^T,
	\end{aligned}
\end{equation*} respectively. Then the $\sigma$-action on a pair $({\bf d},{\bf v})$ is understood as acting on ${\bf d}$ and ${\bf v}$ simultaneously. We also define the action of the $i$-th circle of the torus $\mathbb{T}^n$  as the planar rotation of the space spanned by the $i$-th and $(n+i)$-th entries of ${\bf v}$. This torus action is understood to leave ${\bf d}$ invariant. More concretely, it is the action
\begin{equation*}\begin{aligned}
		&\Gamma_{\theta_i}\big((d_1,\dots,d_n,v_1,\dots,v_{2n})^T\big)\\
		&=(d_1,\dots,d_n,v_1,\dots,v_{i-1},cos\theta_iv_i-sin\theta_iv_{n+i},v_{i+1},\dots,v_{n},\\&~~~~~~~~~~~~~~~~~~v_{n+1},\dots,v_{n+i-1},sin\theta_iv_{i}+cos\theta_iv_{n+i},v_{n+i+1},\dots,v_{2n})^T.
	\end{aligned}
\end{equation*}
With these actions of the groups $S_n$ and $\mathbb{T}^n$ on $\Theta_{CH_{n}}$ we define the map $\Gamma_{(\sigma,(\theta_1,\dots,\theta_n)^T)}:(\mathbb{R}^n_{+}\times\mathbb{R}^{2n})\rightarrow(\mathbb{R}^n_{+}\times\mathbb{R}^{2n})$
as
\begin{multline}
\label{semi_direct_action}
		\Gamma_{(\sigma,(\theta_1,\dots,\theta_n)^T)}({\bf d},{\bf v})
		=\Gamma_{\theta_1}\circ\dots\circ\Gamma_{\theta_n}\circ\Gamma_{\sigma}({\bf d},{\bf v})\\
		=(P_\sigma\cdot {\bf d},\Gamma_{\theta_1}\circ\dots\circ\Gamma_{\theta_n}\bigg(\begin{bmatrix}
			P_\sigma&0\\
			0&P_\sigma
		\end{bmatrix}\cdot {\bf v}\bigg))=(P_\sigma \cdot {\bf d},RP\cdot {\bf v}),
\end{multline}
which constitutes an action of the semi-direct product group $S_n\rtimes_{\phi}\mathbb{T}^{n}$, where $\phi:S_n\rightarrow Aut(\mathbb{T}^{n})$ is given by the permutation $\phi(\sigma)((\theta_1,\dots,\theta_n)^T)=P_{\sigma}\cdot(\theta_1,\dots,\theta_n)^{T}$. 
Note that the matrix of $\Gamma_{\theta_1}\circ\dots\circ\Gamma_{\theta_n}$ is given by (\ref{torus}), $P_\sigma$ is the permutation matrix that corresponds to $\sigma\in S_n$, and $P=\begin{bmatrix}
	P_\sigma&0\\
	0&P_\sigma
\end{bmatrix}$.

\begin{proposition}
\label{action_verification}
	The map $\Gamma_{(\sigma,(\theta_1,\dots,\theta_n)^T)}$ defined as (\ref{semi_direct_action}) for $\sigma\in S_n$ and $(\theta_1,\dots,\theta_n)^T\in\mathbb{T}^{n}$ is a left group action of $(S_n\rtimes_{\phi}\mathbb{T}^{n})$ on $\Theta_{CH_n}$.
\end{proposition}

Using the definition of the $(S_n\rtimes_{\phi}\mathbb{T}^{n})$-action on $\Theta_{CH_n}$, two elements $({\bf d}_1,{\bf v}_1),({\bf d}_2,{\bf v}_2)\in \Theta_{CH_n}$ are in the same orbit if and only if the following conditions hold true for some $\sigma\in S_n$:
\begin{description}\label{Dv_equivalence}
	\item[(i)] ${d}_{2,i}= {d}_{1,\sigma(i)},$
	\item[(ii)] ${v}_{2,i}^2+{v}_{2,n+i}^2={v}_{1,\sigma(i)}^2+{v}_{1,n+\sigma(i)}^2,~i=1,\dots,n.$
\end{description}
By Theorem \ref{equivalence_conditon_on_D,v} {\bf (I)}, parts {\bf (i)} and {\bf (ii)} it is clear that there is a close relation between the  the $(S_n\rtimes_{\phi}\mathbb{T}^{n})$-action and the equivalence relation $\sim_{\star}$ on $\Theta_{CH_n}$. The next proposition demonstrates that the orbits of the $(S_n\rtimes_{\phi}\mathbb{T}^{n})$-action coincide with the equivalence classes of the relation $\sim_{sys} $ when we restrict our attention to the subset $\Theta_{CH_n}^{can}$.

\begin{proposition}[{\bf Characterization of $ {\Theta_{CH_n}^{can}/\sim_{sys}}$ as $ {\Theta_{CH_n}^{can}/} {(S_n\rtimes_{\phi}\mathbb{T}^{n})}$}]
\label{characterization3}
		Given $({\bf d}_1,{\bf v}_1)$ and $({\bf d}_2,{\bf v}_2)$ in $\Theta_{CH_n}^{can}$, then $({\bf d}_1,{\bf v}_1)\sim_{sys}({\bf d}_2,{\bf v}_2)$ if and only if $({\bf d}_1,{\bf v}_1)$ and $({\bf d}_2,{\bf v}_2)$ lie in the same orbit of the $(S_n\rtimes_{\phi}\mathbb{T}^{n})$-action.
\end{proposition}

\subsection{Global Euclidean coordinates for the unique identifiability space of canonical port-Hamiltonian systems  }
\label{sub6}

Recall from Section \ref{sub3} that $\Theta_{CH_n}^{can}$ contains pairs $({\bf d},{\bf v})$ where ${\bf d}\in \mathbb{R}_+^{n}$ and ${\bf v}\in\mathbb{R}^{2n}$ are such that the entries $d_l$'s are all distinct and ${v}_{l}^2+{v}_{n+l}^2> 0$ for all $l=1,\dots,n$. We define for convenience a function $\mathcal{R}:\mathbb{R}^{2n}\rightarrow\mathbb{R}_{\geq0}^{n}$ as $\mathcal{R}((v_{1},\dots,v_{2n})^{T})=\big({v}_{1}^2+{v}_{n+1}^2,\dots,{v}_{n}^2+{v}_{2n}^2\big)$.\par Now observe that the quotient space $\Theta_{CH_n}^{can}/(S_n\rtimes_{\phi}\mathbb{T}^{n})$ naturally has a smooth manifold structure. We briefly prove this in the following lines. Note that the torus $\mathbb{T}^{n}$ is a connected abelian compact Lie group. The symmetry group $S_n$ is a finite group, and hence compact as well. Thus, it is easy to see that the semi-direct product $S_n\rtimes_{\phi}\mathbb{T}^{n}$ is also a compact Lie group, and hence its action on $\Theta_{CH_n}^{can}$ is automatically proper. On the other hand, since $\Theta_{CH_n}^{can}$ is the space of $({\bf d},{\bf v})$ pairs satisfying that ${\bf d}$ contains distinct entries and $\mathcal{R}({\bf v})^{(l)}>0$ for $l=1,\dots,n$, it necessarily holds that the only element in $S_n\rtimes_{\phi}\mathbb{T}^{n}$ that possibly keep any element in $\Theta_{CH_n}^{can}$ invariant is the identity, which implies the $(S_n\rtimes_{\phi}\mathbb{T}^{n})$-action on $\Theta_{CH_n}^{can}$ is free. Classical results in Lie theory \cite[Proposition 2.3.8]{Ortega2004} guarantee that $\Theta_{CH_n}^{can}/(S_n\rtimes_{\phi}\mathbb{T}^{n})$ admits a unique smooth structure such that the quotient map $\pi:\Theta_{CH_n}^{can}\rightarrow\Theta_{CH_n}^{can}/(S_n\rtimes_{\phi}\mathbb{T}^{n})$ is a submersion. With this as a motivation, we try to find the quotient space explicitly in the following. \par

For a fixed ${\bf d}$, we denote by ${\bf d}_{\uparrow}$ the reordered vector constructed out of ${\bf d}$ by placing the entries in increasing order. Denote by $\mathcal{X}^{n}_{\uparrow}$ the set of ${\bf d}\in \mathbb{R}_+^n$ with distinct positive entries in increasing order. We have then the following proposition that explicitly characterizes the quotient space $\Theta_{CH_n}^{can}/(S_n\rtimes_{\phi}\mathbb{T}^{n})$.

\begin{proposition}[{\bf Global Euclidean coordinates for the orbit space} ${\Theta_{CH_n}^{can}/}{(S_n\rtimes_{\phi}\mathbb{T}^{n})}$]
\label{characterization4}
	 The map $f:\Theta_{CH_n}^{can}/(S_n\rtimes_{\phi}\mathbb{T}^{n})\rightarrow\mathcal{X}^{n}_{\uparrow}\times \mathbb{R}_{+}^{n}$ defined by $f([{\bf d},{\bf v}])=({\bf d}_{\uparrow},\mathcal{R}(\Gamma_{\sigma}({\bf v})))$, where $\sigma\in S_n$ is the unique permutation such that $\Gamma_{\sigma}({\bf d})={\bf d}_{\uparrow}$, is an isomorphism.
\end{proposition}

\section{Linear port-Hamiltonian systems in normal form are restrictions of higher dimensional ones} 

In this section, we prove a theorem (Theorem \ref{higher dimension}), inspired by the classical Kalman Decomposition \cite{Jacob2012}, which says the filter induced by any $(Q,B)\in{PH}_n$ can be regarded as that induced by some $(Q^\prime,B^\prime)\in{PH}_m$, where $m$ can be any integer that is at least $n$. The motivation for these considerations is given by the fact that in many practical situations in which an input/ouput system has to be learned, the dimension of the underlying state-space system is not known. In that situation, we may want to have the flexibility of considering the actual system that needs to be learned as a lower-dimensional restriction of a much larger-dimensional one that we have picked for the learning task.

We shall carry this out by producing an explicit injective system morphism between the state space of $(Q,B)$ and that of $(Q^{\prime},B^{\prime})$ in our next Theorem \ref{higher dimension}. In Proposition \ref{low_high_isomorphism}, we show that the quotient space ${PH}_n/\sim_{sys}$ can be characterized as $PH_{m,n}/\sim_{sys}$, where $PH_{m,n}\subset PH_m$ is the space containing all the systems of the form $(Q^{\prime},B^{\prime})$. Motivated by the developments in Section \ref{quotient}, we then characterize the pair $({\bf d}^\prime,{\bf v}^{\prime})$ that corresponds to $(Q^\prime,B^\prime)$ in Proposition \ref{Dv_correspondence}. Eventually, in Proposition \ref{Dv_QB_correspondence}, we show that the isomorphism $PH_n/\sim_{sys}~\cong~\Theta_{CH_n}/\sim_\star$ can be lifted to high dimension as well. We shall comment further at the end of this section on the significance of the above-mentioned results in the context of machine learning.\par
\medskip
The following theorem states that the filter induced by $(Q,B)\in PH_n$ can be reproduced using systems in an arbitrarily higher dimension.
\begin{theorem}
	\label{higher dimension}
	Given any system $(Q,B)\in PH_n$, then 
    \begin{description}
    	\item[(i)] For any $m\geq n$, there exists an orthogonal matrix $O\in O(2m,\mathbb{R})$ such that the filter induced by $(Q^{\prime},B^{\prime})=\bigg(O\begin{bmatrix}
    		Q&0\\
    		0&\mathbb{I}_{2m-2n}
    	\end{bmatrix}O^{T},O\begin{bmatrix}
    		B\\
    		0
    	\end{bmatrix}\bigg)\in{PH}_m$ coincides with that induced by $(Q,B)$.
    	\item[(ii)] The map $f:\mathbb{R}^{2n}\rightarrow \mathbb{R}^{2m}$ defined by $f({\bf z})=O\begin{bmatrix}
    		\mathbb{I}_{2n}\\0
    	\end{bmatrix}\cdot {\bf z}$ is an injective system morphism between the state spaces of $(Q,B)$ and $(Q^{\prime},B^{\prime})$.
    \end{description} 	
\end{theorem}

As it can be seen in the proof (included in Appendix \ref{Proof of Theorem higher dimension}), the matrix $O\in O(2m,\mathbb{R})$ above is constructed so that 
\begin{equation}
\label{condition o matrix}
O\begin{bmatrix}
	\mathbb{J}_{n}&0\\
	0&\mathbb{J}_{m-n}
\end{bmatrix}O^{T}=\mathbb{J}_m. 
\end{equation}
From now on, we denote by ${PH}_{m,n}\subset{PH}_m$ the space of linear port-Hamiltonian systems parametrized by pairs $(Q^{\prime},B^{\prime})$ of the form $\bigg(O\begin{bmatrix}
	Q&0\\
	0&\mathbb{I}_{2m-2n}
\end{bmatrix}O^{T},O\begin{bmatrix}
	B\\
	0
\end{bmatrix}\bigg),$ where $O\in O(2m,\mathbb{R})$ satisfies \eqref{condition o matrix},
and equip it with the system automorphism relation $\sim_{sys}$ defined on $PH_m$. The following proposition states that the space of possible input-output filters induced by ${PH}_n$ is indeed the same as those induced by ${PH}_{m,n}$. This means we can exactly reproduce the filters of $2n$-dimensional port-Hamiltonian systems in higher dimension by simply considering the elements $(Q^{\prime},B^\prime)$ in ${PH}_{m,n}$.
\begin{proposition}
\label{low_high_isomorphism}
	The  function $f:{PH}_{n}/\sim_{sys}\rightarrow{PH}_{m,n}/\sim_{sys}$ defined by \begin{equation*}
		f([Q,B]_{sys})=\left[O\begin{bmatrix}
			Q&0\\
			0&\mathbb{I}_{2m-2n}
		\end{bmatrix}O^{T},O\begin{bmatrix}
			B\\
			0
		\end{bmatrix}\right]_{sys}
	\end{equation*}is an isomorphism, where $O\in O(2m,\mathbb{R})$ is as in Theorem \ref{higher dimension} and hence satisfies \eqref{condition o matrix}.
\end{proposition}

Recall that for a system $(Q,B)\in{PH}_n$, we derive the corresponding object $({\bf d},{\bf v})\in\Theta_{CH_{n}}$ from Williamson's decomposition ${\displaystyle Q=S^{T}\begin{bmatrix}
		D&0\\
		0&D
\end{bmatrix}}S$ and ${\bf v}=S^{-1}B$. We have seen that $(Q^{\prime},B^{\prime})\in{PH}_{m,n}\subset {PH}_m$ is also a linear port-Hamiltonian system in normal form. Therefore, it makes sense to investigate the relation between $({\bf d},{\bf v})$ and the element $({\bf d}^{\prime},{\bf v}^{\prime})$ which corresponds to $(Q^\prime,B^\prime)$. The following proposition asserts that ${\bf d}^{\prime}$ can be obtained from ${\bf d}$ by padding it with ones and, similarly, ${\bf v}^\prime$ can be obtained by splitting ${\bf v}$ and padding each segment with zeros.
\begin{proposition}[{\bf Symplectic eigenvalues for corresponding higher dimensional systems}]
\label{Dv_correspondence}
	Let $(Q,B)$ and $(Q^{\prime},B^{\prime})$ be as in Theorem \ref{higher dimension}, and let ${\bf d}$ and ${\bf d}^{\prime}$ be their corresponding symplectic eigenvalues. Then, up to reordering, ${\bf d}^{\prime}$= $(d_1,\cdots,d_n,1,1,\dots,1)^T$. Even though ${\bf v}$ and ${\bf v}^{\prime}$ are not uniquely determined (See Remark \ref{non_uniqueness}), there exists a choice of ${\bf v}^{\prime}$ that is related to ${\bf v}=(v_1,\cdots,v_n,v_{n+1},\cdots,v_{2n})^T$ via
	\begin{equation*}
		{\bf v}^{\prime}=\big(
		v_1,\cdots,v_n,\underbrace{0,\cdots,0}_{m-n},v_{n+1},\cdots,v_{2n},\underbrace{0~\cdots~0}_{m-n}\big)^{T}.
	\end{equation*}
\end{proposition}

From the above proposition, we call ${\bf d}^{\prime}$ the extended symplectic eigenvalues and ${\bf v}^\prime$ the extended vector. Now we define the space $\Theta_{CH_{m,n}}$ as the set of all pairs of the form $({\bf d}^\prime,{\bf v}^\prime)$ and equip $\Theta_{CH_{m,n}}$ with the equivalence relation $\sim_{\star}$ as in Definition \ref{equivalence_relation} but in dimension $m$ instead of $n$. Recall that we proved $\Theta_{CH_{n}}/\sim_\star~\cong~{PH}_n/\sim_{sys}$. Now we proceed to show that the above isomorphism in dimension $2n$ can be lifted to dimension $2m$ by considering only the restricted parameter spaces with vectors of the form $({\bf d}^\prime,{\bf v}^\prime)$ and $(Q^\prime,B^\prime)$.

\begin{proposition}
\label{Dv_QB_correspondence}
	The function $f:\Theta_{CH_{m,n}}/\sim_{\star}\longrightarrow{PH}_{m,n}/\sim_{sys}$ defined by 
	\begin{equation*}
		f([{\bf d}^{\prime},{\bf v}^{\prime}]_{\sim_{\star}})=\left[\begin{bmatrix}
			D^{\prime}&0\\0&D^{\prime}
		\end{bmatrix},{\bf v}^{\prime}\right]_{sys}
	\end{equation*} where $D^{\prime}={\rm diag}({\bf d}^{\prime})$, is an isomorphism.
\end{proposition}

Note that in general ${\bf d}^{\prime}$ contains repeated symplectic eigenvalues because of all the ones used in the extension and that $ v_l^{\prime2}+v_{m+l}^{\prime2}=0$ for $l>n$. Therefore, it is impossible that $\Theta_{CH_{m,n}}$ contains canonical systems for $m>n$. In other words, lifting ${PH}_n$ to ${PH}_{m,n}$ {\it introduces degeneracies that exclude the possibility of the systems being canonical.} \par

We emphasize that the above-mentioned series of results are crucial in machine learning applications. Very often in practice, the dimension $2n$ of the underlying data-generating process, that is, the latent port-Hamiltonian system (\ref{ph_definition}), is not known, causing a problem when choosing the dimension of the controllable/observable Hamiltonian representation for learning. This issue can be solved by composing the morphism in Theorem \ref{higher dimension} {\bf (ii)} (which is injective) and the one in Theorem \ref{main_theorem} (not necessarily injective). The composition of system morphisms is still a system morphism, this time between the underlying system $\theta_{PH _n}(Q,B)$ and the observable Hamiltonian representation in an arbitrarily higher dimension $2m\geq 2n$. In this way, the observable Hamiltonian representation in dimension $2m$ still has full expressive power to represent any $2n$-dimensional system in ${PH}_n$, and hence can be used for learning. Practically, one can choose a sufficiently large $m$, and parameterize the observable Hamiltonian representation using $({\bf d},{\bf v})$ (we use the notation $({\bf d},{\bf v})$ instead of $({\bf d}^\prime,{\bf v}^\prime)$ because practically we do not know what $n$ is) and then estimate them. We emphasize that the higher-dimensional port-Hamiltonian systems are in general not canonical, hence the $({\bf d},{\bf v})$-pair that corresponds to the data-generating process is not guaranteed to be unique. Still, we always know there is at least one choice of $({\bf d},{\bf v})$ that works no matter how large an $m$ we choose, and which is constructed using the recipe in Proposition \ref{Dv_correspondence}.

\section{Practical implementation of the results}
\label{learnability}

We start with a diagram that summarizes the results that we have proved.

\begin{theorem}
\label{diagram}
	The following diagram holds true using the isomorphisms explicitly constructed in all the preceeding results.
		\[
\begin{tikzcd}
	&{\Theta_{CH_{m,n}}/\sim_{\star}} \arrow{r}{\cong} \arrow{d}{\cong} & {{PH}_{m,n}/\sim_{sys}} \arrow{l}{} \arrow{d}{\cong}\\
	\Theta_{CH_{n}}/\mathcal{H}_n \arrow{r}{\cong} & \Theta_{CH_{n}}/\sim_\star \arrow{u}{} \arrow{r}{\cong} \arrow{l}{}  & {{PH}_{n}/\sim_{sys}} \arrow{u}{} \arrow{l}[swap]{} \arrow{r}{\cong} &  {{PH}_{n}/\mathcal{G}_n} \arrow{l}{} \\
	&\Theta_{CH_{n}}^{can}/\sim_\star \arrow{r}{\cong} \arrow{d}{\cong} \arrow{u}{}& {{PH}^{can}_{n}/\sim_{sys}}  \arrow{l}[swap]{} \arrow{d}{\cong} \arrow{u}{} \arrow{r}[swap]{} & PH^{can}_n/\sim_{filter} \arrow{l}[swap]{\cong} \arrow{r}{} & \mathcal{PH}^{can}_n \arrow{l}[swap]{\cong}\\
	&\Theta^{can}_{CH_{n}}/(S_n\rtimes_{\phi}\mathbb{T}^{n}) \arrow{r}{\cong} \arrow{u}[swap]{} \arrow{d}{\cong}& \Theta^{can}_{CH_{n}}/\sim_{sys}  \arrow{l}[swap]{}  \arrow{u}[swap]{}\\
	&{ \mathcal{X}^{n}_{{\uparrow}}\times \mathbb{R}_{+}^{n}} \arrow{u}[swap]{}\\
\end{tikzcd}
\]
\end{theorem}

We now comment on how to use the results contained in the diagram above depending on the different learning situations that we may encounter. Indeed, we can use our statements to tackle three different learning scenarios:


\medskip

\begin{itemize}
	\item Case 1: The target port-Hamiltonian system (the data generating process that we want to learn) is canonical and its state-space dimension is known, that is,  $\theta_{{PH}_n}(Q,B)\in {PH}^{can}_n$ with $n$ known. This is the most favorable situation in the sense that we can exactly represent the system $\theta_{{PH}_n}(Q,B)$ by either the controllable or the observable Hamiltonian representation, which are both isomorphic to the original system. Furthermore, since in this case $\sim_{filter}$ coincides with $\sim_{sys}$, the filter in $\mathcal{PH}_n^{can}$ induced by $\theta_{{PH}_n}(Q,B)$ can be uniquely identified with an element $({\bf d}_\uparrow,R)\in\mathcal{X}^{n}_{{\uparrow}}\times \mathbb{R}_{+}^{n}$, which are  the unique parameters that need to be estimated.
	\item Case 2: The target port-Hamiltonian system is not guaranteed to be canonical but its dimension is known, that is, $\theta_{{PH}_n}(Q,B)\in {PH}_n$ with $n$ known. In this case, there is a trade-off between the controllable Hamiltonian representation and the observable one. As mentioned before, the controllable one will be structure-preserving but its expressive power depends on the controllability of the target system $\theta_{{PH}_n}(Q,B)$. On the other hand, the observable one always possesses full expressive power but does not always guarantee the port-Hamiltonian structure of the induced filter.
	\item Case 3: We are agnostic about the dimension of the target port-Hamiltonian system, that is, given $\theta_{{PH}_n}(Q,B)\in {PH}_n$ with $n$ unknown. In this case, we need to choose a sufficiently large $m$ so that $m\geq n$, then based on composition of system morphisms, it suffices to learn some $({\bf d},{\bf v})\in\Theta_{CH_{m}}$ and use the $2m$-dimensional observable Hamiltonian representation to reproduce the input-output dynamics of $(Q,B)$. Due to the loss of the canonical property, such a $({\bf d},{\bf v})$ pair may not be unique. Additionally, we do not know the dimension $2n$ of the data generating process we ignore how many ones are used to pad ${\bf d}$ (and similarly, how many zeros are padded into the vector ${\bf v}$). However, we do know that an element $({\bf d},{\bf v})$ exists in some $\Theta_{CH_{m,n}}\subset\Theta_{CH_{m}}$ which is given by Proposition \ref{Dv_correspondence}. 
\end{itemize}

An important special case is when there is no input to the port-Hamiltonian system, that is, $u(t)=0$. In this case, the port-Hamiltonian system reduces to a linear Hamiltonian system with an arbitrary linear readout matrix because $Q$ is positive-definite and hence invertible. We emphasize that the observable Hamiltonian representation in a higher dimension is totally independent of $B$ since it is simply given by 
\begin{equation}\label{autonumous_case}
	\left\{
	\begin{aligned}
		\dot{{\bf s}}&=g^{obs}_1({\bf d})\cdot {\bf s},\\
		y&=\left(
		0 , 0 ,\cdots, 0, 1\right)\cdot {\bf s},
	\end{aligned} \right.
\end{equation} 

 In other words, Hamiltonian systems with linear readout can be learned by adjusting the initial state ${\bf s}_0$ and symplectic eigenvalues $d_i$, without even knowing the linear readout function that yields the observations.

\section{Numerical illustrations}
\label{numerics}

In this section, we present two numerical examples to demonstrate the effectiveness of our representation results from a learning point of view.

\subsection{Non-dissipative circuit}

\begin{figure}[h]\centering
	\begin{circuitikz} 
		\draw
		(0,0) to[battery, l=V, color=red] (5,0) to (5,-1)
		to [C, l=$C_1$, color=blue] (2,-1) to [cute inductor, l = $L_1$,color=green] (0,-1) to (0,0);
		\draw (5,-1) to (5,-2.5) to [C, l=$C_2$, color=blue] (2,-2.5) to [cute inductor, l = $L_2$,color=green] (0,-2.5) to (0,0);
		\draw (5,-2.5) to (5,-4) to  (2,-4);
		\node[draw=none] (ellipsis1) at (1.5,-4) {$\cdots$};
		\draw (1,-4) to (0,-4) to (0,-4);
		\draw (5,-4) to (5,-5) to [C, l=$C_5$,color=blue] (2,-5) to [cute inductor, l = $L_5$,color=green] (0,-5) to (0,0);
	\end{circuitikz}
	\caption{Lossless circuit port-Hamiltonian system}\label{circuit} 
\end{figure}
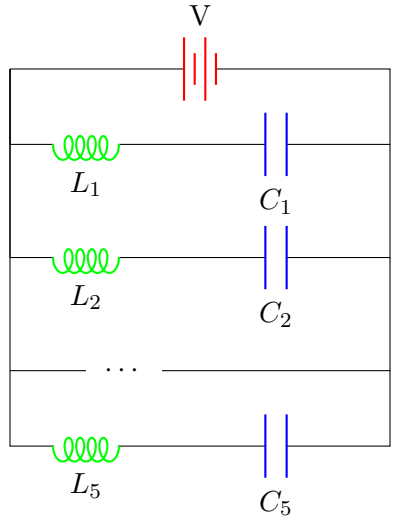

Similar to an example in \cite{Medianu2021}, we consider a circuit consisting of a power source with voltage $V=u(t)$, together with ten parallelizations, each of them containing a capacitor $C_i$ with charge $Q_i$ and an inductor $L_i$ with magnetic flux linkage $\phi_i$ for $i=1,\dots,10$ (see Figure \ref{circuit}). Using Kirchhoff laws, we obtain the following port-Hamiltonian system in normal form (\ref{circuit_pch1}) and (\ref{circuit_pch2}), where the Hamiltonian of the system is \begin{equation*}
	H(Q_1,\dots,Q_5,\phi_1,\dots,\phi_5)=\frac{Q^2_1}{2C_1}+\dots+\frac{Q^2_5}{2C_5}+\frac{\phi^2_1}{2L_1}+\dots+\frac{\phi^2_5}{2L_5}.
\end{equation*}

\begin{equation}\label{circuit_pch1}
	\begin{bmatrix}
		\dot{Q}_1 \\ \vdots \\ \dot{Q}_5 \\[3pt]	\dot{\phi}_1 \\ \vdots \\ \dot{\phi}_5
	\end{bmatrix}=\begin{bmatrix}
	0&\mathbb{I}_5\\-\mathbb{I}_5&0 \end{bmatrix} \cdot \begin{bmatrix}
		\frac{\partial H}{\partial Q_1} \\ \vdots \\ \frac{\partial H}{\partial Q_5} \\[5pt] 	\frac{\partial H}{\partial \phi_1} \\ \vdots \\ \frac{\partial H}{\partial \phi_5}
\end{bmatrix}+\begin{bmatrix}
0\\\vdots\\0\\1\\\vdots\\1
\end{bmatrix}\cdot u
\end{equation}
\begin{equation}\label{circuit_pch2}
	y = \frac{\partial H}{\partial \phi_1}+\frac{\partial H}{\partial \phi_2}+\dots+\frac{\partial H}{\partial \phi_5},
\end{equation}

This port-Hamiltonian system treats the power supply $V=u$ as input and the current through the power supply, that is $y$, as output. One verifies that such a system is {\it non-canonical}. Our purpose is to learn the input-output behavior of this system without any access to the internal physical state and training only with input-output observations. 

In our implementation, we choose for simplicity $C_i=1$ and $L_i=1$ for $i=1,\dots,5$. We choose to learn with a 10-dimensional observable Hamiltonian representation to show that the dynamics can be captured even in the non-canonical case. (Indeed, with our choice of $C_i$s and $L_i$s, the system is readily checked to be noncanonical). We randomly generate an initial condition for the ground-truth system and integrate it using Euler's method (see Appendix \ref{integrator_section} for more sophisticated structure-preserving integration methods) with a discretization step of $0.01$ for $1000 $ time steps. The input is chosen as $u(t)=\sin (t)$. The $1000$ pairs of input and output data will be used as training data. During the training phase, we estimate the initial state ${\bf x}\in \mathbb{R}^{10}$ as well as the parameters ${\bf d}\in \mathbb{R}_{+}^{5}$ and ${\bf v}\in \mathbb{R}^{10}$. This is carried out via gradient descent using a learning rate of $\lambda = 0.1$ for $500$ epochs. At each gradient descent iteration we integrate the state-space equations corresponding to the current parameter values over $1000$ times steps with Euler's method and then we compute the squared error with respect to the training set. 

\begin{figure*}[!htb]
	\centering
	\begin{subfigure}[t]{0.5\textwidth}
		\centering
		\includegraphics[height=2in]{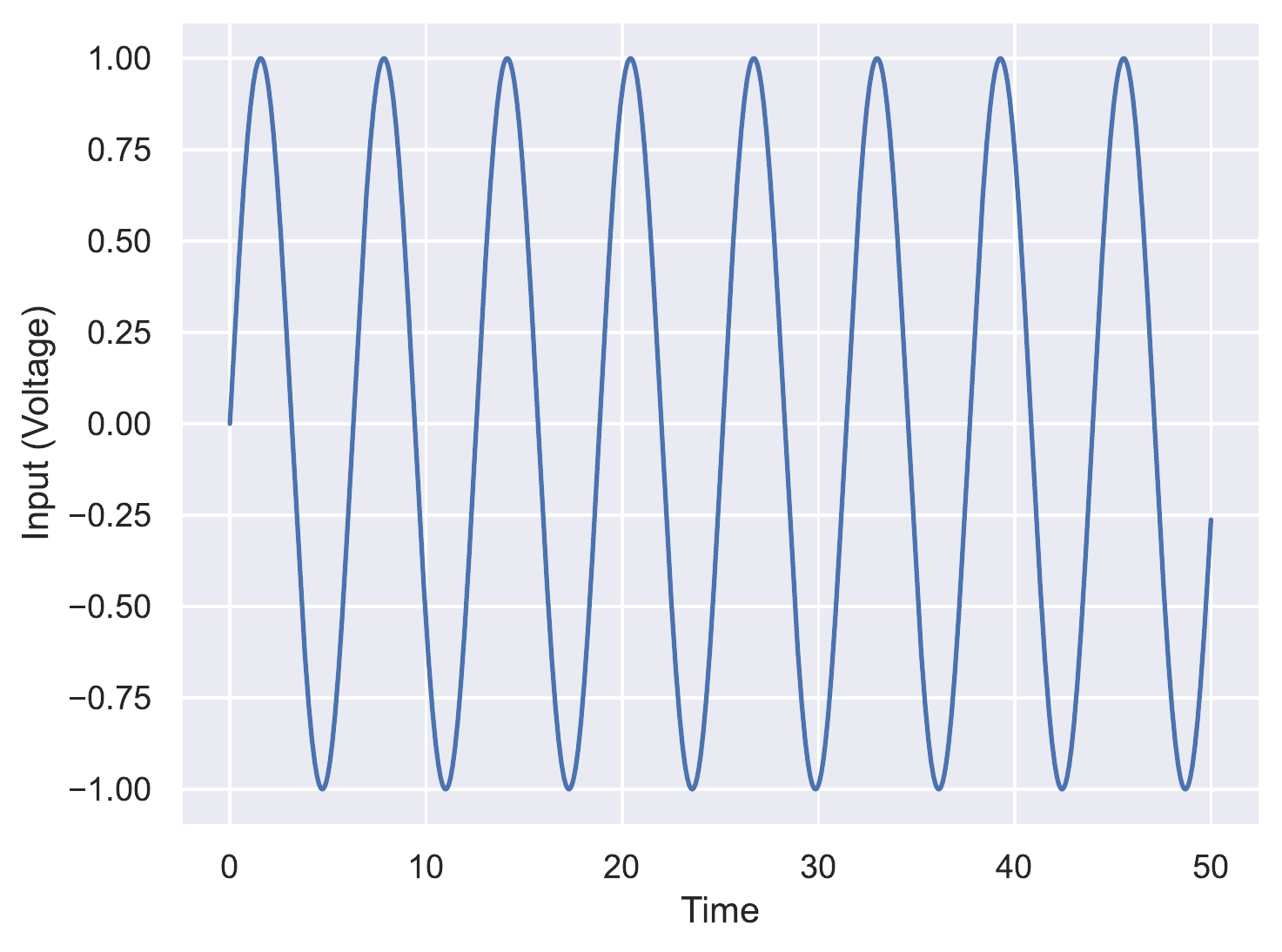}
		\caption{Input signal $u(t)$}
	\end{subfigure}%
	~ 
	\begin{subfigure}[t]{0.5\textwidth}
		\centering
		\includegraphics[height=2in]{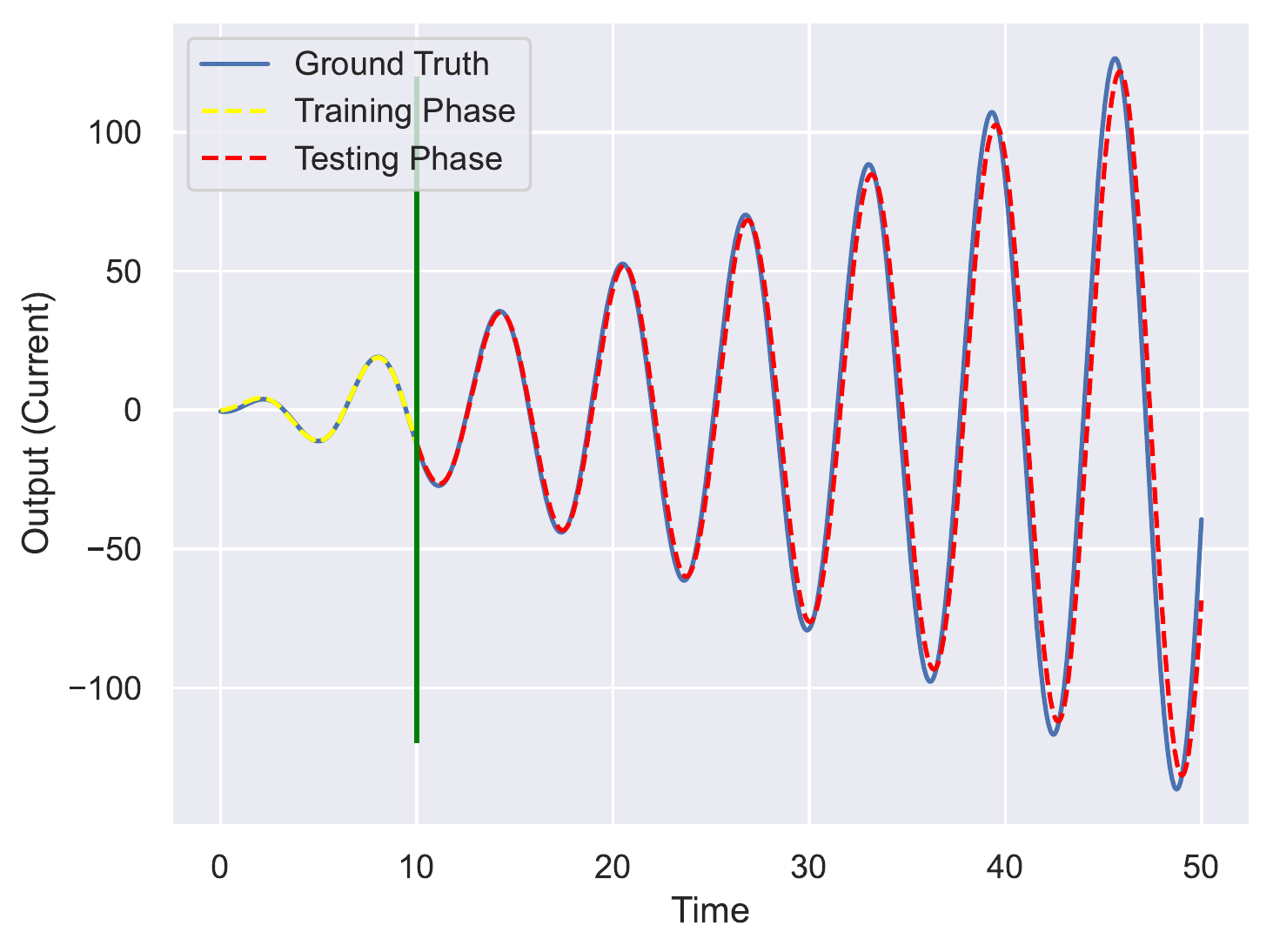}
		\caption{output $y(t)$}
	\end{subfigure}
	\caption{Training and testing on a sinusoidal signal.}
	\label{signal1}
\end{figure*}

\begin{figure*}[!htb]
	\centering
	\begin{subfigure}[t]{0.5\textwidth}
		\centering
		\includegraphics[height=2in]{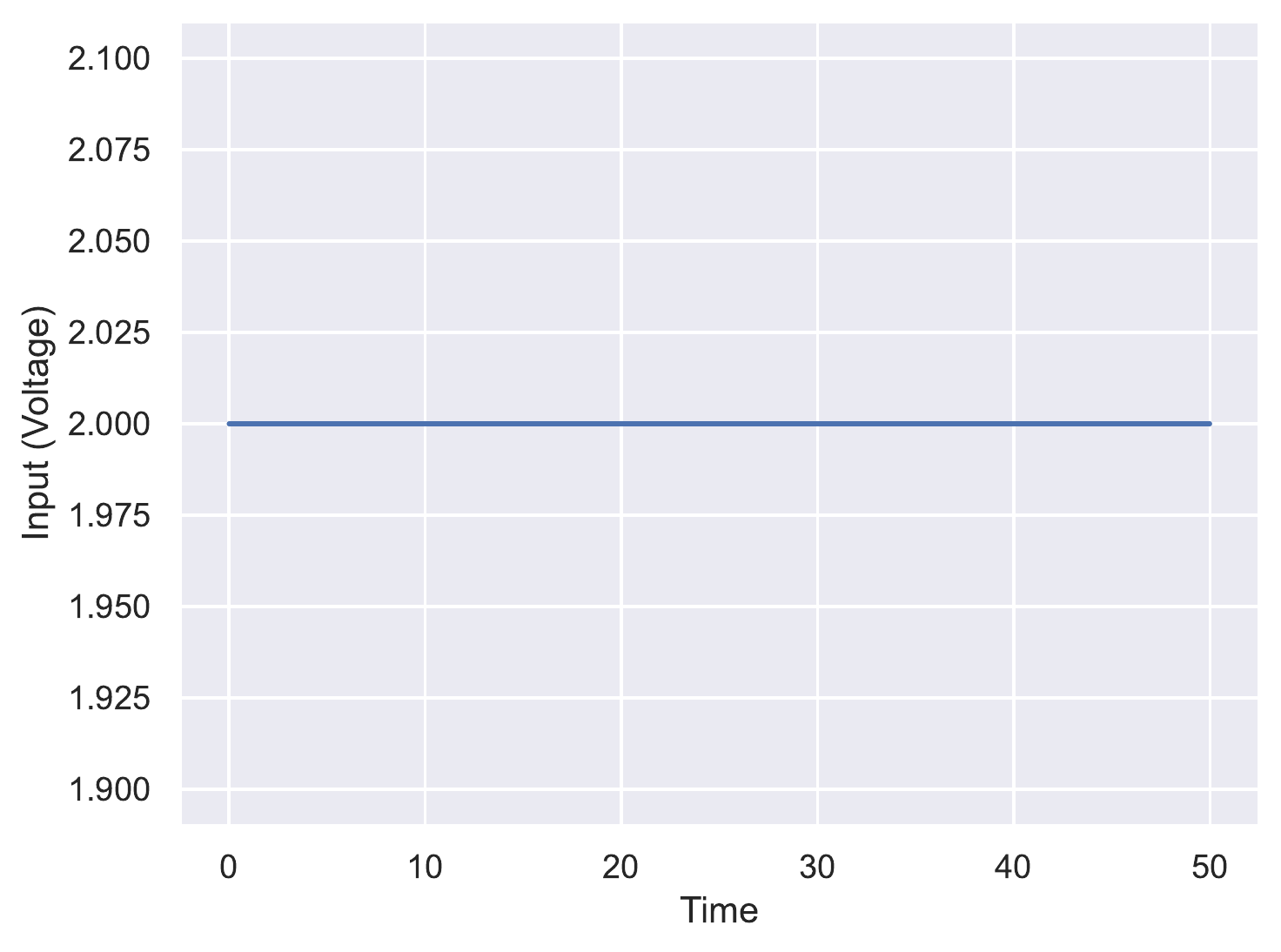}
		\caption{Input signal $u(t)$}
	\end{subfigure}%
	~ 
	\begin{subfigure}[t]{0.5\textwidth}
		\centering
		\includegraphics[height=2in]{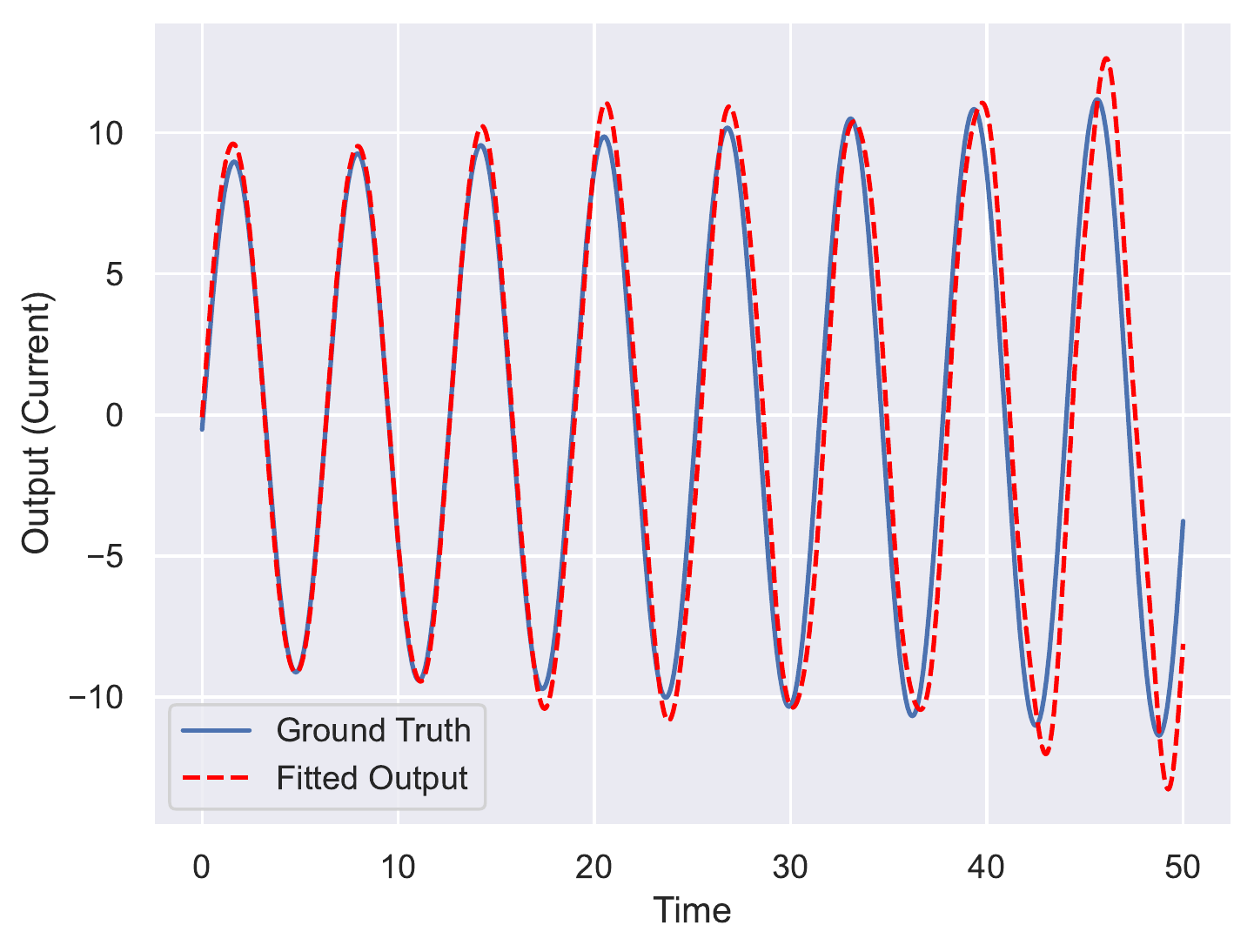}
		\caption{output $y(t)$}
	\end{subfigure}
	\caption{Testing on a constant signal. The training had been carried out using a sinusoidal signal. See Figure \ref{signal1}.}
	\label{signal2}
\end{figure*}

\begin{figure*}[!htb]
	\centering
	\begin{subfigure}[t]{0.5\textwidth}
		\centering
		\includegraphics[height=2in]{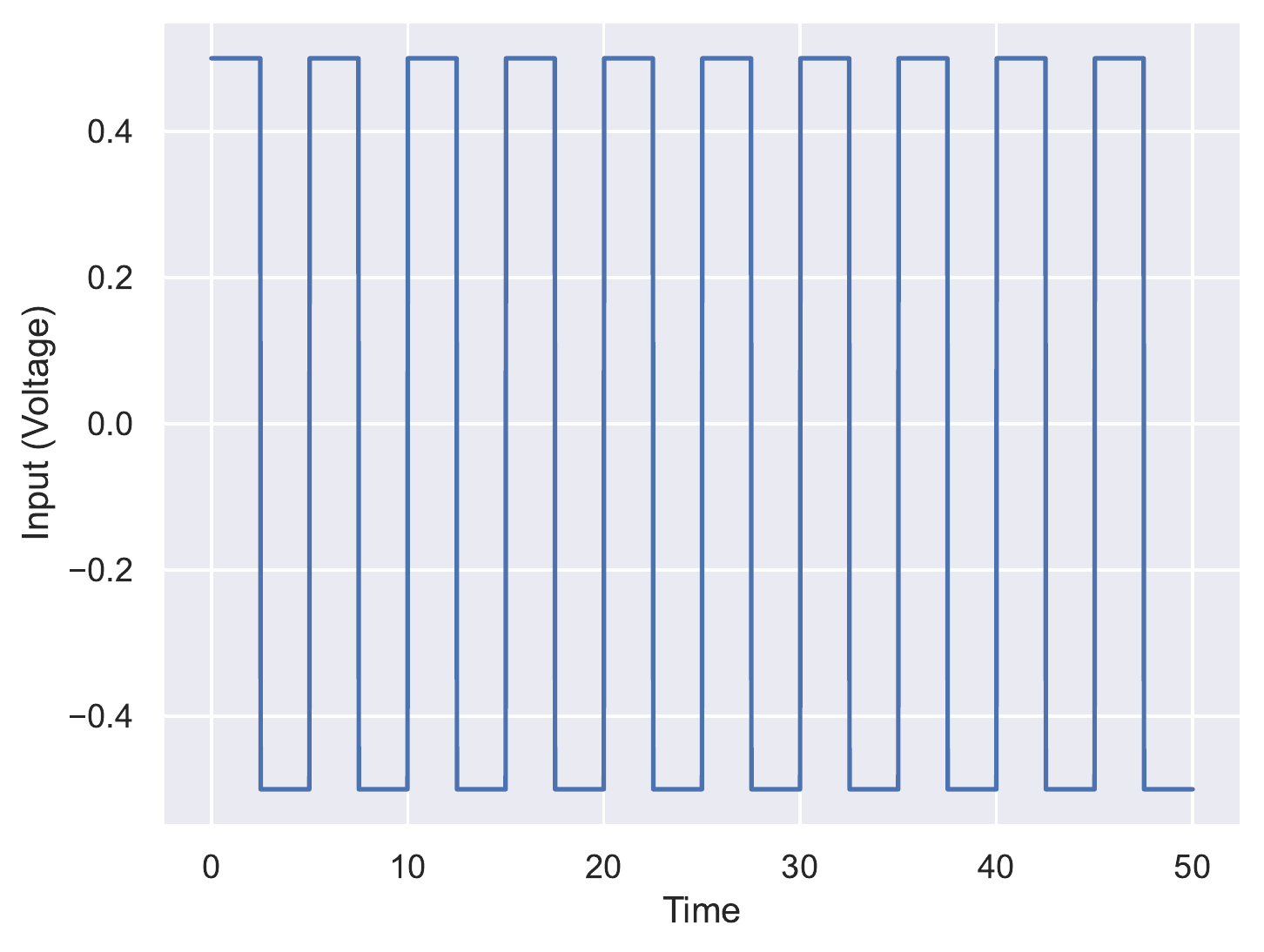}
		\caption{Input signal $u(t)$}
	\end{subfigure}%
	~ 
	\begin{subfigure}[t]{0.5\textwidth}
		\centering
		\includegraphics[height=2in]{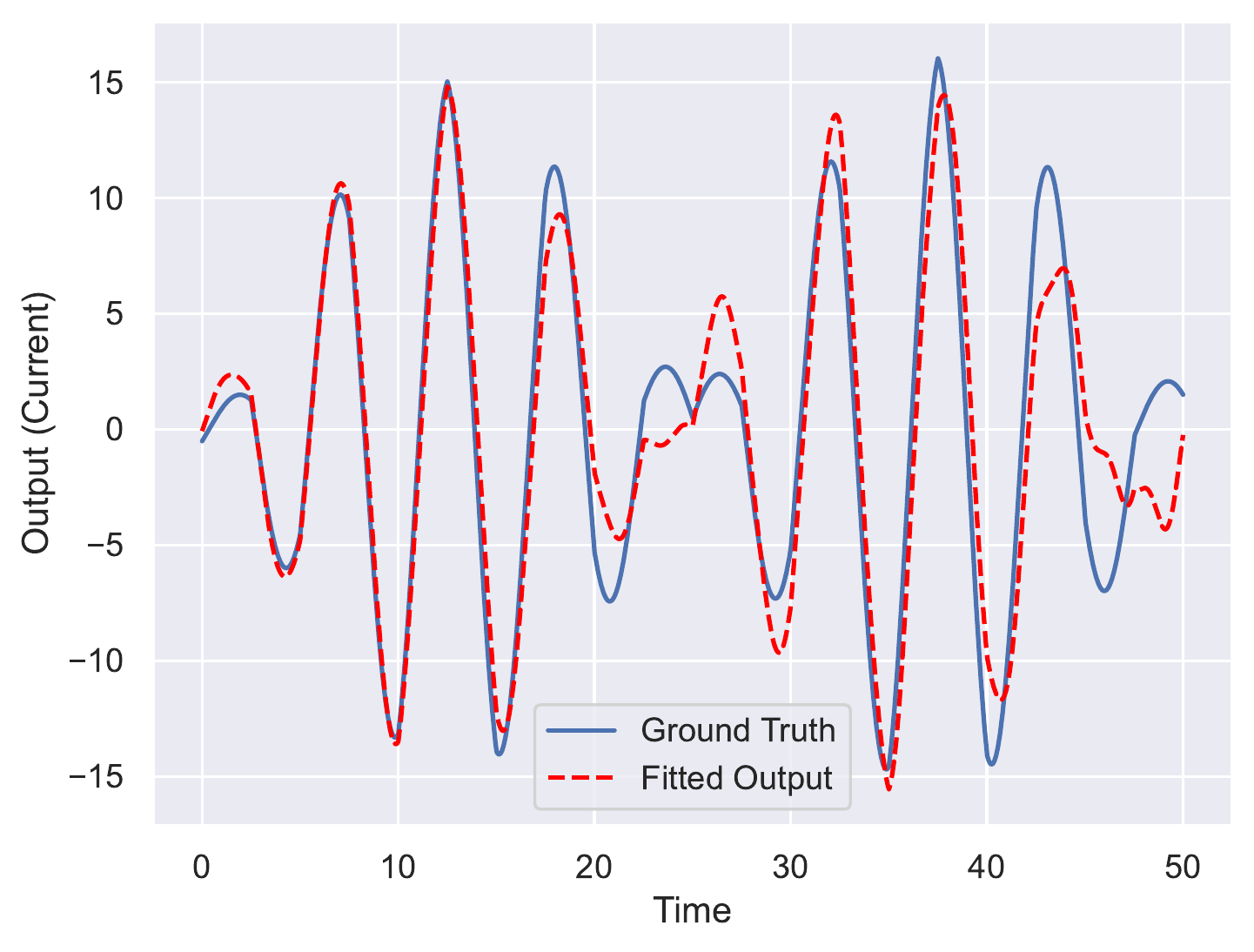}
		\caption{output $y(t)$}
	\end{subfigure}
	\caption{Testing on a square signal. The training had been carried out using a sinusoidal signal. See Figure \ref{signal1}.}
	\label{signal3}
\end{figure*}

\begin{figure*}[!htb]
	\centering
	\begin{subfigure}[t]{0.5\textwidth}
		\centering
		\includegraphics[height=2in]{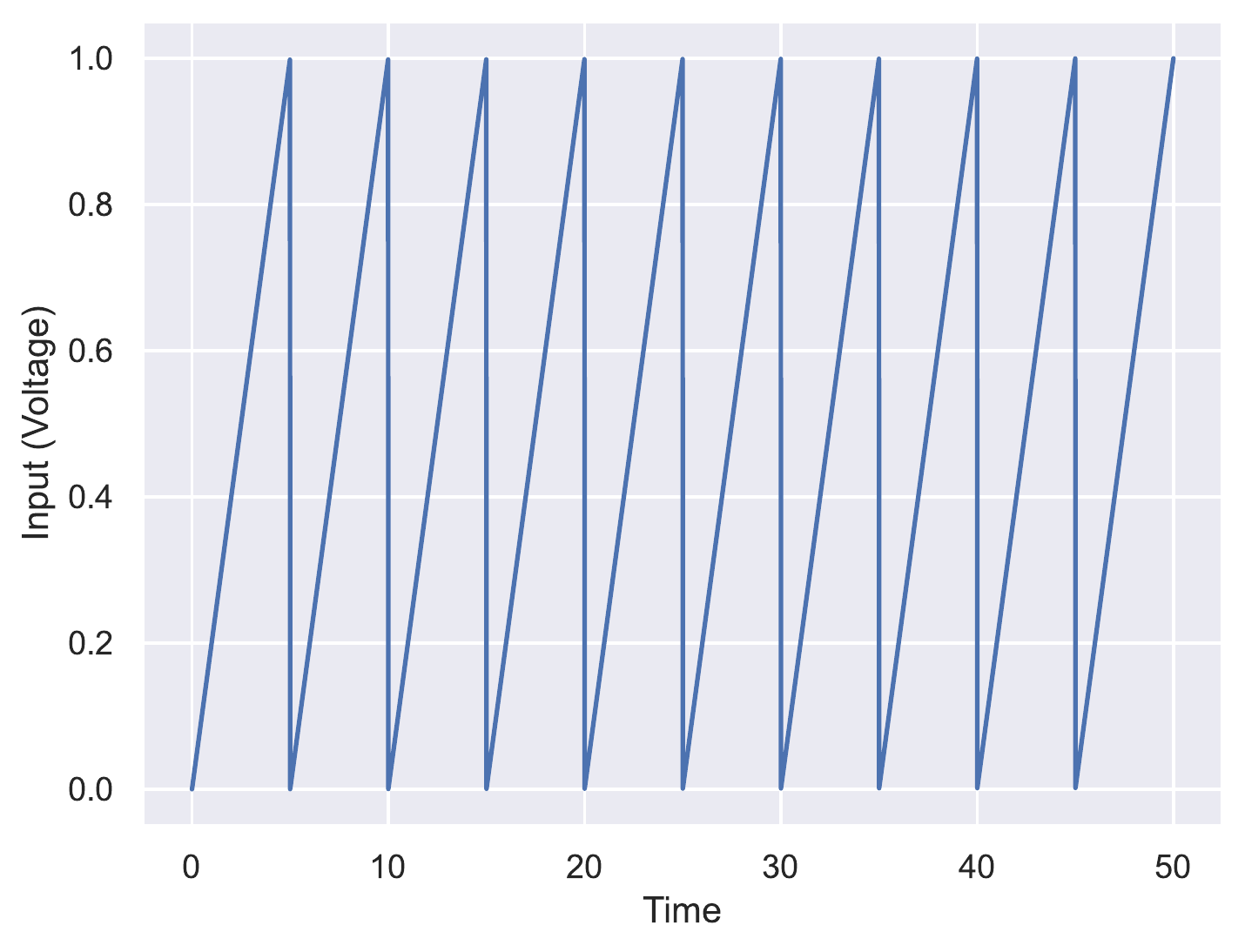}
		\caption{Input signal $u(t)$}
	\end{subfigure}%
	~ 
	\begin{subfigure}[t]{0.5\textwidth}
		\centering
		\includegraphics[height=2in]{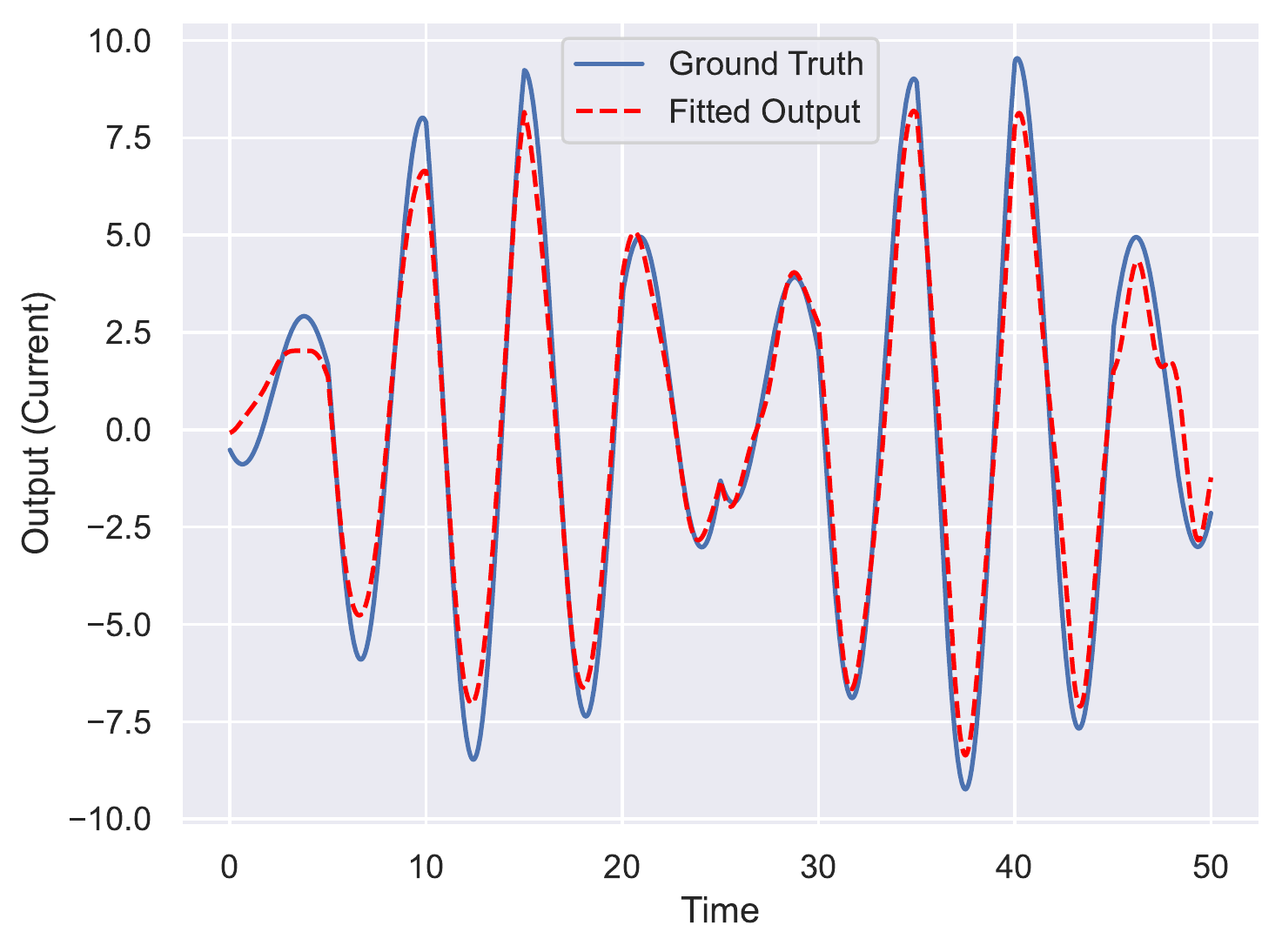}
		\caption{output $y(t)$}
	\end{subfigure}
	\caption{Testing on a ramp signal. The training had been carried out using a sinusoidal signal. See Figure \ref{signal1}.}
	\label{signal4}
\end{figure*}

We set a testing period of $4000$ time steps and demonstrate the robustness of our approach by not only testing our trained model on the original input $u(t)=\sin (t)$ but evaluating on other three commonly used input signals (see Figure \ref{signal1}, \ref{signal2}, \ref{signal3} and \ref{signal4}). The numerical experiments provide a strong indication that the underlying system is learned independently of the input signal and is robust with respect to various forms of inputs.

\subsection{Positive definite Frenkel-Kontorova model}
As a second example, we consider a modification of the well-known Frenkel-Kontorova model such that it becomes a linear port-Hamiltonian system with a positive-definite Hamiltonian function. Recall that the general form of Frenkel-Kontorova model describes the motion of classical particles with nearest neighbor interactions using periodic potentials. The Hamiltonian function can be written as 
\begin{equation*}
	H=\sum_{n=1}^N\left[\frac{1}{2}\cdot \dot{q}^2_n+\left(1-\cos q_n+\frac{1}{2}g\cdot(q_{n+1}-q_n-a_0)^2\right)\right].
\end{equation*}Since we are dealing with linear systems, we remove the periodic potential and rescale the potential coefficient. By fixing $a_0=0$, we obtain the Hamiltonian 
\begin{equation*}
H=\frac{1}{2}\cdot\sum_{n=1}^N\left[ \dot{q}^2_n+(q_{n+1}-q_n)^2\right].
\end{equation*}
 
 In order to consider a Hamiltonian that is strictly positive definite, we add a term $\frac{1}{2}q^2_1$ to the Hamiltonian, which carries the physical meaning that the particle $q_1$ interacts with the origin via a spring. In summary, our model of interest now has the positive-definite Hamiltonian 
\begin{equation*}
	H=\frac{1}{2}\cdot\sum_{n=1}^N\left[ \dot{q}^2_n+(q_{n+1}-q_n)^2\right]+\frac{1}{2}\cdot q^2_1=\frac{1}{2}\cdot\sum_{n=1}^N\left[ p^2_n+(q_{n+1}-q_n)^2\right]+\frac{1}{2}\cdot q^2_1.
\end{equation*}

For the sake of simplicity, consider a Hamiltonian systems with two unit mass particles (so that $p_i=\dot{q}_i$) and an external force $F=u$ that is imposed on the first particle. This gives a linear port-Hamiltonian system in normal form as below with the output being the velocity of the first particle.
\begin{equation}\label{fk_pch1}
	\begin{bmatrix}
		\dot{q}_1 \\ \dot{q}_2 \\[3pt]	\dot{p}_1  \\ \dot{p}_2
	\end{bmatrix}=\begin{bmatrix}
		0&\mathbb{I}_2\\-\mathbb{I}_2&0 \end{bmatrix} \cdot \begin{bmatrix}
		\frac{\partial H}{\partial q_1} \\[3pt] \frac{\partial H}{\partial q_2} \\[5pt] 	\frac{\partial H}{\partial p_1} \\[3pt] \frac{\partial H}{\partial p_2}
	\end{bmatrix}+\begin{bmatrix}
		0\\0\\1\\0
	\end{bmatrix}\cdot u
\end{equation}
\begin{equation}\label{fk_pch2}
	y = \frac{\partial H}{\partial p_1}.
\end{equation}

In contrast to the first example, this system is {\it canonical}. Therefore, based on our theoretical results, any input-output dynamics can be captured by either a controllable or an observable Hamiltonian representation and furthermore, it is possible to uniquely identify the system by learning the parameters in the quotient space $\mathcal{X}^{2}_{\uparrow}\times \mathbb{R}_{+}^{2}$. 

For the sake of the numerical illustration, we choose the initial state condition $\mathbf{x}=(2,1,-3,-3 )^T$ for the ground-truth system and integrate it $1000$ time steps times using Euler's method with step of $0.01$ (see Appendix \ref{integrator_section} for more sophisticated structure-preserving integration methods), where the input is chosen as $u(t)=\sin(t)$. The $1000$ pairs of input and output data are then used as training data. 

As motivated above, we apply two different training mechanisms in which we learn the initial state condition and the parameter values of the model using both the natural parameters from $\Theta_{CH_{n}}$  of the observable Hamiltonian representation and those in the unique identifiability space $\mathcal{X}^{2}_{\uparrow}\times \mathbb{R}_{+}^{2}$. As in the previous example, we carry out the training using gradient descent with a learning rate of $\lambda = 0.02$ over $1500$ epochs out of randomly chosen initial values for the initial state condition and the model parameters in $\Theta_{CH_{n}}$ and $\mathcal{X}_{\uparrow}^{2}\times \mathbb{R}_{+}^{2}$.

\begin{figure}[!htb]\centering
	\includegraphics[width=7cm]{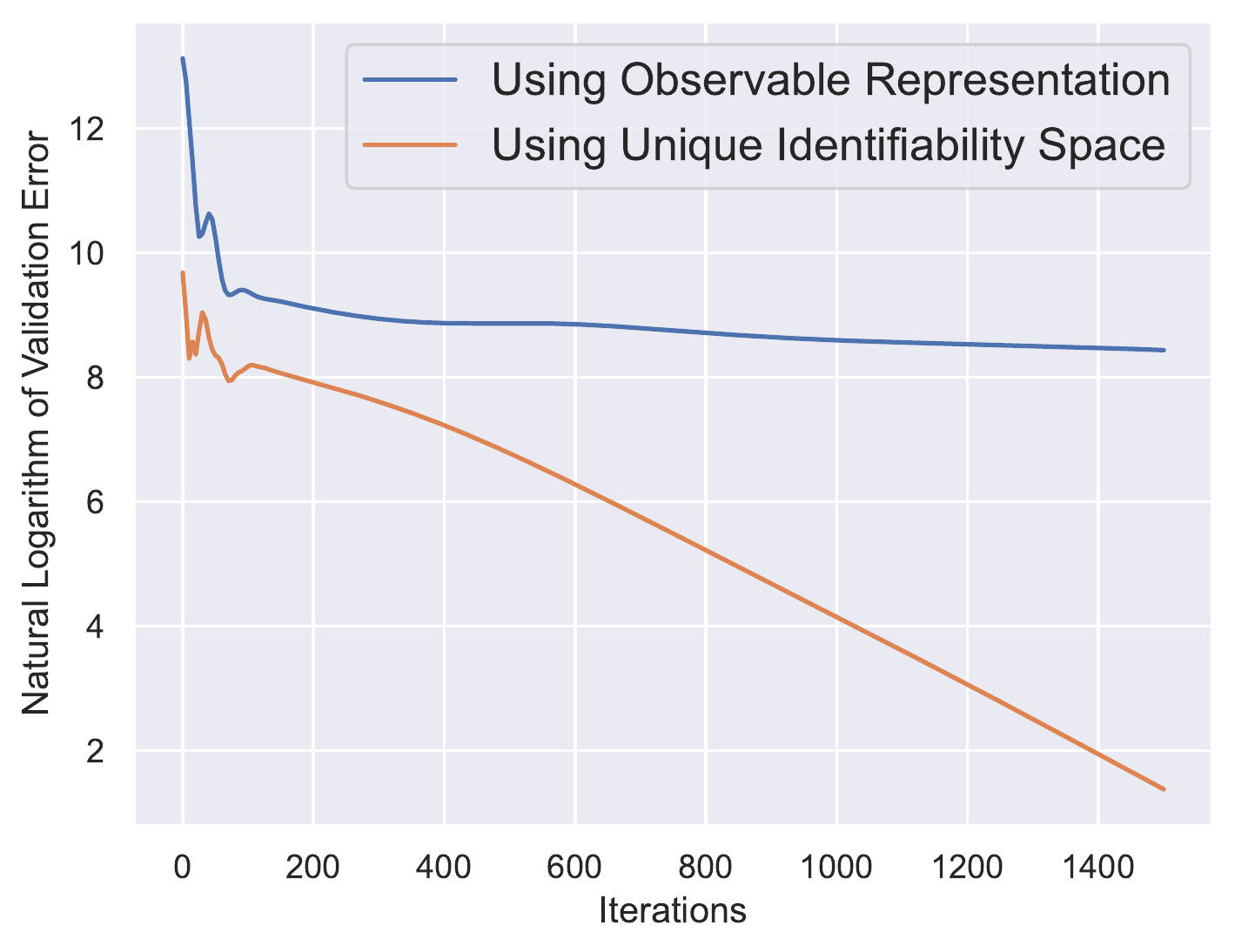}
	\caption{Logarithm of validation errors of the two training mechanisms based on using the natural parameters of the observable representation and the unique identifiability space}
	\label{validation_error}
\end{figure}

\begin{figure*}[!htb]
	\centering
	\begin{subfigure}[t]{0.5\textwidth}
		\centering
		\includegraphics[height=2in]{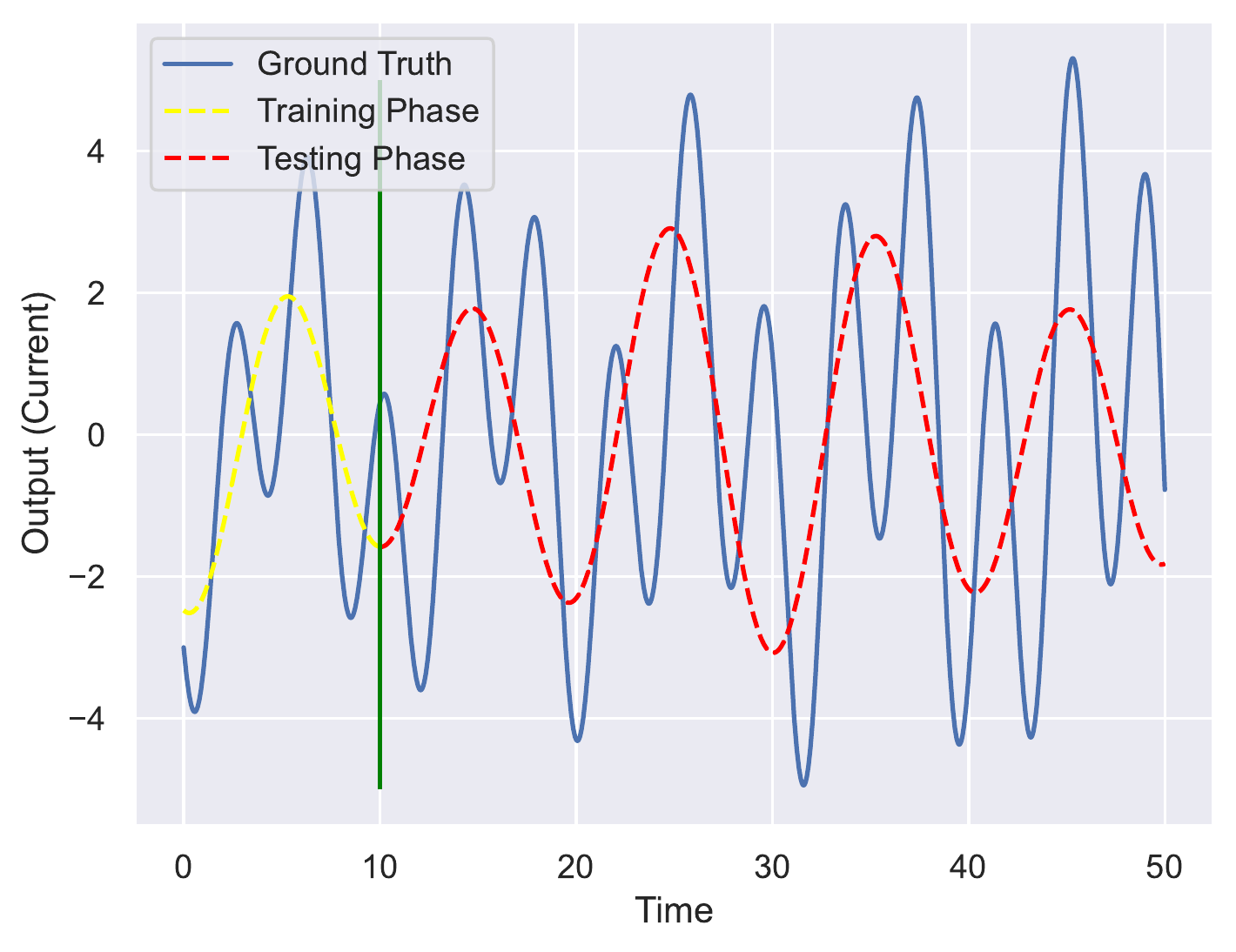}
		\caption{Input signal $u(t)$}
	\end{subfigure}%
	~ 
	\begin{subfigure}[t]{0.5\textwidth}
		\centering
		\includegraphics[height=2in]{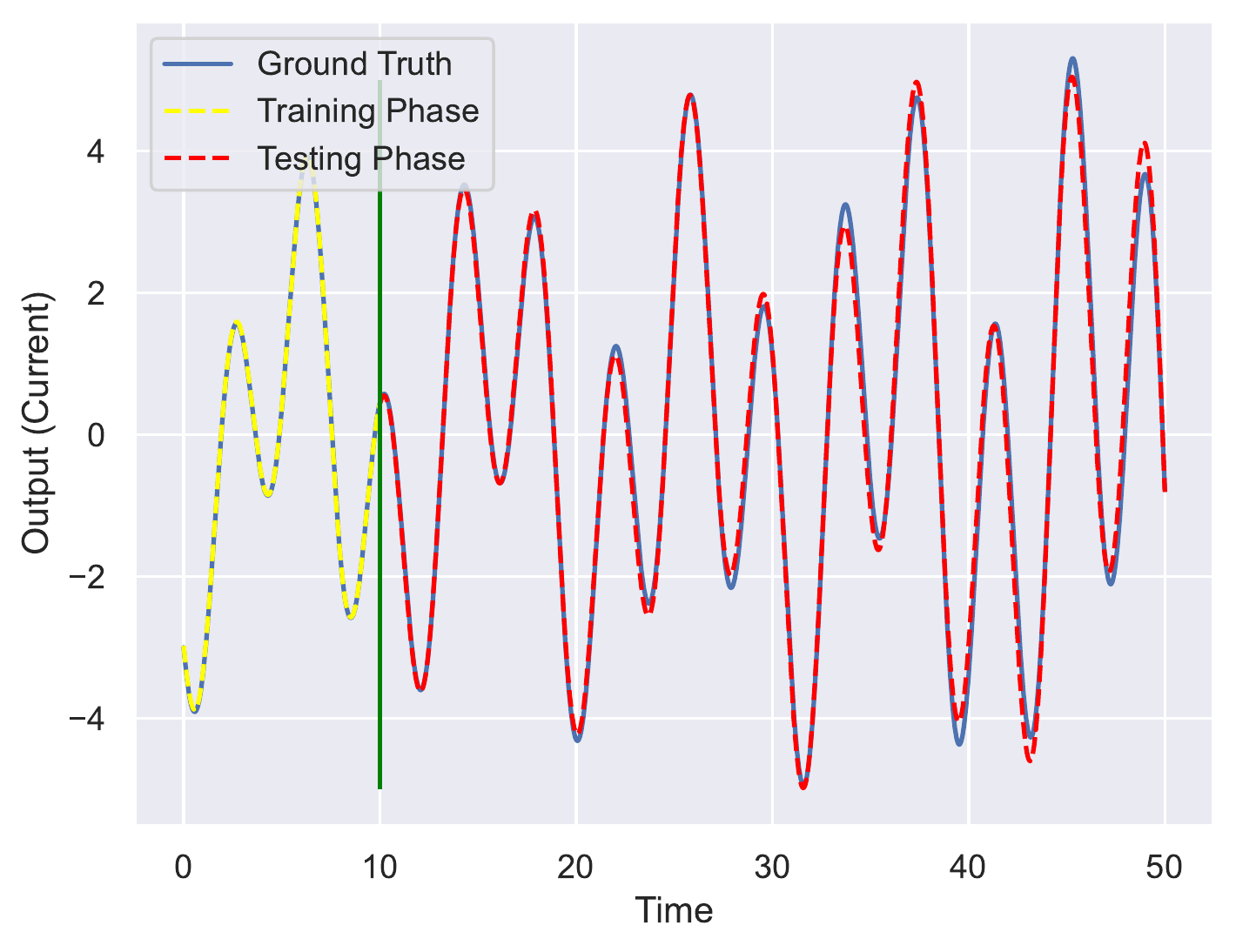}
		\caption{Input signal $u(t)$}
	\end{subfigure}
	\caption{Training and testing performance of the two training mechanisms after 1500 gradient descent iterations based on using the natural parameters of the observable representation (pane (a)) and the unique identifiability space (pane (b))}
	\label{FK_model}
\end{figure*}

We record the validation error during $1500$ gradient descent iterations of both training mechanisms to compare their convergence rates. Heuristically, it should be expected that the rate of convergence is faster when the models are trained using the coordinates that provide unique identifiability. This is empirically confirmed in Figure \ref{validation_error} (indeed, unique identifiability provides exponentially faster convergence). After $1500$ iterations, the prediction accuracy when training was carried out using the unique identifiability space significantly outperforms the other setting, as can be seen in Figure \ref{FK_model}. Moreover, we found that the learned parameters ${\bf d}\in \mathcal{X}^{2}_{\uparrow}$ are exactly the same as the eigenvalues of the Hamiltonian matrix, which is theoretically guaranteed by unique identifiability. It is worth emphasizing that despite the difference in the convergence rates both mechanisms eventually lead to perfect path continuations of the input-output dynamics after enough training iterations. 

\section{Conclusions}

In this paper we have introduced a complete structure-preserving learning scheme for single-input/single-output (SISO) linear port-Hamiltonian systems. The construction is based on the solution, when possible, of the unique identification problem for these systems, in ways that reveal fundamental relationships between classical notions in control theory and crucial properties in the machine learning context, like structure-preservation and expressive power. 

The main building block in our construction is a representation result that we introduced for linear port-Hamiltonian systems in normal form that provides two subfamilies of linear systems that are by construction controllable and observable (Definition \ref{definition_of_controllable}). We showed that morphisms can be established between the elements in these families and those in the category of normal form port-Hamiltonian systems. The existence of these morphisms immediately guarantees that the complexity of the family of port-Hamiltonian filters is actually not $\mathcal{O}(n^2)$, as it could be guessed from the standard parametrization of this family, but $\mathcal{O}(n)$. We showed that the expressive power of our proposed representations is limited for non-canonical port-Hamiltonian systems. Indeed, we saw that the observable representation is guaranteed to capture all possible input-output dynamics of port-Hamiltonian systems (full expressive power), but it does not always produce port-Hamiltonian dynamics (fails to be structure-preserving). In the controllable case, structure preservation is guaranteed, but there is, in general, no full expressive power. For canonical port-Hamiltonian systems, these representations are both structure-preserving and have full expressive power. 

We saw that even in the canonical situation, the availability of the controllable/observable representations did not yet provide a well-specified learning problem for this category since the invariance of these systems under system automorphisms implies the existence of symmetries (or degeneracies) in those parametrizations. We tackled this problem by solving the unique identifiability of input-output dynamics of linear port-Hamiltonian systems in normal form by characterizing the quotient space by system automorphisms as a {Lie groupoid} orbit space. Moreover, we showed that in the {canonical} case the corresponding quotient spaces can be characterized as orbit spaces with respect to an explicit  {group}  action and that can be explicitly endowed with a smooth manifold structure that has global Euclidean coordinates that can be used at the time of constructing estimation algorithms. 
 Consequently, we showed that canonical port-Hamiltonian dynamics can be identified fully and explicitly in either the controllable or the observable Hamiltonian representations and learned by estimating a unique set of parameters in a smooth manifold that is obtained as a group orbit space. Additionally, we complemented this learning scheme with results that allow us to extend it to situations where we remain agnostic as to the dimension of the underlying data-generating port-Hamiltonian system. 

We concluded the paper with some numerical examples that illustrate the viability of the method that we propose in systems with various levels of complexity and dimensions as well as the computational advantages associated with the use of the parameter space in which unique identification is guaranteed.

\addcontentsline{toc}{section}{Acknowledgments}
\section*{Acknowledgments}
The authors thank Lyudmila Grigoryeva for helpful discussions and remarks and acknowledge partial financial support from the Swiss National Science Foundation (grant number 175801/1) and the School of Physical and Mathematical Sciences of the Nanyang Technological University. DY is funded by the Nanyang President's Graduate Scholarship of Nanyang Technological University.

\bibliographystyle{wmaainf}
\addcontentsline{toc}{section}{Bibliography}
\bibliography{literature}

\section{Appendices}

\footnotesize

\subsection{Proof of Theorem~\ref{main_theorem} (i)}

	Let $({\bf d},{\bf v})\in \Theta_{CH_{n}}$ and let 
	\begin{equation}\label{proposed_system_learner for theorem}
		\left\{
		\begin{aligned}
			\dot{{\bf s}}&=g^{ctr}_1({\bf d})\cdot {\bf s}+\left(
			0 , 0 ,\cdots, 0, 1\right)^{T}\cdot u,\\
			y&=g^{ctr}_2({\bf d}, {\bf v})\cdot {\bf s},
		\end{aligned} \right.
	\end{equation} 
	be the corresponding linear controllable state-space system. In the following paragraphs, we construct for every $S\in Sp(2n,\mathbb{R})$, a linear system morphism $f_S^{({\bf d},{\bf v})}:\mathbb{R}^{2n}\rightarrow\mathbb{R}^{2n}$ between \eqref{proposed_system_learner for theorem} and the port-Hamiltonian system $(Q,B)=\varphi_S(\theta_{CH _n}({\bf d},{\bf v}))\in{PH}_n$ in the statement. Notice, first of all that $Q$ is by construction symmetric and positive-definite. Let now $L \in \mathbb{M}_{2n}$ be the matrix implementing the linear map $f_S^{({\bf d},{\bf v})} $, that is, $f_S^{({\bf d},{\bf v})}({\bf s})=L {\bf s} $, ${\bf s} \in {\Bbb R}^{2n} $. We now explicitly construct $L$ and prove that it provides a system morphism.
	We start by denoting $A:=\mathbb{J}_{n}\begin{bmatrix}
		D&0\\
		0&D
	\end{bmatrix}$, and define for each $k=1,\dots,2n$, a matrix $L_k\in \mathbb{M}^{2n}$ as
	\begin{equation*}
		L_{2n-k}:=A^{k}+a_{2n-1}\cdot A^{k-1}+\dots+a_{2n-k}\cdot \mathbb{I}_{2n}.
	\end{equation*}
	In particular, $L_{2n}=\mathbb{I}_{2n}$. Then, $L$ is constructed as $
	L^{\prime}:=\begin{bmatrix}
		L_1{\bf v}&L_2{\bf v}&\cdots&L_{2n}{\bf v}
	\end{bmatrix}$, and $L:=S^{-1}L^{\prime}$.
	
	We now check that $f_{S}^{({\bf d},{\bf v})}({\bf s})=L{\bf s}$ is indeed a system morphism between (\ref{proposed_system_learner for theorem}) and the port-Hamiltonian system (\ref{ph_definition}) with $Q=S^{T}\begin{bmatrix}
		D&0\\
		0&D
	\end{bmatrix}S$ and $B=S^{-1}{\bf v}$. This amounts to checking that \begin{description}
		\item[(i)] $L\cdot g^{ctr}_1({\bf d})=\mathbb{J}_nQL$
		\item[(ii)] $L\cdot \left(
		0 , 0 ,\cdots, 0, 1\right)^{T}=B$
		\item[(iii)] $g^{ctr}_2({\bf d},{\bf v})=B^TQL$.
	\end{description}
	We note that {\bf (ii)} trivially holds. Now, {\bf (i)} is equivalent to 
	\begin{align*}
		&S^{-1}L^{\prime}g^{ctr}_1({\bf d})=\mathbb{J}_nS^{T}\begin{bmatrix}
			D&0\\
			0&D
		\end{bmatrix}SS^{-1}L^{\prime}=S^{-1}\mathbb{J}_n\begin{bmatrix}
			D&0\\
			0&D
		\end{bmatrix}SS^{-1}L^{\prime}\\
		&\iff L^{\prime}g^{ctr}_1({\bf d})=\mathbb{J}_n\begin{bmatrix}
			D&0\\
			0&D
		\end{bmatrix}L^{\prime}\\
		&\iff \begin{bmatrix}
			L_1{\bf v}&L_2{\bf v}&\cdots&L_{2n}{\bf v}
		\end{bmatrix} \begin{bmatrix}
			0&1&0&\dots&0 \\
			0&0&1&\dots&0 \\
			\vdots&\vdots&\ddots&\vdots&\vdots\\
			0&0&0&\dots&1 \\
			-a_0&-a_1&-a_2&\dots&-a_{2n-1}
		\end{bmatrix}=A\begin{bmatrix}
			L_1{\bf v}&L_2{\bf v}&\cdots&L_{2n}{\bf v}
		\end{bmatrix}
	\end{align*}
	We compare the $k$-th columns of the left and the right-hand sides in this equality. When $k=1$, the difference between the first columns in the left and the right-hand side is 
	\begin{align}
		AL_1{\bf v}+a_0{\bf v}&=A(A^{2n-1}+a_{2n-1}\cdot A^{2n-2}+\dots+a_{1}\cdot \mathbb{I}){\bf v}+a_0{\bf v}\notag\\
		&=(A^{2n}+a_{2n-1}A^{2n-1}+\cdots+a_1A+a_0){\bf v}=0. \label{almost cayley}
	\end{align}
	The last equality holds as a consequence of the Cayley-Hamilton theorem. Indeed, by the definition of the entries $\left\{a _0, a _1, \ldots , a_{2n-1}\right\}$  we have that the characteristic polynomial of $A$ is
	\begin{align*}
		\det\left(\lambda \mathbb{I}_{2n}-A\right)&=\det\left(\lambda \mathbb{I}_{2n}-\begin{bmatrix}
			0&D\\
			-D&0
		\end{bmatrix}\right)\notag
		=\det\left(\begin{bmatrix}
			\lambda \mathbb{I}_{n}&-D\\
			D&\lambda \mathbb{I}_{n}
		\end{bmatrix}\right)\\
		&=\det\left(\lambda \mathbb{I}_{2n}\right)\cdot \det\left(\lambda \mathbb{I}_{2n}-\left(-D\right)\left(\frac{1}{\lambda} \mathbb{I} _n\right)D\right)\\
		&=({\lambda}^2+{d^2_1})({\lambda}^2+{d^2_2})\dots({\lambda}^2+{d^2_n})={\lambda}^{2n}+\sum_{i=0}^{2n-1}a_i\cdot{\lambda}^{i}.\notag
	\end{align*}
	Consequently, since by the Cayley-Hamilton theorem, $A$ solves its characteristic polynomial, we can conclude that $A^{2n}+a_{2n-1}A^{2n-1}+\cdots+a_1A+a_0 = 0$ and hence \eqref{almost cayley} follows. 
	When $1<k\leq 2n$, the difference between the $k$-th columns in the left and the right-hand side is
	\begin{align*}
		(L_{k-1}{\bf v}-a_{k-1}{\bf v})-AL_{k}{\bf v}=(L_{k-1}-a_{k-1}\mathbb{I}_{2n}-AL_k){\bf v}=0,
	\end{align*}
	since \begin{align*}
		L_{k-1}-a_{k-1}\mathbb{I}_{2n}-AL_k&=(A^{2n-k+1}+a_{2n-1}\cdot A^{2n-k}+\dots+a_{k-1}\cdot \mathbb{I}_{2n})-a_{k-1}\cdot \mathbb{I}_{2n}\\
		&-A(A^{2n-k}+a_{2n-1}\cdot A^{2n-k-1}+\dots+a_{k}\cdot \mathbb{I}_{2n})=0.
	\end{align*} 
	We have hence proved that {\bf (i)} holds. We now proceed to check {\bf (iii)}. This amounts to computing
	\begin{align}
		B^{T}QL&=(S^{-1}{\bf v})^TS^{T}\begin{bmatrix}
			D&0\\
			0&D
		\end{bmatrix}SS^{-1}L^{\prime}={\bf v}^T\begin{bmatrix}
			D&0\\
			0&D
		\end{bmatrix}L^{\prime} \nonumber \\
		&={\bf v}^T\begin{bmatrix}
			D&0\\
			0&D
		\end{bmatrix}\begin{bmatrix}
			L_1{\bf v}&L_2{\bf v}&\cdots&L_{2n}{\bf v}\label{coefficients_formula}
		\end{bmatrix}.
	\end{align}
	Let us denote \begin{equation}\label{coefficients_definition}
		B^{T}QL=\begin{bmatrix}
			c_{2n}\:\:c_{2n-1}\:\:c_{2n-2}\:\hdots\:\:c_{2}\:\:c_{1}
		\end{bmatrix}.
	\end{equation} Then we observe that for $k=1,\dots,n$, \begin{align*}
		c_{2k}&={\bf v}^{T}\begin{bmatrix}
			D&0\\
			0&D
		\end{bmatrix}L_{2n-2k+1}{\bf v}\\
		&={\bf v}^{T}\begin{bmatrix}
			D&0\\
			0&D
		\end{bmatrix}\left[\left(\mathbb{J}_n\begin{bmatrix}
			D&0\\
			0&D
		\end{bmatrix}\right)^{2k-1}+a_{2n-1}\left(\mathbb{J}_n\begin{bmatrix}
			D&0\\
			0&D
		\end{bmatrix}\right)^{2k-2}+\cdots+a_{2n-2k+1}\cdot \mathbb{I}_{2n}\right]{\bf v}\\
		&={\bf v}^{T}\begin{bmatrix}
			D&0\\
			0&D
		\end{bmatrix}\left[\left(\mathbb{J}_n\begin{bmatrix}
			D&0\\
			0&D
		\end{bmatrix}\right)^{2k-1}+a_{2n-2}\left(\mathbb{J}_n\begin{bmatrix}
			D&0\\
			0&D
		\end{bmatrix}\right)^{2k-3}+\cdots+a_{2n-2k+2}\left(\mathbb{J}_n\begin{bmatrix}
			D&0\\
			0&D
		\end{bmatrix}\right)\right]{\bf v}\\
		&={\bf v}^{T}\left[\mathbb{J}^{2k-1}_n\cdot\begin{bmatrix}
			D&0\\
			0&D
		\end{bmatrix}^{2k}+a_{2n-2}\cdot \mathbb{J}^{2k-3}_n\cdot\begin{bmatrix}
			D&0\\
			0&D
		\end{bmatrix}^{2k-2}+\cdots+a_{2n-2k+2}\cdot \mathbb{J}_n\cdot \begin{bmatrix}
			D&0\\
			0&D
		\end{bmatrix}^2\right]{\bf v}=0,
	\end{align*}
	The last equation follows from the fact that each summand is a skew-symmetric matrix. On the other hand, for $k=0,\dots,n-1$, \begin{align*}
		c_{2k+1}&={\bf v}^{T}\begin{bmatrix}
			D&0\\
			0&D
		\end{bmatrix}L_{2n-2k}{\bf v}\\
		&={\bf v}^{T}\begin{bmatrix}
			D&0\\
			0&D
		\end{bmatrix}\left[\left(\mathbb{J}_n\begin{bmatrix}
			D&0\\
			0&D
		\end{bmatrix}\right)^{2k}+a_{2n-1}\left(\mathbb{J}_n\begin{bmatrix}
			D&0\\
			0&D
		\end{bmatrix}\right)^{2k-1}+\cdots+a_{2n-2k}\cdot \mathbb{I}_{2n}\right]{\bf v}\\
		&={\bf v}^{T}\begin{bmatrix}
			D&0\\
			0&D
		\end{bmatrix}\left[\left(\mathbb{J}_n\begin{bmatrix}
			D&0\\
			0&D
		\end{bmatrix}\right)^{2k}+a_{2n-2}\left(\mathbb{J}_n\begin{bmatrix}
			D&0\\
			0&D
		\end{bmatrix}\right)^{2k-2}+\cdots+a_{2n-2k}\cdot \mathbb{I}_{2n}\right]{\bf v}\\
		&={\bf v}^{T}\left[(-1)^{k}\cdot \begin{bmatrix}
			D&0\\
			0&D
		\end{bmatrix}^{2k+1}+a_{2n-2}\cdot (-1)^{k-1}\cdot \begin{bmatrix}
			D&0\\
			0&D
		\end{bmatrix}^{2k-1}+\cdots+a_{2n-2k}\cdot \begin{bmatrix}
			D&0\\
			0&D
		\end{bmatrix}\right]{\bf v}.
	\end{align*}
	
	Substitute the values of coefficients $a_{2k}$ as expressions in terms of $d_i$'s, we obtain that \begin{equation*}
		c_{2k+1}={\bf v}^{T}\begin{bmatrix}
			F_k&0\\
			0&F_k
		\end{bmatrix} {\bf v},
	\end{equation*} for $k=0,\dots,n-1$, and
	\begin{equation*}
		F_k=\begin{bmatrix}
			f_1&&\\
			&f_2&&\giantzero\\
			&&\ddots\\
			&\giantzero&&f_{n-1}\\
			&&&&f_n\\
		\end{bmatrix}
	\end{equation*}
	with $f_l=d_l\cdot\sum_{\substack{j_1,\dots,j_k\neq l\\1\leq j_1<\dots<j_k\leq n}}\big(d_{j_1}d_{j_2}\cdots d_{j_k}\big)^2$, $l=1,\dots,n$. This is exactly how we define $g^{ctr}_2({\bf d},{\bf v})$. Hence, {\bf (iii)} is also verified.

\subsection{Proof of Theorem~\ref{main_theorem} (ii)}

	Let $(Q,B)\in{PH}_n$. Obtain ${\bf d}$ and ${\bf v}$ from $(Q,B)$ as in the statement of the theorem. We aim to construct a linear system morphism $f_S^{(Q,B)}:\mathbb{R}^{2n}\rightarrow\mathbb{R}^{2n}$ between the port-Hamiltonian system $(Q,B)\in{PH}_n$ and the observable Hamiltonian representation associated to $({\bf d},{\bf v})\in \Theta_{OH _n}$, that is,  
	\begin{equation}
		\label{proposed_system_learner2_for_theorem}
		\left\{
		\begin{aligned}
			\dot{{\bf s}}&=g^{obs}_1({\bf d})\cdot {\bf s}+g^{obs}_2({\bf d},{\bf v})\cdot u,\\
			y&=\left(
			0 , 0 ,\cdots, 0, 1\right)\cdot {\bf s}.
		\end{aligned} \right.
	\end{equation} 
	
	Denote by $L \in \mathbb{M}_{2n}$ the matrix implementing the linear map $f_S^{(Q,B)}$, that is, $f_S^{(Q,B)}({\bf s})=L {\bf s} $, ${\bf s} \in {\Bbb R}^{2n} $.  We now construct a $L$ which yields a system morphism.
	We start by writing $A:=\mathbb{J}_{n}\begin{bmatrix}
		D&0\\
		0&D
	\end{bmatrix}$ and define, for each $k=1,\dots,2n$, a matrix $L_k\in \mathbb{M}^{2n}$ as
	\begin{align*}
		L_{2n-k}&:=\left(\mathbb{J}_{n}Q\right)^{k}+a_{2n-1}\cdot \left(\mathbb{J}_{n}Q\right)^{k-1}+\dots+a_{2n-k}\cdot \mathbb{I}_{2n}\\
		&=\left(S^{-1}AS\right)^{k}+a_{2n-1}\cdot \left(S^{-1}AS\right)^{k-1}+\dots+a_{2n-k}\cdot \mathbb{I}_{2n}\\
		&=S^{-1}(A^{k}+a_{2n-1}\cdot A^{k-1}+\dots+a_{2n-k}\cdot \mathbb{I}_{2n})\cdot S.
	\end{align*}
	In particular, $L_{2n}=\mathbb{I}_{2n}$. Then, define $L$ is as $
	L:=\begin{bmatrix}
		B^TQL_1\\B^TQL_2\\\vdots\\B^TQL_{2n}
	\end{bmatrix}_{2n\times 2n}$.
	
	We now check that $f_S^{(Q,B)}({\bf s})=L{\bf s}$ is indeed a system morphism between the port-Hamiltonian system (\ref{ph_definition}) and the observable Hamiltonian representation (\ref{proposed_system_learner2_for_theorem}) with $Q=S^{T}\begin{bmatrix}
		D&0\\
		0&D
	\end{bmatrix}S$ and $B=S^{-1}{\bf v}$. This amounts to checking that
	\begin{description}
		\item[(i)] $g^{obs}_1({\bf d})\cdot L=L\mathbb{J}_nQ$
		\item[(ii)] $LB=g^{obs}_2({\bf d},{\bf v})$
		\item[(iii)] $B^TQ=\left(
		0 , 0 ,\cdots, 0, 1\right)\cdot L$.
	\end{description}
	
	We note that {\bf (iii)} is straightforward. Now, {\bf (i)} is equivalent to 
	\begin{align*}
		\begin{bmatrix}
			0&0&\dots&0&-a_0 \\
			1&0&\dots&0&-a_1 \\
			0&1&\ddots&0&-a_2\\
			\vdots&\vdots&\cdots&\vdots&\vdots\\
			0&0&\dots&1&-a_{2n-1}
		\end{bmatrix}\begin{bmatrix}
			B^TQL_1\\B^TQL_2\\\vdots\\B^TQL_{2n}
		\end{bmatrix}=\begin{bmatrix}
			B^TQL_1\\B^TQL_2\\\vdots\\B^TQL_{2n}
		\end{bmatrix}\cdot S^{-1}AS.
	\end{align*}
	
	Compare now the $k$-th rows of the left and the right-hand sides of this equality. When $k=1$, the difference between the first rows in the left and the right-hand sides are 
	\begin{align*}
		B^{T}QL_1S^{-1}AS+a_0{B^{T}QL_{2n}}&=B^{T}Q(L_1S^{-1}AS+a_0\mathbb{I}_{2n})\\
		&=B^{T}QS^{-1}(A^{2n}+a_{2n-1}A^{2n-1}+\cdots+a_1A+a_0\cdot \mathbb{I}_{2n})S=B^{T}Q\cdot 0=0.
	\end{align*} 
	The last equality follows, as in the proof of Theorem \ref{main_theorem}, from the Cayley-Hamilton theorem.
	
	When $1<k\leq 2n$, the difference between the $k$-th rows in the left and the right-hand sides are:
	\begin{align*}
		B^TQL_{k-1}&-a_{k-1}B^TQL_{2n}-B^TQL_kS^{-1}AS=B^{T}Q(L_{k-1}-a_{k-1}\cdot\mathbb{I}_{2n}-L_kS^{-1}AS)\\
		&=B^{T}QS^{-1}\bigg[(A^{2n-k+1}+a_{2n-1}\cdot A^{2n-k}+\dots+a_{k-1}\cdot \mathbb{I}_{2n})-a_{k-1}\cdot\mathbb{I}_{2n}\\
		&~~~~~~~~~~~-(A^{2n-k}+a_{2n-1}\cdot A^{2n-k-1}+\dots+a_{k}\cdot \mathbb{I}_{2n})A\bigg]S=0,
	\end{align*} 
	which shows that {\bf (i)} holds. We now proceed to check {\bf (ii)}. This is equivalent to computing 
	\begin{align*}
		LB=\begin{bmatrix}
			B^TQL_1\\B^TQL_2\\\vdots\\B^TQL_{2n}
		\end{bmatrix}B.
	\end{align*} Let us denote $LB=\begin{bmatrix}
		c_{2n}\:\:c_{2n-1}\:\:c_{2n-2}\:\hdots\:\:c_{2}\:\:c_{1}
	\end{bmatrix}^T$. Then we have, for $k=1,\dots,2n$, \begin{align*}
		c_{2n-k+1}&=B^TQL_kB=(S^{-1}{\bf v})^TS^{T}\begin{bmatrix}
			D&0\\
			0&D
		\end{bmatrix}SS^{-1}(A^{2n-k}+a_{2n-1}\cdot A^{2n-k-1}+\dots+a_{k}\cdot \mathbb{I}_{2n})SS^{-1}{\bf v}\\
		&={\bf v}^T\begin{bmatrix}
			D&0\\
			0&D
		\end{bmatrix}(A^{2n-k}+a_{2n-1}\cdot A^{2n-k-1}+\dots+a_{k}\cdot \mathbb{I}_{2n}){\bf v},
	\end{align*} which coincides exactly with the expression of $c_{2n-k+1}$ in the equations (\ref{coefficients_formula}) and (\ref{coefficients_definition}) that we provided in the controllable Hamiltonian case. Thus, for {\bf (iii)} to hold, we simply need to require that $g^{obs}_2({\bf d},{\bf v})=(g^{ctr}_2({\bf d},{\bf v}))^T$.

\subsection{Proof of Theorem~\ref{equivalence_conditon_on_D,v}}

\noindent {\bf  Proof of part (i)}.
We have that $({\bf d}_1,{\bf v}_1)\sim_{sys}({\bf d}_2,{\bf v}_2)$ implies the existence of an invertible matrix $L$ such that \begin{equation*}
	\left\{
	\begin{aligned}
		&L\cdot g^{ctr}_1({\bf d}_1)=g^{ctr}_1({\bf d}_2)\cdot L\\
		&L\cdot \left(
		0 , 0 ,\cdots, 0, 1\right)^{T}=\left(
		0 , 0 ,\cdots, 0, 1\right)^{T}\\
		&g^{ctr}_2({\bf d}_1,{\bf v}_1)=g^{ctr}_2({\bf d}_2,{\bf v}_2)\cdot L
	\end{aligned} \right.
\end{equation*}
The first condition implies that $\det(\lambda\mathbb{I}-g^{ctr}_1({\bf d}_1))=\det(\lambda\mathbb{I}-g^{ctr}_1({\bf d}_2))$, meaning that \begin{equation*}
	\begin{aligned}
		&\big(\lambda^2+d_{1,1}^2\big)\dots\big(\lambda^2+d_{1,n}^2\big)=\big(\lambda^2+d_{2,1}^2\big)\dots\big(\lambda^2+d_{2,n}^2\big).
	\end{aligned}
\end{equation*} 
Therefore, (i) is clear. With the symmetry in (i), it is clear that $g^{ctr}_1({\bf d}_1)=g^{ctr}_2({\bf d}_2)$. Note that the second condition says the last column of $L$ is $\left(
0 , 0 ,\cdots, 0, 1\right)^{T}$. Bring both facts into the first condition $L\cdot g^{ctr}_1({\bf d}_1)=g^{ctr}_1({\bf d}_2)\cdot L$ and compare both sides. This will deduce $L$ can only be the identity. Thus the third condition becomes $g^{ctr}_2({\bf d}_1,{\bf v}_1)=g^{ctr}_2({\bf d}_2,{\bf v}_2)$, which is exactly (ii).\\
Conversely, with (i) and (ii) hold, we can check $L$ being identity works. Thus, $({\bf d}_1,{\bf v}_1)\sim_{sys}({\bf d}_2,{\bf v}_2)$.\\

\noindent {\bf  Proof of part (ii)}.
Since $\theta_{CH_n}({\bf d}_i,{\bf v}_i)$ ($i=1,2$) are linear systems, we can explicitly write down the filters with the initial conditions set to zero as
\begin{align*}
	(y_i(u))_t &= g^{ctr}_2({\bf d}_i,{\bf v}_i)\int_{0}^{t}e^{g^{ctr}_1({\bf d}_i)(t-s)}\cdot \left(
	0 , 0 ,\cdots, 0, 1\right)^{T}u(s)ds\\
	&= g^{ctr}_2({\bf d}_i,{\bf v}_i)\int_{0}^{t}\left[\mathbb{I}+g^{ctr}_1({\bf d}_i)(t-s)+(g^{ctr}_1({\bf d}_i))^2\frac{(t-s)^2}{2!}+\cdots\right]\cdot \left(
	0 , 0 ,\cdots, 0, 1\right)^{T}u(s)ds.
\end{align*}
By differentiating the above and using the fact that the input $u(t)$ is arbitrary, we see that $y_1$ and $y_2$ coincides as filters if and only if $g^{ctr}_2({\bf d}_1,{\bf v}_1)(g^{ctr}_1({\bf d}_1))^k\cdot \left(
0 , 0 ,\cdots, 0, 1\right)^{T}=g^{ctr}_2({\bf d}_2,{\bf v}_2)(g^{ctr}_1({\bf d}_2))^k\cdot \left(
0 , 0 ,\cdots, 0, 1\right)^{T}$ for all $k\in\mathbb{N}$. Moreover, one verifies that $g^{ctr}_2({\bf d},{\bf v})(g^{ctr}_1({\bf d}))^k\cdot \left(
0 , 0 ,\cdots, 0, 1\right)^{T}=0$ for odd $k$. Thus, if we define $e_i({\bf d},{\bf v})=g^{ctr}_2({\bf d},{\bf v})(g^{ctr}_1({\bf d}))^{2i}\cdot \left(
0 , 0 ,\cdots, 0, 1\right)^{T}$, then $y_1=y_2$ is equivalent to $e_i({\bf d}_1,{\bf v}_1)=e_i({\bf d}_2,{\bf v}_2)$ for all $i\in \mathbb{N}$.

Now, one finds a recursion in the values of $e_i$'s. More precisely, for all $m\geq n$
\begin{equation*}
	e_m({\bf d},{\bf v}) = -a_{2n-2}\cdot e_{m-1}({\bf d},{\bf v})-a_{2n-4}\cdot e_{m-2}({\bf d},{\bf v})-\cdots-a_2\cdot e_1({\bf d},{\bf v})-a_0\cdot e_{m-n}({\bf d},{\bf v}).
\end{equation*}
On the other hand, $(a_{2n-2},\dots,a_2,a_0)$ happens to be the coefficients of the characteristic polynomial of $g^{ctr}_1({\bf d},{\bf v})$, therefore, by Cayley-Hamilton Theorem, $e_m({\bf d},{\bf v})=0$ for all $m\geq n$.

  In conclusion, $\theta_{CH_n}({\bf d}_1,{\bf v}_1)$ and $\theta_{CH_n}({\bf d}_2,{\bf v}_2)$ induce the same filter if and only if $e_i({\bf d}_1,{\bf v}_1)=e_i({\bf d}_2,{\bf v}_2)$ for all $0\leq i\leq n-1$, where it can be checked that the following recursion holds true
	\begin{align*}
	\begin{split}
		e_1 &= c_1\\
		e_2 &= c_3-a_{2n-2}\cdot e_1\\
		e_3&=c_5-a_{2n-2}\cdot e_2-a_{2n-4}\cdot e_1\\
		\vdots\\
		e_n &= c_{2n-1}-a_{2n-2}\cdot e_{n-1}-a_{2n-4}\cdot e_{n-2}-\cdots-a_2\cdot e_1.
	\end{split}
\end{align*}

\subsection{Proof of Theorem~\ref{characterization}}

	\noindent {\bf  $\varphi_S\circ \theta_{CH _n}$ is compatible with $\sim_{\star}$ and $\sim_{sys}$.} Fix a choice of $S\in Sp(2n,\mathbb{R})$. We need to show that $({\bf d}_1,{\bf v}_1)\sim_\star({\bf d}_2,{\bf v}_2)$ if and only if $\bigg(S^{T}\begin{bmatrix}
		D_1&0\\
		0&D_1
	\end{bmatrix}S,S^{-1}{\bf v}_1\bigg):=(Q_1,B_1)\sim_{sys}(Q_2,B_2):=\bigg(S^{T}\begin{bmatrix}
		D_2&0\\
		0&D_2
	\end{bmatrix}S,S^{-1}{\bf v}_2\bigg)$. This means there exists an invertible $L$ such that (\ref{automorphism_condition}) holds. We claim that $L=S^{-1}PAS$ does the job, where $P$ and $A$ are given by Definition \ref{equivalence_relation}. \par
	
	The first condition is \begin{equation*}
		\begin{aligned}
			&L\mathbb{J}Q_1=\mathbb{J}Q_2L\\
			&\iff L\mathbb{J}Q_1L^{-1}=\mathbb{J}Q_2\\
			&\iff
			S^{-1}PAS\mathbb{J}S^{T}\begin{bmatrix}
				D_1&0\\
				0&D_1
			\end{bmatrix}SS^{-1}A^{-1}P^{-1}S=\mathbb{J}S^{T}\begin{bmatrix}
				D_2&0\\
				0&D_2
			\end{bmatrix}S\\
			&\iff
			S^{-1}PAS\mathbb{J}S^{T}\begin{bmatrix}
				D_1&0\\
				0&D_1
			\end{bmatrix}A^{-1}P^{-1}S=S^{-1}\mathbb{J}\begin{bmatrix}
				D_2&0\\
				0&D_2
			\end{bmatrix}S\\
			&\iff PA\mathbb{J}\begin{bmatrix}
				D_1&0\\
				0&D_1
			\end{bmatrix}A^{-1}P^{-1}=\mathbb{J}\begin{bmatrix}
				D_2&0\\
				0&D_2
			\end{bmatrix}\\
			&\iff A\mathbb{J}\begin{bmatrix}
				D_1&0\\
				0&D_1
			\end{bmatrix}A^{-1}=P^{T}\mathbb{J}\begin{bmatrix}
				D_2&0\\
				0&D_2
			\end{bmatrix}P\\
			&\iff A\mathbb{J}\begin{bmatrix}
				D_1&0\\
				0&D_1
			\end{bmatrix}=\mathbb{J}\begin{bmatrix}
				D_1&0\\
				0&D_1
			\end{bmatrix}A.
		\end{aligned}
	\end{equation*}
	The second condition is true by construction, namely
	\begin{equation*}
		LB_1=B_2\iff S^{-1}PASS^{-1}{\bf v}_1=S^{-1}{\bf v}_2\\
		\iff {\bf v}_2=PA{\bf v}_1.
	\end{equation*}
	
	The third condition is 
	\begin{equation*}
		\begin{aligned}
			&B^{T}_1Q_1=B^{T}_2Q_2L\\
			&\iff {\bf v}^{T}_1S^{-T}S^{T}\begin{bmatrix}
				D_1&0\\
				0&D_1
			\end{bmatrix}S\\&~~~~~~~~~~={\bf v}^{T}_2S^{-T}S^{T}\begin{bmatrix}
				D_2&0\\
				0&D_2
			\end{bmatrix}SS^{-1}PAS\\
			&\iff {\bf v}^{T}_1\begin{bmatrix}
				D_1&0\\
				0&D_1
			\end{bmatrix}={\bf v}^{T}_2\begin{bmatrix}
				D_2&0\\
				0&D_2
			\end{bmatrix}PA\\
			&\iff {\bf v}^{T}_1\begin{bmatrix}
				D_1&0\\
				0&D_1
			\end{bmatrix}={\bf v}^{T}_1A^{T}P^{T}\begin{bmatrix}
				D_2&0\\
				0&D_2
			\end{bmatrix}PA\\
			&\iff {\bf v}^{T}_1\begin{bmatrix}
				D_1&0\\
				0&D_1
			\end{bmatrix}={\bf v}^{T}_1A^{T}\begin{bmatrix}
				D_1&0\\
				0&D_1
			\end{bmatrix}A.
		\end{aligned}
	\end{equation*}
	
	Based on the compatibility result above, we know that $\varphi_S\circ \theta_{CH_n}$ induces a unique map $\Phi_S:\Theta_{CH_{n}}/\sim_{\star}\rightarrow{PH}_n/\sim_{sys}$ defined as $\Phi_S([{\bf d},{\bf v}]_{\star})=\left[S^T\begin{bmatrix}
		D&0\\
		0&D
	\end{bmatrix}S,S^{-1}{\bf v}\right]_{sys}$. We now verify that $\Phi_S$ does not depend on the choice of $S\in Sp(2n,\mathbb{R})$.\\\par
	
\noindent {\bf   $\Phi_S$ is independent of $S$.} It suffices to check that, for $S_1\neq S_2$, we have $\bigg(S_1^{T}\begin{bmatrix}
		D&0\\
		0&D
	\end{bmatrix}S_1,S_1^{-1}{\bf v}\bigg):=(Q^{\prime}_1,B^{\prime}_1)\sim_{sys}(Q^{\prime}_2,B^{\prime}_2):=\bigg(S_2^{T}\begin{bmatrix}
		D&0\\
		0&D
	\end{bmatrix}S_2,S_2^{-1}{\bf v}\bigg)$, which again goes back to checking (\ref{automorphism_condition}) holds for some invertible $L$. We claim that $L=S_2^{-1}S_1$ does the job. The first condition is 
\begin{equation*}
		\begin{aligned}
			&L\mathbb{J}Q^{\prime}_1=\mathbb{J}Q^{\prime}_2L\\
			&\iff L\mathbb{J}Q^{\prime}_1L^{-1}=\mathbb{J}Q^{\prime}_2\\
			&\iff
			S_2^{-1}S_1\mathbb{J}S_1^{T}\begin{bmatrix}
				D&0\\
				0&D
			\end{bmatrix}S_1S_1^{-1}S_2=\mathbb{J}S_2^{T}\begin{bmatrix}
				D&0\\
				0&D
			\end{bmatrix}S_2\\
			&\iff \begin{bmatrix}
				D&0\\
				0&D
			\end{bmatrix}=\begin{bmatrix}
				D&0\\
				0&D
			\end{bmatrix}.
		\end{aligned}
	\end{equation*}
The second condition is
	\begin{equation*}
		LB^{\prime}_1=B^{\prime}_2\iff S_2^{-1}S_1S^{-1}_1{\bf v}=S_2^{-1}{\bf v}\\
		\iff {\bf v}={\bf v}.
	\end{equation*}
The third condition is 
	\begin{equation*}
		\begin{aligned}j
			&B^{{\prime}{T}}_1Q^{\prime}_1=B^{{\prime}T}_2Q^{\prime}_2L\\
			&\iff {\bf v}^{T}S_1^{-T}S_1^{T}\begin{bmatrix}
				D&0\\
				0&D
			\end{bmatrix}S_1\\&~~~~~~~~~~={\bf v}^{T}S_2^{-T}S_2^{T}\begin{bmatrix}
				D&0\\
				0&D
			\end{bmatrix}S_2S_2^{-1}S_1\\
			&\iff {\bf v}^{T}\begin{bmatrix}
				D&0\\
				0&D
			\end{bmatrix}={\bf v}^{T}\begin{bmatrix}
				D&0\\
				0&D
			\end{bmatrix}.
		\end{aligned}
	\end{equation*}
Since $\Phi_S$ does not depend on $S\in Sp(2n,\mathbb{R})$, we may as well choose $S=\mathbb{J}_n$ and call it $\Phi$. Then $\Phi$ has the expression $\Phi([{\bf d},{\bf v}]_{\star})=\left[\begin{bmatrix}
		D&0\\
		0&D
	\end{bmatrix},{\bf v}\right]_{sys}$. We now verify that $\Phi$ is injective and surjective, and hence an isomorphism.\\\par
\noindent {\bf  $\Phi$ is surjective.} For an arbitrary choice $[Q,B]_{sys}$ of equivalence class, we take a representative $Q$ and $B$. Since $Q$ is symmetric positive-definite, by Williamson's theorem, $Q=S^{T}\begin{bmatrix}
		D&0\\
		0&D
	\end{bmatrix}S$ for some $S\in Sp(2n,\mathbb{R})$ and the diagonal entries of $D$ are nonnegative and can be identified with ${\bf d}$. Let ${\bf v}=S\cdot B$. Then we have $\Phi_S([{\bf d},{\bf v}]_{\star})=[Q,B]_{sys}$. Given that $\Phi_S=\Phi$ for any $S$, it holds that $\Phi([{\bf d},{\bf v}]_{\star})=[Q,B]_{sys}.$ This concludes $\Phi$ being surjective.\\\par
	
	\noindent {\bf  $\Phi$ is injective.} For $\left(\begin{bmatrix}
		D_1&0\\
		0&D_1
	\end{bmatrix},{\bf v}_1\right)\sim_{sys}\left(\begin{bmatrix}
		D_2&0\\
		0&D_2
	\end{bmatrix},{\bf v}_2\right)$, it means there exists some invertible $L$ such that the conditions in $(\ref{automorphism_condition})$ are all satisfied. We aim to show that $({\bf d}_1,{\bf v}_1)\sim_{\star} ({\bf d}_2,{\bf v}_2)$. The first condition gives \begin{equation*}
		\begin{aligned}
			&L\mathbb{J}\begin{bmatrix}
				D_1&0\\
				0&D_1
			\end{bmatrix}=\mathbb{J}\begin{bmatrix}
				D_2&0\\
				0&D_2
			\end{bmatrix}L\Rightarrow L\mathbb{J}\begin{bmatrix}
				D_1&0\\
				0&D_1
			\end{bmatrix}L^{-1}=\mathbb{J}\begin{bmatrix}
				D_2&0\\
				0&D_2
			\end{bmatrix}\\
			&\Rightarrow \det\left(\lambda \mathbb{I}-\mathbb{J}\begin{bmatrix}
				D_1&0\\
				0&D_1
			\end{bmatrix}\right)=\det\left(\lambda \mathbb{I}-\mathbb{J}\begin{bmatrix}
				D_2&0\\
				0&D_2
			\end{bmatrix}\right)\\
			&\Rightarrow \big(\lambda^2+d_{1,1}^2\big)\dots\big(\lambda^2+d_{1,n}^2\big)=\big(\lambda^2+d_{2,1}^2\big)\dots\big(\lambda^2+d_{2,n}^2\big).
		\end{aligned}
	\end{equation*} 
	Therefore, $\{d_{1,i}|i=1,\dots,n\}$ is the same as $\{d_{2,i}|i=1,\dots,n\}$ as a set, and this implies the existence of some $\sigma\in S_n$ such that $d_{2,i} = d_{1,\sigma(i)}$ for $i=1,\dots,n.$ In other words, there exists some permutation matrix $P_{\sigma}$ such that $P\begin{bmatrix}
		D_{1}&0\\0&D_{1}
	\end{bmatrix}P^{T}=\begin{bmatrix}
		D_{2}&0\\0&D_{2}
	\end{bmatrix}$. Thus, (i) of Definition \ref{equivalence_relation} holds. Further, we have 
	\begin{equation*}
		\begin{aligned}
			&L\mathbb{J}\begin{bmatrix}
				D_1&0\\
				0&D_1
			\end{bmatrix}=\mathbb{J}\begin{bmatrix}
				D_2&0\\
				0&D_2
			\end{bmatrix}L\\
			&\iff L\mathbb{J}\begin{bmatrix}
				D_1&0\\
				0&D_1
			\end{bmatrix}=\mathbb{J}P\begin{bmatrix}
				D_1&0\\
				0&D_1
			\end{bmatrix}P^TL\\
			&\iff L\mathbb{J}\begin{bmatrix}
				D_1&0\\0&D_1
			\end{bmatrix}=P\mathbb{J}\begin{bmatrix}
				D_1&0\\0&D_1
			\end{bmatrix}P^TL\\
			&\iff P^{T}L\mathbb{J}\begin{bmatrix}
				D_1&0\\0&D_1
			\end{bmatrix}=\mathbb{J}\begin{bmatrix}
				D_1&0\\0&D_1
			\end{bmatrix}P^{T}L\\
			&\iff A\mathbb{J}\begin{bmatrix}
				D_1&0\\0&D_1
			\end{bmatrix}=\mathbb{J}\begin{bmatrix}
				D_1&0\\0&D_1
			\end{bmatrix}A,
		\end{aligned}
	\end{equation*}
	if we denote $A:=P^{T}L$. Thus, (iii) of Definition \ref{equivalence_relation} holds true. The second condition of (\ref{automorphism_condition}) says ${\bf v}_2=L{\bf v}_1=PA{\bf v}_1$. Thus, (iv) of Definition \ref{equivalence_relation} holds true. Lastly, the third condition implies
	\begin{equation*}
		\begin{aligned}
			&{\bf v}^{T}_1\begin{bmatrix}
				D_1&0\\0&D_1
			\end{bmatrix}={\bf v}^{T}_2\begin{bmatrix}
				D_2&0\\0&D_2
			\end{bmatrix}L\\
			&\iff {\bf v}_1^{T}\begin{bmatrix}
				D_1&0\\
				0&D_1
			\end{bmatrix}=(PA{\bf v}_1)^T\left(P\begin{bmatrix}
				D_1&0\\
				0&D_1
			\end{bmatrix}P^T\right)(PA)\\
			&\iff {\bf v}^{T}_1\begin{bmatrix}
				D_1&0\\
				0&D_1
			\end{bmatrix}={\bf v}^{T}_1A^{T}\begin{bmatrix}
				D_1&0\\
				0&D_1
			\end{bmatrix}A,
		\end{aligned}
	\end{equation*}
	Thus, {\bf (ii)} in Definition \ref{equivalence_relation} holds. We conclude that $\Phi$ is injective.\\\par
	
\noindent {\bf   $\Phi$  is a homeomorphism with respect to the quotient topology.}  
Before we prove this statement, we first quote a lemma (see, for instance, \cite{mta}). 
    \begin{lemma}\label{pass_homeomorphism_to_quotient}
    	Let $X$ and $Y$ be sets equipped with equivalence relations $\sim_X$ and $\sim_Y$ respectively. If $\phi:X\rightarrow Y$ is a map such that, for any $x_1$, $x_2\in X$, $x_1\sim_X x_2$ if and only if $\phi(x_1)\sim_Y\phi(x_2)$, then $\phi$ projects to a unique map $\tilde{\phi}:X/\sim_X\rightarrow Y/\sim_Y$ between the quotient spaces given by $\tilde{\phi}([x]_{\sim_X})=[\phi(x)]_{\sim_Y}$ and such that the following diagram commutes. In particular, if $\phi$ is a homeomorphism between two topological spaces $X$ and $Y$, then $\tilde{\phi}$ is also a homeomorphism.
    	\[
    	\begin{tikzcd}
    		X \arrow{r}{\phi} \arrow{d}{\pi_{X}}& Y \arrow{d}{\pi_{Y}}  \\
    		X/\sim_{X} \arrow{r}{\tilde{\phi}} & Y/\sim_{Y} 
    	\end{tikzcd}
    	\]
    	
    \end{lemma}
  \noindent We now proceed with the proof.
\begin{description}
		\item[(i)] If $(Q_1,B_1)$ and $(Q_2,B_2)\in{PH}_n$ are linked by some linear symplectic map $S\in Sp(2n,\mathbb{R})$ by $(Q_2,B_2)=(S^{-T}Q_1S^{-1},SB_1)$, then $(Q_1,B_1)\sim_{sys}(Q_2,B_2)$. Therefore, as an immediate consequence of Williamson's normal form, we have ${PH}_n/\sim_{sys}={PH}^{diag}_n/\sim_{sys}$, where $PH^{diag}_n:=\left\{\left(\begin{bmatrix}D&0\\0&D\end{bmatrix},{\bf v}\right)\bigg|D={\rm diag}({\bf d}),d_i>0, {\bf v}\in\mathbb{R}^{2n}\right\}$.
		\item[(ii)] There is an obvious homeomorphism $\varphi:\Theta_{CH_{n}}\rightarrow{PH}_n^{diag}$ given by $\varphi({\bf d},{\bf v})=\left(\begin{bmatrix}D&0\\0&D\end{bmatrix},{\bf v}\right)$. Therefore, by identifying ${PH}_n/\sim_{sys}$ with ${PH}^{diag}_n/\sim_{sys}$, the induced map of $\varphi$ on the quotients is exactly $\Phi$. By Lemma \ref{pass_homeomorphism_to_quotient}, $\Phi$ is also a homeomorphism.
\end{description}
	
	To summarize, we have that the following diagram commutes.
	\[
	\begin{tikzcd}
		\Theta_{CH_{n}} \arrow{r}{\varphi_S\circ \theta_{CH _n}} \arrow{d}{\pi_{\star}}& {PH}_{n} \arrow{d}{\pi_{sys}}  \\
		\Theta_{CH_{n}}/\sim_{\star} \arrow{r}{\cong}[swap]{\Phi_S=\Phi} & {PH}_{n}/\sim_{sys} \arrow{l}{} 
	\end{tikzcd}
	\]

\subsection{Proof of Proposition~\ref{groupoid_orbit}}

	The axioms of being a groupoid mostly follow from the definition. Here, we only check the closure of the multiplication operation $m$, i.e. $(L_1L_2, (Q_2,B_2))\in \mathcal{G}_n$. Note that 
	\begin{equation*}
		\mathbb{J}^T(L_1L_2)\mathbb{J}Q_2(L_1L_2)^{-1}=\mathbb{J}^TL_1\mathbb{J}(\mathbb{J}^TL_2\mathbb{J}Q_2L_2^{-1})L^{-1}_1=\mathbb{J}^TL_1\mathbb{J}Q_1L^{-1}_1
	\end{equation*} is symmetric positive-definite. On the other hand, we have 
	\begin{align*}
		\mathbb{J}^T(L_1L_2)^T&\mathbb{J}(L_1L_2)B_2=\mathbb{J}^TL^T_2L^{T}_1\mathbb{J}L_1L_2B_2=\mathbb{J}^TL^T_2L^{T}_1\mathbb{J}L_1B_1\\
		&=\mathbb{J}^TL^T_2\mathbb{J}(\mathbb{J}^TL^{T}_1\mathbb{J}L_1B_1)=\mathbb{J}^TL^T_2\mathbb{J}B_1=\mathbb{J}^TL^T_2\mathbb{J}L_2B_2=B_2.
	\end{align*} Thus, closure of multiplication is proved. We also need to show that $\alpha$ and $\beta$ are submersions. Indeed, for $(L,(Q,B))\in \mathcal{G}_n$ and $(N,(P,C))\in T_{(L,(Q,B))}\mathcal{G}_n$, it holds that
	\begin{align*}
		T_{(L,(Q,B))}\alpha(N,(P,C))&=\frac{d}{dt}\bigg|_{t=0}\bigg(\mathbb{J}^T(L+tN)\mathbb{J}(Q+tP)(L+tN)^{-1},(L+tN)(B+tC)\bigg)\\
		&=(\mathbb{J}^TN\mathbb{J}QL^{-1}+\mathbb{J}^TL\mathbb{J}PL^{-1}-\mathbb{J}^TL\mathbb{J}QL^{-1}NL^{-1},LC+NB).
	\end{align*}
	Obviously, $LC+NB$ can traverse $\mathbb{R}^{2n}$ with varying $N\in GL(2n,\mathbb{R})$ and $C\in \mathbb{R}^{2n}$. For the first component, we can take $N=L$ such that it becomes $\mathbb{J}^TL\mathbb{J}PL^{-1}$. Since the tangent space of an open submanifold can be identified with the tangent space of the whole manifold, plus the fact that the tangent space of a vector space can be identified with itself, we naturally conclude that $T_{(L,(Q,B))}\alpha$ is surjective and hence $\alpha$ is a submersion. Similarly, one check that $\beta$ is a submersion. \\
	Then, the orbit of the groupoid containing $(Q,B)$ is given by 
	\begin{align*}
		\alpha(\beta^{-1}(Q,B))&=\alpha(\{(L,(Q,B))| L \text{ satisfies 1.(i) and 1.(ii) in Definition \ref{groupoid} }\})\\
		&=\left\{(\mathbb{J}^TL\mathbb{J}QL^{-1},LB)| L \text{ satisfies 1.(i) and 1.(ii) in Definition \ref{groupoid} }\right\}\\
		&=\left\{(Q^{\prime},B^{\prime})| (Q^{\prime},B^{\prime})\sim_{sys}(Q,B)\right\}
	\end{align*}

\subsection{Proof of Proposition~\ref{characterization2}}

	\noindent {\bf  $f$ is well-defined.} If $({\bf d}_1,{\bf v}_1)\sim_{sys}({\bf d}_2,{\bf v}_2)$, then there exists an invertible matrix $L_0$ such that \begin{equation*}
		\left\{
		\begin{aligned}
			&L_0\cdot g^{ctr}_1({\bf d}_1)=g^{ctr}_1({\bf d}_2)\cdot L_0\\
			&L_0\cdot \left(
			0 , 0 ,\cdots, 0, 1\right)^{T}=\left(
			0 , 0 ,\cdots, 0, 1\right)^{T}\\
			&g^{ctr}_2({\bf d}_1,{\bf v}_1)=g^{ctr}_2({\bf d}_2,{\bf v}_2)\cdot L_0
		\end{aligned} \right.
	\end{equation*}
	Since we are restricting on canonical systems, we apply the representation theorem to deduce the existence of some invertible matrices $L_i$, $i=1,2$ such that \begin{equation*}
		\left\{
		\begin{aligned}
			&L_i\cdot g^{ctr}_1({\bf d}_i)=\mathbb{J}Q_i\cdot L_i\\
			&L_i\cdot \left(
			0 , 0 ,\cdots, 0, 1\right)^{T}=B_i\\
			&g^{ctr}_2({\bf d}_i,{\bf v}_i)=B^{T}_iQ_i\cdot L_i
		\end{aligned} \right.
	\end{equation*}
	Now, check $L=L_2L_0L^{-1}_1$ is invertible and satisfies
	\begin{equation*}
		\left\{
		\begin{aligned}
			L\mathbb{J}Q_1&=\mathbb{J}Q_2L\\
			LB_1&=B_2\\
			B^{T}_1Q_1&=B^{T}_2Q_2L.
		\end{aligned} \right.
	\end{equation*}
	Therefore, $(Q_1,B_1)\sim_{sys}(Q_2,B_2)$.\\\par
	
\noindent {\bf   $f$ is surjective.} This is obvious, see the proof above.\\\par
	
\noindent {\bf   $f$ is injective.} Given all the matrices are invertible, this can be shown by essentially reversing the proof of $f$ being well-defined.

\subsection{Proof of Proposition~\ref{action_verification}}
	We directly verify that
	\begin{equation*}
		\begin{aligned}
			&\Gamma_{(\sigma,(\theta_1,\dots,\theta_n)^T)\circ(\bar{\sigma},(\bar{\theta}_1,\dots,\bar{\theta}_n)^T)}({\bf d},{\bf v})\\
			&=\Gamma_{(\sigma\bar{\sigma},(\theta_1,\dots,\theta_n)^T+P_{{\sigma}}\cdot(\bar{\theta}_1,\dots,\bar{\theta}_n)^T)}({\bf d},{\bf v})\\
			&=(P_{\sigma\bar{\sigma}} \cdot {\bf d}, \Gamma_{(\theta_1,\dots,\theta_n)^T}\circ\Gamma_{P_{\sigma}\cdot(\bar{\theta}_1,\dots,\bar{\theta}_n)^T}\bigg(\begin{bmatrix}
				P_{\sigma\bar{\sigma}}&0\\
				0&P_{\sigma\bar{\sigma}}
			\end{bmatrix}{\bf v}\bigg))\\
			&=(P_\sigma P_{\bar{\sigma}} \cdot {\bf d},\\ &~~~~~\Gamma_{(\theta_1,\dots,\theta_n)^T}\circ{\left[
				\begin{array}{ccc|ccc}
					\cos\bar{\theta}_{\sigma(1)}&&\bigzero&-\sin\bar{\theta}_{\sigma(1)}&&\bigzero\\
					&\ddots&&&\ddots&\\
					\bigzero&&\cos\bar{\theta}_{\sigma(n)}&	\bigzero&&-\sin\bar{\theta}_{\sigma(n)}\\
					\hline
					\sin\bar{\theta}_{\sigma(1)}&&\bigzero&\cos\bar{\theta}_{\sigma(1)}&&\bigzero\\
					&\ddots&&&\ddots&\\
					\bigzero&&\sin\bar{\theta}_{\sigma(n)}&	\bigzero&&\cos\bar{\theta}_{\sigma(n)}\\
				\end{array}
				\right]}\begin{bmatrix}
				P_\sigma &0\\
				0&P_\sigma
			\end{bmatrix}\begin{bmatrix}
				P_{\bar{\sigma}} &0\\
				0&P_{\bar{\sigma}}
			\end{bmatrix}{\bf v})\\
			&=(P_\sigma P_{\bar{\sigma}} \cdot {\bf d},\Gamma_{(\theta_1,\dots,\theta_n)^T}\begin{bmatrix}
				P_\sigma &0\\
				0&P_\sigma
			\end{bmatrix}\Gamma_{(\bar{\theta}_1,\dots,\bar{\theta}_n)^T}\begin{bmatrix}
				P_{\bar{\sigma}} &0\\
				0&P_{\bar{\sigma}}
			\end{bmatrix}{\bf v})\\
			&=\Gamma_{(\sigma,(\theta_1,\dots,\theta_n)^T)}(\Gamma_{(\bar{\sigma},(\bar{\theta}_1,\dots,\bar{\theta}_n)^T)}({\bf d},{\bf v})).
		\end{aligned}
	\end{equation*}

\subsection{Proof of Proposition~\ref{characterization3}}

	Recall that $({\bf d}_1,{\bf v}_1)$ and $({\bf d}_2,{\bf v}_2)$ lie in the same $(S_n\rtimes_{\phi}\mathbb{T}^{n})$-orbit if and only if for some $\sigma\in S_n$ and $(\theta_1,\dots,\theta_n)\in {\mathbb{T}}^n$.
	\begin{description}
		\item[(i)] $d_{2,i} = d_{1,\sigma(i)},~i=1,\dots,n.$
		\item[(ii)] ${v}_{2,i}^2+{v}_{2,n+i}^2={v}_{1,\sigma(i)}^2+{v}_{1,n+\sigma(i)}^2$.
	\end{description}	
	Clearly, {\bf (i)} above is equivalent to Proposition \ref{equivalence_conditon_on_D,v} {\bf (i)}. Moreover, Proposition \ref{equivalence_conditon_on_D,v} {\bf (ii)} implies that for $k=0,\dots,n-1$,
	
	\begin{equation*}
		\begin{aligned}
			&{\bf v}_1^{T}\begin{bmatrix}
				F_{1,k}&0\\
				0&F_{1,k}
			\end{bmatrix} {\bf v}_1=(P^{T}{\bf v}_2)^{T}\begin{bmatrix}
				F_{1,k}&0\\
				0&F_{1,k}
			\end{bmatrix} P^{T}{\bf v}_2\\
			&\iff\sum_{i=1}^{n}F_{1,k}^{(i)}\cdot(v_{1,i}^2+v_{1,n+i}^2)\\
			&~~~~~~~~~~~~~~~~~~~~~=\sum_{i=1}^{n}F_{1,k}^{(i)}\cdot (v_{2,\sigma^{-1}(i)}^2+v_{2,n+\sigma^{-1}(i)}^2)
		\end{aligned}
	\end{equation*}
	
	Now, let $\bar{R}_1=(R_{1,1},\dots,R_{1,n})^{T}$, where $R_{1,i}=v_{1,i}^2+v_{1,n+i}^2$. Let $\bar{R}_2=(R_{2,1},\dots,R_{2,n})^{T}$, where $R_{2,i}=v_{2,\sigma^{-1}(i)}^2+v_{2,n+\sigma^{-1}(i)}^2$. Identify the diagonal matrix $F_{1,k}$ as a row vector in $\mathbb{R}^{n}$. Then, the above is equivalent to saying that the inner product of $F_{1,k}$ with $\bar{R}_1$ and $\bar{R}_2$ are the same for all $k=0,\dots,n-1$. Rewrite these inner products as matrix multiplication gives
	
	\begin{equation*}
		\begin{bmatrix}F_{1,0} \\ F_{1,1} \\ \rowsvdots \\ F_{1,n-1}\end{bmatrix}\cdot(\bar{R}_1-\bar{R}_2)=0. 
	\end{equation*}
	The determinant of this matrix is $\big(\prod_{i=1}^{n}d_i\big)\cdot\big(\prod_{1\leq j<k\leq n}^{}d_j^2-d_k^2\big)$. Since there are no repeated symplectic eigenvalues, we must have $\bar{R}_1=\bar{R}_2$, namely ${v}_{1,i}^2+{v}_{1,n+i}^2={v}_{2,\sigma^{-1}(i)}^2+{v}_{2,n+\sigma^{-1}(i)}^2$ for all $i=1,\dots,n$. Thus, (ii) holds by inversing the permutation $\sigma$. The converse is clearly true. 

\subsection{Proof of Proposition~\ref{characterization4}}

\noindent {\bf  $f$ is well-defined.} Let $({\bf d}_1,{\bf v}_1)$ and $({\bf d}_2,{\bf v}_2)$ be in the same orbit of the $(S_n\rtimes_{\phi}\mathbb{T}^{n})$-action. This means there exists $\sigma\in S_n$ and ${\bf \Theta}\in \mathbb{T}^{n}$ such that $\Gamma_{\sigma}({\bf d}_1)={\bf d}_2$ and $\Gamma_{{\bf \Theta}}(\Gamma_{\sigma}({\bf v}_1))={\bf v}_2$. This immediately implies $({\bf d}_1)_{\uparrow}=({\bf d}_2)_{\uparrow}$, as well as $\mathcal{R}({\bf v}_2)=\mathcal{R}(\Gamma_{\sigma}({\bf v}_1))$. Moreover, let $\sigma_i\in S_n$ be the unique permutation such that $\Gamma_{\sigma_i}({\bf d}_i)=({\bf d}_i)_{\uparrow}$, $i=1,2$. Then we have, \begin{equation*}
		\begin{aligned}
			{\bf d}_2&=\Gamma_{\sigma^{-1}_2}(({\bf d}_2)_{\uparrow})=\Gamma_{\sigma^{-1}_2}(({\bf d}_1)_{\uparrow})\\
			&=(\Gamma_{\sigma^{-1}_2}\circ\Gamma_{\sigma_1})({\bf d}_1)\\
			&=\Gamma_{\sigma^{-1}_2 \sigma_1}({\bf d}_1).
		\end{aligned}
	\end{equation*}
	Since all the entries of ${\bf d}$ are distinct, we necessarily have $\sigma=\sigma^{-1}_2\sigma_1$. We want to show $\mathcal{R}(\Gamma_{\sigma_1}({\bf v}_1))=\mathcal{R}(\Gamma_{\sigma_2}({\bf v}_2))$, but since $\mathcal{R}$ and $\Gamma_{\sigma}$ commutes for any $\sigma\in S_n$, this is equivalent to
	\begin{equation*}
		\begin{aligned}
			&~~~~~~\Gamma_{\sigma_1}(\mathcal{R}({\bf v}_1))=\Gamma_{\sigma_2}(\mathcal{R}({\bf v}_2))\\
			&\iff \Gamma_{\sigma_1}(\mathcal{R}({\bf v}_1))=\Gamma_{\sigma_2}(\mathcal{R}(\Gamma_{\sigma}({\bf v}_1)))\\
			&\iff \Gamma_{\sigma_1}(\mathcal{R}({\bf v}_1))=\Gamma_{\sigma_2}(\mathcal{R}(\Gamma_{\sigma^{-1}_2\sigma_1}({\bf v}_1)))\\
			&\iff \Gamma_{\sigma_1}(\mathcal{R}({\bf v}_1))=\mathcal{R}(\Gamma_{\sigma_1}({\bf v}_1)),
		\end{aligned}
	\end{equation*}
	which is clearly true.\par
\noindent {\bf  $f$ is surjective.} This is obvious.\par
\noindent {\bf  $f$ is injective.} Now suppose $(({\bf d}_1)_{\uparrow},\mathcal{R}(\Gamma_{\sigma_1}({\bf v}_1)))=(({\bf d}_2)_{\uparrow},\mathcal{R}(\Gamma_{\sigma_2}({\bf v}_2)))$. This immediately implies the existence of some $\sigma\in S_n$ such that $\Gamma_{\sigma}({\bf d}_1)={\bf d}_2$. On the other hand, since ${\bf d}_i=\Gamma_{\sigma^{-1}_i}({\bf d}_i)_{\uparrow}, i=1,2$, we have $\sigma=\sigma^{-1}_2\sigma_1$ and hence ${\bf d}_2=\Gamma_{\sigma^{-1}_2\sigma_1}({\bf d}_1)$. On the other hand, $\mathcal{R}(\Gamma_{\sigma_1}({\bf v}_1))=\mathcal{R}(\Gamma_{\sigma_2}({\bf v}_2))$ implies $\mathcal{R}(\Gamma_{\sigma^{-1}_2\sigma_1}({\bf v}_1))=\mathcal{R}({\bf v}_2)$, which further implies the existence of some ${\bf \Theta}\in \mathbb{T}^{n}$ such that ${\bf v}_2=\Gamma_{\bf \Theta}(\Gamma_{\sigma^{-1}_2\sigma_1}({\bf v}_1))$. This concludes that $({\bf d}_1,{\bf v}_1)$ and $({\bf d}_2,{\bf v}_2)$ lie in the same orbit.

\subsection{Proof of Theorem~\ref{higher dimension}}
\label{Proof of Theorem higher dimension}

\noindent {\bf   Proof of part (i).}
	Say we are given a latent system \begin{equation}\label{small_dimension}
		\left\{
		\begin{aligned}
			\dot{{\bf z}}&=\mathbb{J}_{n}Q{\bf z}+Bu\\
			y&={B}^{T}Q{\bf z},
		\end{aligned} \right.
	\end{equation} where $\mathbb{J}_{n}=\begin{bmatrix}
		0 & \mathbb{I}_{n}\\
		-\mathbb{I}_{n} & 0
	\end{bmatrix}$, $B\in\mathbb{R}^{2n}$ and $Q$ a $2n$ by $2n$ symmetric, positive-definite matrix. Consider the matrix \begin{equation*}\begin{bmatrix}
			\mathbb{J}_{n}&0\\
			0&\mathbb{J}_{m-n}
		\end{bmatrix}=
		\begin{bmatrix}
			0 & \mathbb{I}_{n}&0&0\\
			-\mathbb{I}_{n} & 0&0&0\\
			0&0&0&\mathbb{I}_{m-n}\\
			0&0&-\mathbb{I}_{m-n}&0
		\end{bmatrix}.\end{equation*} There exists a conjugate transform by an orthogonal matrix that turns this matrix into $\mathbb{J}_m$, since only elementary row(column) permutation matrices are involved, and these elementary matrices are themselves orthogonal. That is, there exists $OO^{T}=O^{T}O=\mathbb{I}_{2m}$ such that $O\begin{bmatrix}
		\mathbb{J}_{n}&0\\
		0&\mathbb{J}_{m-n}
	\end{bmatrix}O^{T}=\mathbb{J}_m$. Now, consider the following linear port-Hamiltonian system in normal form\\ \begin{equation*}\footnotesize\label{large_dimension}
		\left\{
		\begin{aligned}
			\dot{\tilde{{\bf z}}}&=\left(O\begin{bmatrix}
				\mathbb{J}_{n}&0\\
				0&\mathbb{J}_{m-n}
			\end{bmatrix}O^{T}\right)\left(O\begin{bmatrix}
				Q&0\\
				0&\mathbb{I}_{2m-2n}
			\end{bmatrix}O^{T}\right)\tilde{{\bf z}}+O\begin{bmatrix}
				B\\0
			\end{bmatrix}u\\
			&=\mathbb{J}_m\left(O\begin{bmatrix}
				Q&0\\
				0&\mathbb{I}_{2m-2n}
			\end{bmatrix}O^{T}\right)\tilde{{\bf z}}+O\begin{bmatrix}
				B\\0
			\end{bmatrix}u\\
			y&=\left(O\begin{bmatrix}
				B\\0
			\end{bmatrix}\right)^{T}\left(O\begin{bmatrix}
				Q&0\\
				0&\mathbb{I}_{2m-2n}
			\end{bmatrix}O^{T}\right)\tilde{{\bf z}}\\
			&=\begin{bmatrix}
				{B}^{T}Q&0
			\end{bmatrix}O^{T}\tilde{{\bf z}}
		\end{aligned} \right.
	\end{equation*} with the change of variable ${{\bf z}}=O^{T}\tilde{{\bf z}}$, which is equivalent to \begin{equation*}\label{large_dimension}
		\left\{
		\begin{aligned}
			\dot{{\bf z}}&=\begin{bmatrix}
				\mathbb{J}_{n}Q&0\\
				0&\mathbb{J}_{m-n}
			\end{bmatrix}{{\bf z}}+\begin{bmatrix}B\\0
			\end{bmatrix}u\\
			y&=\begin{bmatrix}
				{B}^{T}Q&0
			\end{bmatrix}{{\bf z}},
		\end{aligned} \right.
	\end{equation*} which, restricted to the upper subspace, coinsides with (\ref{small_dimension}). Moreover, the matrix $O\begin{bmatrix}
		Q&0\\
		0&\mathbb{I}_{2m-2n}
	\end{bmatrix}O^{T}$ is again symmetric positive-definite by construction.

\medskip
\noindent {\bf   Proof of part (ii).}
	According to the system morphism conditions, we just need to check 
	\begin{equation*}
		\left\{
		\begin{aligned}
			L\mathbb{J}_nQ&=\mathbb{J}_mO\begin{bmatrix}
				Q&0\\
				0&\mathbb{I}_{2m-2n}
			\end{bmatrix}O^{T}L\\
			LB&=O\begin{bmatrix}
				B\\0
			\end{bmatrix}\\
			B^{T}Q&=\begin{bmatrix}
				B^{T}Q&0
			\end{bmatrix}O^{T}L.
		\end{aligned} \right.
	\end{equation*}
The first condition is 
	\begin{equation*}
		\begin{aligned}
			&L\mathbb{J}_nQ=\mathbb{J}_mO\begin{bmatrix}
				Q&0\\
				0&\mathbb{I}_{2m-2n}
			\end{bmatrix}\begin{bmatrix}
				\mathbb{I}_{2n}\\0
			\end{bmatrix}\\
			&\iff O^{T}L\mathbb{J}_nQ=O^{T}\mathbb{J}_mO\begin{bmatrix}
				Q&0\\
				0&\mathbb{I}_{2m-2n}
			\end{bmatrix}\begin{bmatrix}
				\mathbb{I}_{2n}\\0
			\end{bmatrix}\\
			&\iff \begin{bmatrix}
				I_{2n}\\0
			\end{bmatrix}\mathbb{J}_nQ=\begin{bmatrix}
				\mathbb{J}_n&0\\
				0&\mathbb{J}_{m-n}
			\end{bmatrix}\begin{bmatrix}
				Q&0\\
				0&\mathbb{I}_{2m-2n}
			\end{bmatrix}\begin{bmatrix}
				\mathbb{I}_{2n}\\0
			\end{bmatrix}\\
			&\iff \begin{bmatrix}
				\mathbb{I}_{2n}\\0
			\end{bmatrix}\mathbb{J}_nQ=\begin{bmatrix}
				\mathbb{J}_nQ\\0
			\end{bmatrix}.
		\end{aligned}
	\end{equation*} 
The second and third conditions are clear with $L=O\begin{bmatrix}
		\mathbb{I}_{2n}\\0
	\end{bmatrix}$

\subsection{Proof of Proposition~\ref{low_high_isomorphism}}
\noindent {\bf  $f$ is well-defined.} Given $(Q_1,B_1)\sim_{sys}(Q_2,B_2)$, there exists an invertible $L\in\mathbb{R}^{2n}$ such that (\ref{automorphism_condition}) is satisfied. Let $L^{\prime}=O\begin{bmatrix}L&0\\0&\mathbb{I}_{2m-2n}\end{bmatrix}O^{T}$. Check that $L^{\prime}$ satisfies the conditions (\ref{automorphism_condition}) together with $(Q^{\prime}_1,B^{\prime}_1)$ and $(Q^{\prime}_2,B^{\prime}_2)$. Therefore, $(Q^{\prime}_1,B^{\prime}_1)\sim_{sys} (Q^{\prime}_2,B^{\prime}_2)$.\\\par
\noindent {\bf  $f$ is surjective.} This is clear from definition of $(Q^{\prime},B^{\prime})$.\\\par
\noindent {\bf  $f$ is injective.} Given $(Q^{\prime}_1,B^{\prime}_1)\sim_{sys}(Q^{\prime}_2,B^{\prime}_2)$, it means there exists an invertible $L^{\prime}\in \mathbb{R}^{2m}$ such that $L^{\prime}$ satisfies the conditions in (\ref{automorphism_condition}) together with $(Q^{\prime}_1,B^{\prime}_1)$ and $(Q^{\prime}_2,B^{\prime}_2)$. Write the matrix $O^{T}L^{\prime}O$ in the form $\begin{bmatrix}
		L_1&L_2\\L_{3}&L_4
	\end{bmatrix}$, where $L_1\in\mathbb{R}^{2n}$. Then check $L_1$ satisfies the conditions (\ref{automorphism_condition}) together with $(Q_1,B_1)$ and $(Q_2,B_2)$. Therefore, $(Q_1,B_1)\sim_{sys} (Q_2,B_2)$.

\subsection{Proof of Proposition~\ref{Dv_correspondence}}
	Clearly, $Q^{\prime}$ is also symmetric and positive-definite. Thus, again by Williamson's theorem, $Q^{\prime}=(S^{\prime})^{T}\begin{bmatrix}
		D^{\prime}&0\\
		0&D^{\prime}
	\end{bmatrix}S^{\prime}$. As before, we have
	\begin{equation*}
		\begin{aligned}
			&\big(\lambda^2+d_1^{\prime2}\big)\cdots\big(\lambda^2+d_m^{\prime2}\big)\\
			&=\det\bigg(\lambda \mathbb{I}_{2m}-\begin{bmatrix}
				0&D^{\prime}\\
				-D^{\prime}&0
			\end{bmatrix}\bigg)\\
			&=\det\bigg(\lambda \mathbb{I}_{2m}-(S^{\prime})^{-1}\begin{bmatrix}
				0&D^{\prime}\\
				-D^{\prime}&0
			\end{bmatrix}S^{\prime}\bigg)\\
			&=\det\bigg(\lambda \mathbb{I}_{2m}-\mathbb{J}_{m}(S^{\prime})^{T}\begin{bmatrix}
				D^{\prime}&0\\
				0&D^{\prime}
			\end{bmatrix}S^{\prime}\bigg)\\
			&=\det(\lambda \mathbb{I}_{2m}-\mathbb{J}_{m}Q^{\prime})\\
			&=\det\bigg(\lambda \mathbb{I}_{2m}-\mathbb{J}_{m}O\begin{bmatrix}
				Q&0\\
				0&\mathbb{I}_{2m-2n}
			\end{bmatrix}O^{T}\bigg)\\
			&=\det\bigg(\lambda \mathbb{J}_{m}+\mathbb{J}_{m}O\begin{bmatrix}
				Q&0\\
				0&\mathbb{I}_{2m-2n}
			\end{bmatrix}O^{T}\mathbb{J}_{m}^{T}\bigg)\\
			&=\det\bigg(\lambda \mathbb{J}_{m}+\mathbb{J}_{m}O\begin{bmatrix}
				Q&0\\
				0&\mathbb{I}_{2m-2n}
			\end{bmatrix}(\mathbb{J}_{m}O)^{-1}\bigg)\\
			&=\det\bigg(\lambda (\mathbb{J}_{m}O)^{-1}\mathbb{J}_{m}(\mathbb{J}_{m}O)+\begin{bmatrix}
				Q&0\\
				0&\mathbb{I}_{2m-2n}
			\end{bmatrix}\bigg)\\
			&=\det\bigg(\lambda O^{T}\mathbb{J}_{m}O+\begin{bmatrix}
				Q&0\\
				0&\mathbb{I}_{2m-2n}
			\end{bmatrix}\bigg)\\
			&=\det\bigg(\begin{bmatrix}
				\lambda \mathbb{J}_{n}&0\\
				0&\lambda \mathbb{J}_{m-n}
			\end{bmatrix}+\begin{bmatrix}
				Q&0\\
				0&\mathbb{I}_{2m-2n}
			\end{bmatrix}\bigg)\\
			&=\det(\lambda \mathbb{J}_{m-n}+\mathbb{I}_{2m-2n})\cdot \det(\lambda \mathbb{J}_{n}+Q)\\
			&=(\lambda^2+1)^{m-n}\cdot \det(\lambda \mathbb{I}_{2n}-\mathbb{J}_{n}Q)\\
			&=(\lambda^2+1)^{m-n}\big(\lambda^2+d_1^2\big)\cdots\big(\lambda^2+d_n^2\big)
		\end{aligned}
	\end{equation*}

	If we fixed the order of symplectic eigenvalues ${\bf d}^{\prime}$ according to ${\bf d}^{\prime}=(d_1,\dots,d_n,1,\dots,1)$, then 
	\begin{equation*}
		\begin{aligned}
			&Q^{\prime}=O\begin{bmatrix}
				Q&0\\
				0&\mathbb{I}_{2m-2n}
			\end{bmatrix}O^{T}\\
			&=O\begin{bmatrix}
				S^{T}\begin{bmatrix}
					D&0\\
					0&D
				\end{bmatrix}S&0\\
				0&\mathbb{I}_{2m-2n}
			\end{bmatrix}O^{T}\\
			&=O\begin{bmatrix}
				S^{T}&0\\
				0&\mathbb{J}_{m-n}
			\end{bmatrix}\begin{bmatrix}
				\begin{matrix}
					D&0\\
					0&D
				\end{matrix}&\rvline&\bigzero\\ \hline \bigzero& \rvline &
				\begin{matrix}\mathbb{I}_{m-n}&0\\
					0&\mathbb{I}_{m-n}
			\end{matrix}\end{bmatrix}\begin{bmatrix}
				S&0\\
				0&\mathbb{J}_{m-n}
			\end{bmatrix}O^{T}\\
			&=O\begin{bmatrix}
				S^{T}&0\\
				0&\mathbb{J}^T_{m-n}
			\end{bmatrix}O^{T}\begin{bmatrix}\begin{matrix}
					D&0\\
					0&\mathbb{I}_{m-n}\end{matrix} &\rvline &\bigzero\\ \hline \bigzero &\rvline & \begin{matrix}
					D&0\\
					0&\mathbb{I}_{m-n}\end{matrix}
			\end{bmatrix}O\begin{bmatrix}
				S&0\\
				0&\mathbb{J}_{m-n}
			\end{bmatrix}O^{T}\\
		\end{aligned}
	\end{equation*}
	Now, we check the matrix $O\begin{bmatrix}
		S&0\\
		0&\mathbb{J}_{m-n}
	\end{bmatrix}O^{T}$ is symplectic, i.e.
	\begin{equation*}
		\begin{aligned}
			&\bigg(O\begin{bmatrix}
				S&0\\
				0&\mathbb{J}_{m-n}
			\end{bmatrix}O^{T}\bigg)^{T}\mathbb{J}_{m}\bigg(O\begin{bmatrix}
				S&0\\
				0&\mathbb{J}_{m-n}
			\end{bmatrix}O^{T}\bigg)\\
			&=O\begin{bmatrix}
				S^{T}&0\\
				0&\mathbb{J}^T_{m-n}
			\end{bmatrix}\big(O^{T}\mathbb{J}_{m}O\big)\begin{bmatrix}
				S&0\\
				0&\mathbb{J}_{m-n}
			\end{bmatrix}
			O^{T}\\
			&=O\begin{bmatrix}
				S^{T}&0\\
				0&\mathbb{J}^T_{m-n}
			\end{bmatrix}\begin{bmatrix}
				\mathbb{J}_{n}&0\\
				0&\mathbb{J}_{m-n}
			\end{bmatrix}\begin{bmatrix}
				S&0\\
				0&\mathbb{J}_{m-n}
			\end{bmatrix}
			O^{T}\\
			&=O\begin{bmatrix}
				\mathbb{J}_{n}&0\\
				0&\mathbb{J}_{m-n}
			\end{bmatrix}O^{T}=\mathbb{J}_{m}.
		\end{aligned}
	\end{equation*}
	Therefore, $O\begin{bmatrix}
		S&0\\
		0&\mathbb{J}_{m-n}
	\end{bmatrix}O^{T}$ is a symplectic matrix diagonalizing $Q^{\prime}$ in Williamson's theorem. Then, we deduce
	\begin{equation*}
		{\bf v}^{\prime}=S^{\prime}B^{\prime}=O\begin{bmatrix}
			S&0\\
			0&\mathbb{J}_{m-n}
		\end{bmatrix}O^{T}O\begin{bmatrix}
			B\\0
		\end{bmatrix}
		=O\begin{bmatrix}
			SB\\0
		\end{bmatrix}=O\begin{bmatrix}
			{\bf v}\\0
		\end{bmatrix}.
	\end{equation*}

\subsection{Proof of Proposition~\ref{Dv_QB_correspondence}}

	Similar to the proof of Theorem \ref{characterization}, simply replace $Q$ with $O\begin{bmatrix}
		Q&0\\
		0&\mathbb{I}_{2m-2n}
	\end{bmatrix}O^{T}$, $B$ with $O\begin{bmatrix}
		B\\
		0
	\end{bmatrix}$, $S$ with $O\begin{bmatrix}
		S&0\\
		0&\mathbb{J}_{m-n}
	\end{bmatrix}O^{T}$, $D$ with $\begin{bmatrix}
		D&0\\0&\mathbb{I}_{m-n}
	\end{bmatrix}$, ${\bf v}=\begin{bmatrix}
		{\bf v}_{upper}\\{\bf v}_{lower}
	\end{bmatrix}$ with $\bar{{\bf v}}=\begin{bmatrix}
		{\bf v}_{upper}^T&0_{m-n}&{\bf v}_{lower}^T&0_{m-n}
	\end{bmatrix}^{T}$.

\subsection{A note on the design of discrete integrators on the transformed space}\label{integrator_section}
Even though in the numerical illustration we used just a simple Euler integration scheme,  structure-preserving integrations algorithms could have been used. In particular we could have used an implicit midpoint rule which is symplectic \cite{marsden_west_2001}, that is, it preserves the symplectic form $d{\bf q}\wedge d{\bf p}$. Recall that if $L_{Lag}({\bf q}, \dot{\bf q})$ is the Lagrangian function of the system of interest, then the midpoint integrator is obtained by using the discrete Lagrangian  
\begin{equation*}
	L^\alpha_d({\bf q}_0,{\bf q}_1,h)=hL_{Lag}((1-\alpha){\bf q}_0+\alpha {\bf q}_1,\frac{{\bf q}_1-{\bf q}_0}{h}),
\end{equation*}with $\alpha=\frac{1}{2}$ to approximate the exact discrete Lagrangian \begin{equation*}
	L^E_d({\bf q}_0,{\bf q}_1,h)=\int^{h}_0L_{Lag}({\bf q}_{0,1}(t),\dot{{\bf q}}_{0,1}(t))dt.
\end{equation*} 
Explicitly, the midpoint integrator for a linear autonomous Hamiltonian system is 
\begin{equation*}
	{\bf z}_{n+1}-{\bf z}_{n}=h\cdot \mathbb{J}Q\left(\frac{{\bf z}_{n+1}+{\bf z}_{n}}{2}\right), 
\end{equation*} 
which in terms of the controllable Hamiltonian representation reads \begin{equation}\label{integrator}
	L({\bf s}_{n+1}-{\bf s}_n)=\frac{h}{2}\mathbb{J}QL({\bf s}_{n+1}+{\bf s}_{n})=\frac{h}{2}L\cdot g^{ctr}_1({\bf d})({\bf s}_{n+1}+{\bf s}_n),
\end{equation}
where the second equality holds by the construction of $L$ in the proof of Theorem \ref{main_theorem} part {\bf (i)}.

Thus, for the symplectic structure to be preserved in the original space, we can merely integrate by requiring ${\bf s}_{n+1}-{\bf s}_n=\frac{h}{2}g^{ctr}_1({\bf d})({\bf s}_{n+1}+{\bf s}_n)$, where $g^{ctr}_1({\bf d})$ as we have seen, takes the form \begin{equation*}
	\begin{bmatrix}
		0&1&0&\dots&0 \\
		0&0&1&\dots&0 \\
		\vdots&\vdots&\ddots&&\vdots\\
		0&0&0&\dots&1 \\
		-a_0&-a_1&-a_2&\dots&-a_{2m-1}
	\end{bmatrix}_{2m\times 2m}.
\end{equation*} 
Therefore, the integrator is given by
\begin{equation*} {\bf s}_{n+1}=(\mathbb{I}_{2n}-\frac{h}{2}g^{ctr}_1({\bf d}))^{-1} (\mathbb{I}_{2n}+\frac{h}{2}g^{ctr}_1({\bf d}))\cdot {\bf s}_n,
\end{equation*}
where the matrix inverse is well-defined for sufficiently small time step $h$. Indeed, the integrator can be defined on the quotient space of $L$, since by (\ref{integrator}), we may as well choose $s_{n+1}$ such that \begin{equation*}
	{\bf s}_{n+1}-{\bf s}_n=\frac{h}{2}g^{ctr}_1({\bf d})({\bf s}_{n+1}+{\bf s}_n)+{\bf s}_{ker}
\end{equation*} for an arbitrary ${\bf s}_{ker}\in \ker(L)$.

By a similar argument, the midpoint rule in terms of observable Hamiltonian representation reads \begin{equation}\label{integrator2}
	{\bf s}_{n+1}-{\bf s}_n=L({\bf z}_{n+1}-{\bf z}_n)=\frac{h}{2}L\mathbb{J}Q({\bf z}_{n+1}+{\bf z}_{n})=\frac{h}{2} g^{obs}_1({\bf d})({\bf s}_{n+1}+{\bf s}_n),
\end{equation}
where the last equality holds by construction of $L$ from Theorem \ref{main_theorem} Part (ii).\par
Therefore, the integrator is \begin{equation*} {\bf s}_{n+1}=(\mathbb{I}_{2n}-\frac{h}{2}g^{obs}_1({\bf d}))^{-1} (\mathbb{I}_{2n}+\frac{h}{2}g^{obs}_1({\bf d}))\cdot {\bf s}_n.
\end{equation*}

In the case of port-Hamiltonian system, if the system is driven by some fiber-preserving external force $f_H$, i.e. some input as in our case, then the discrete Lagrange-d'Alembert Principle can be used to construct variational integrators so that all the correspondence relationships and error analysis of standard variational integrators still hold\cite{marsden_west_2001}. For example, the midpoint rule applied to the controllable Hamiltonian representation becomes 
\begin{align*}
	&{\bf z}_{n+1}-{\bf z}_{n}=h\cdot \mathbb{J}Q\left(\frac{{\bf z}_{n+1}+{\bf z}_{n}}{2}\right)+\begin{bmatrix}
		0\\
		f_H\left(\frac{{\bf z}_n+{\bf z}_{n+1}}{2}\right)
	\end{bmatrix}\\
	\Rightarrow &L({\bf s}_{n+1}-{\bf s}_n)=\frac{h}{2}L\cdot g^{ctr}_1({\bf d})({\bf s}_{n+1}+{\bf s}_n)+\begin{bmatrix}
		0\\
		f_H\left(\frac{L({\bf s}_{n+1}+{\bf s}_n)}{2}\right)
	\end{bmatrix}.
\end{align*}
Note that this structure-preserving integrator is not explicit in general.
\end{document}